\documentclass{amsart}
\usepackage[all]{xy}
\usepackage{latexsym}
\usepackage{amssymb}
\usepackage{amsmath}
\usepackage{amsthm}
\usepackage{amscd}
\usepackage{fancyhdr}
\usepackage[dvips]{graphicx}
\usepackage{fullpage}
\usepackage{mathrsfs}
\input cyracc.def
\xyoption{arc}

 \newfam\cyrfam

  \font\tencyr=wncyr10

  \font\sevencyr=wncyr7

  \font\fivecyr=wncyr5

  \textfont\cyrfam=\tencyr \scriptfont\cyrfam=\sevencyr

    \scriptscriptfont\cyrfam=\fivecyr

  \newfam\cyifam

  \font\tencyi=wncyi10

  \font\sevencyi=wncyi7

  \font\fivecyi=wncyi5

  \textfont\cyifam=\tencyi \scriptfont\cyifam=\sevencyi

    \scriptscriptfont\cyifam=\fivecyi

\def\id{{\mbox{1 \hskip -8pt 1}}}
\newcommand{\sgn}{{\mathit s  \mathit g\mathit  n}}
 \newcommand{\lon}{\longrightarrow}
 \newcommand{\bu}{\bullet}
 \newcommand{\ad}{{\mathrm a\mathrm d}}
 \newcommand{\rar}{\rightarrow}
 \newcommand{\hook}{\hookrightarrow}
 \newcommand{\Proof}{{\bf Proof}.\, }

\newcommand{\p}{{\partial}}
\newcommand{\Id}{{\mathrm I\mathrm d}}

 \newcommand{\Z}{{\mathbb Z}}
 \newcommand{\bS}{{\mathbb S}}

 \renewcommand{\P}{{\mathbb P}}
 \newcommand{\C}{{\mathbb C}}
 \newcommand{\R}{{\mathbb R}}
 \newcommand{\N}{{\mathbb N}}
 \newcommand{\K}{{\mathbb K}}
 \newcommand{\bbH}{{\mathbb H}}
\newcommand{\Conf}{{\mathit C\mathit o \mathit n\mathit f}}
 \newcommand{\ot}{\otimes}
 
  \newcommand{\Poly}{{\mathcal P}{\mathit o}{\mathit l}{\mathit y}}

 \newcommand{\Beq}{\begin{equation}}
 \newcommand{\Eeq}{\end{equation}}
 \newcommand{\Beqr}{\begin{eqnarray}}
 \newcommand{\Eeqr}{\end{eqnarray}}
 \newcommand{\Beqrn}{\begin{eqnarray*}}
 \newcommand{\Eeqrn}{\end{eqnarray*}}
 \newcommand{\Ba}{\begin{array}}
 \newcommand{\Ea}{\end{array}}
 \newcommand{\Bi}{\begin{itemize}}
 \newcommand{\Ei}{\end{itemize}}
 \newcommand{\Bc}{\begin{center}}
 \newcommand{\Ec}{\end{center}}
 \newcommand{\fg}{{\mathfrak g}}

 \newcommand{\fl}{{\mathfrak l}}

 \newcommand{\fG}{{\mathfrak G}}
 \newcommand{\fB}{{\mathfrak B}}

\newcommand{\fp}{{\mathfrak p}}
\newcommand{\ii}{{\mathfrak i}}

 \newcommand{\f}{{\mathcal O}}
 
 \newcommand{\cB}{{\mathcal B}}
 \newcommand{\cC}{{\mathcal C}}
 \newcommand{\caD}{{\mathcal D}}
 \newcommand{\cE}{{\mathcal E}}
 \newcommand{\cF}{{\mathcal F}}
 \newcommand{\cG}{{\mathcal G}}

 \newcommand{\caL}{{\mathcal L}}
 \newcommand{\cM}{{\mathcal M}}
 
 \newcommand{\cP}{{\mathcal P}}

 \newcommand{\cT}{{\mathcal T}}

 \newcommand{\cU}{{\mathcal U}}

 \newcommand{\al}{\alpha}
 
 \newcommand{\ga}{\gamma}
 
 \newcommand{\Ga}{\Gamma}

 \newcommand{\var}{\varepsilon}
 \newcommand{\la}{\lambda}
 \newcommand{\om}{\omega}

 \newcommand{\Img}{{\mathsf I\mathsf m}\, }

 \newcommand{\sip}{\smallskip}
 \newcommand{\bip}{\bigskip}
 \newcommand{\mip}{\vspace{2.5mm}}

\theoremstyle{plain}
\swapnumbers

\newtheorem{prop-def}[theorem]{Proposition-definition}

\newtheorem{f-theorem}{Formality Theorem}[section]
\newtheorem{main-theorem}{Main~Theorem}[section]
\newtheorem{section-theorem}{Theorem}[section]
\newtheorem{section-corollary}{Corollary}[section]

\theoremstyle{definition}

\newtheorem{fact-me}{Fact \cite{Me1}}[subsection]

\begin{document}

 \sloppy

 \newenvironment{proo}{\begin{trivlist} \item{\sc {Proof.}}}
  {\hfill $\square$ \end{trivlist}}

\long\def\symbolfootnote[#1]#2{\begingroup%
\def\thefootnote{\fnsymbol{footnote}}\footnote[#1]{#2}\endgroup}

 \title{Exotic  automorphisms of the \\Schouten algebra
 of polyvector fields}
 \author{ S.A.\ Merkulov}
\address{Sergei~A.~Merkulov: Department of Mathematics, Stockholm University, 10691 Stockholm, Sweden}
\email{sm@math.su.se}

 \begin{abstract}
Using  a new  compactification of the (braid) configuration space of $n$ points in
the upper half plane we give a construction of exotic $\caL ie_\infty$ automorphisms
of the Schouten algebra of polyvector fields on an affine space depending on the choice of a Kontsevich type propagator.

\sip

\end{abstract}
 \maketitle
\markboth{S.A.\ Merkulov}{Exotic automorphisms of the Schouten algebra
 of polyvector fields}

\begin{center}
\sc Contents
\end{center}
{\Small
1.  Introduction 
\vspace{1mm}

2. {\bf Configuration space $C_n$}
\vspace{-0.1mm}

2.1. A Fulton-MacPherson type compactification of $C_n$ 
\vspace{-0.1mm}

2.2.  The face complex of $\{\overline{C}_n\}$ as an operad of Leibniz$_\infty$ algebras 
\vspace{-0.1mm}

2.3.  A  semialgebraic structure on $\overline{C}_n$

\vspace{1mm}

3. {\bf Configuration space $C_{n,0}$ and its new compactifications}
\vspace{-0.1mm}

 3.1.  $C_{n,0}$ as a magnified ${C}_n$ 
\vspace{-0.1mm}

3.2.  New compactifications of $C_{n,0}$ 
\vspace{-0.1mm}

3.3. The face complex of $\widehat{C}_{n,0}$ as an operad $\cM or(\caL eib_\infty)$ 
\vspace{-0.1mm}

3.4. Semialgebraic atlas on $\widehat{C}_{n,0}$
\vspace{-0.1mm}

3.5. Angle functions on  $\widehat{C}_{2,0}$
\vspace{-0.1mm}

3.6.  Renormalized forgetful map

\vspace{1mm}

4. {\bf De Rham field theories on configuration spaces}
\vspace{-0.1mm}

4.1. Families of graphs
\vspace{-0.1mm}

4.2. De Rham field theory on $\overline{C}$ 
\vspace{-0.1mm}

4.3. De Rham field theory on $\overline{C}\sqcup \widehat{C}\sqcup \overline{C}$ 
\vspace{-0.1mm}

4.4. De Rham field theories from angular functions on $\widehat{C}_{2,0}$
\vspace{-0.1mm}

4.5. Gauge equivalence propagators and a proof of the homotopy equivalence theorem 
\vspace{-0.1mm}

4.6. Exotic transformations of Poisson structures 
\vspace{-0.1mm}

4.7. Example: symmetrized Kontsevich's propagators 
\vspace{-0.1mm}

4.8. Example: Kontsevich's (anti)propagator 
\vspace{-0.1mm}

4.9. De Rham field theory of Duflo's strange automorphism 
\vspace{-0.1mm}

5. {\bf Braid configuration spaces
}
\vspace{-0.1mm}

5.1. Compactified braid configuration spaces as operads $\caL ie_\infty$ and $\cM or(\caL ie_\infty)$
\vspace{-0.1mm}

5.2. De Rham field theories on braid configuration spaces 

\vspace{1mm}

6. {\bf Towards a new differential geometry
}
\vspace{-0.1mm}

\vspace{1mm}

{{\sc Appendix 1}: \bf Wheels and zeta function
}
\vspace{-0.1mm}

{{\sc Appendix 2}: \bf Example of a boundary strata
}
\vspace{-0.1mm}

{{\sc Appendix 3}: \bf  $\caL eib_\infty$ automorphisms of Maurer-Cartan sets
}
\vspace{-0.1mm}

{{\sc Appendix 4}:
\bf Weights of all 4-vertex graphs in a de Rham field theory on $\overline{C}_\bu$}

\vspace{-0.6mm}
{{\sc Appendix 5}: \bf Wheeled prop of polyvector fields
}
}

\bip

\section{Introduction}
{\bf 1.1. Statement of the result}. Let  $\cT_{poly}(\R^d)$ be the Lie algebra  of  polyvector fields on $\R^d$
equipped with the grading in which the Schouten brackets have degree $-1$.
This paper gives an explicit construction for a family of
 $\caL ie_\infty$-automorphisms,
$$
F^{\caL ie}=\{ F_n^{\caL ie}: \wedge^n \cT_{poly}(\R^d) \rar \cT_{poly}(\R^d)[2-2n]\}_{n\geq 1},
$$
of  $\cT_{poly}(\R^d)$ parameterized by PA differential 1-forms on the Kontsevich compactified configuration space $\overline{C}_{2,0}$ equipped with a certain semialgebraic structure. The
formulae are universal, i.e.\ independent of the dimension $d<\infty$, have the first
component $F_1^{\caL ie}$  equal to the identity map, and all the other components
are given by sums,
 $$
 F^{\caL ie}_n = \sum_{\Ga\in \fG_{n,2n-2}} \mbox{\sc C}_\Ga
 \Phi_\Ga, \ n\geq 2,
 $$
running over a family of graphs, $\fG_{n,2n-2}$, with $n$ vertices and $2n-2$ edges, where
\Bi
\item $\Phi_\Ga:\otimes^n \cT_{poly}(\R^d) \rar \cT_{poly}(\R^d)[2-2n]$ is a linear map
constructed from the graph $\Gamma$
via a certain simple procedure explained in \S 4.1,
\item  the numerical coefficient, $\mbox{\sc C}_\Ga$,  is given by an integral,
\Beq\label{Intro_c_Gamma}
\mbox{\sc C}_\Ga = \int_{\widehat{C}_{n,0}} \bigwedge_{e\in Edges(\Ga)} \frac{{\fp}_e^*(\om)}{2\pi},
\Eeq
  over a compactified configuration space, $\widehat{C}_{n,0}$, of certain equivalence classes of
   $n$  pairwise distinct points in the upper half plane
$$
\bbH :=\{x+i y\in \C\,|\, y\geq 0\}.
$$
\Ei
The big open cell, $C_{n,0}$, of $\widehat{C}_{n,0}$
is exactly the same as in \cite{Ko},
$$
C_{n,0}:={\{z_1, \ldots, z_n\in \bbH\, |\, z_i\neq z_j\ \mbox{for}\ i\neq j\}}/{G^{(1)}},
$$
$$
G^{(1)}:= \{z\rar az+b\ |\ a,b\in \R, a> 0\},
$$
but our compactification, $\widehat{C}_{n,0}$,  of $C_{n,0}$ is different from Kontsevich's
one, $\overline{C}_{n,0}$ for all $n$ except $n=1$ and $n=2$.
In the above formula the symbol $\fp_e$ stands for a  surjection ({\em not}\, equal to
the ordinary forgetful map, see \S 3.6)
 $\widehat{C}_{n,0}\rar \widehat{C}_{2,0}\simeq \overline{C}_{2,0}$ associated with an edge
$e$ of a graph $\Ga\in \fG_{n,2n-2}$  , and $\om$ stands for an arbitrary closed minimal differential
form on $\widehat{C}_{2,0}$ whose restriction to the boundary
$\p \widehat{C}_{2,0}\simeq S^1\sqcup S^1$ coincides with the standard {\em homogeneous}\, volume form
on each of the two boundary topological circles.  If we drop the requirement of homogeneity, then our formula describes
a $\caL ie_\infty$ quasi-isomorphism between certain $\caL ie_\infty$-extensions of the Schouten bracket
canonically associated with the values of $\om$ on the first and, respectively, the second boundary circle
of $\p \widehat{C}_{2,0}$; we show in \S 3 an explicit formula for such a
$\caL ie_\infty$-extension which looks as the one above except that it involves a family of graphs, $\fG_{n,2n-3}$,
with $n$ vertices and $2n-3$ edges
and a  different
compactified configuration space; the resulting family of (homotopy trivial) $\caL ie_\infty$-extensions of the Schouten
bracket is parameterized by semialgebraic functions on $S^1$ and includes, for example, the
one constructed by Shoikhet in \cite{Sh2}.

\mip

The family of new compactifications, $\{\widehat{C}_{n,0}\}_{n\geq 1}$, as well as its braid version,
 $\{\widehat{B}_{n,0}\}_{n\geq 1}$,
discussed in \S 5,  have nice operadic interpretations:
the fundamental chain complex of the first one is naturally the 2-coloured dg operad of $\caL eib_\infty$-morphisms
of $\caL eib_\infty$-algebras while the fundamental chain complex of the second has a natural structure of
 the 2-coloured dg operad of $\caL ie_\infty$
morphisms of $\caL ie_\infty$-algebras.  Here $\caL eib$ stands for the
operad of  Leibniz algebras introduced by J.-L.\ Loday in \cite{Lo}, and $\caL eib_\infty$ for its minimal resolution. Thus
the face structure underlying the compactification  $\widehat{C}_{n,0}$ suggests that there
might exist a generalization
of the above construction
producing more general $\caL eib_\infty$-automorphisms,
$$
F^{Leib}=\{ F_n^{Leib}: \otimes^n \cT_{poly}(\R^d) \rar \cT_{poly}(\R^d)[2-2n]\}_{n\geq 1},
$$
of the Schouten algebra. Any $\caL ie_\infty$-automorphism is, of course, a $\caL eib_\infty$-automorphism
but not vice versa. Though the symmetrization of a generic $\caL ieb_\infty$ automorphism  does {\em not}\,
give a $\caL ie_\infty$ automorphism in general,  they both
induce automorphisms,
$$
\al \rar F^{Leib}(\al):=\sum_{n\geq 1}\frac{\hbar^{n-1}}{n!}F^{Leib}_n(\al,\ldots,\al),\ \
\al \rar F^{Lie}(\al):=\sum_{n\geq 1}\frac{\hbar^{n-1}}{n!}F^{Lie}_n(\al,\ldots,\al),
$$
 of one and the same set
$$
\cM\cC(\cT_{poly}(\R^d)[[\hbar]]):=\left\{\al\in \cT_{poly}(\R^d)\ot\C[[\hbar]]:\, |\al|=2\ \mbox{and}\ [\al,\al]_{Schouten}=0\right\},
$$
of Poisson structures on  $\R^d$ depending on a formal parameter $\hbar$. The space $\R^d$ can, in general,
  be equipped with a non-trivial $\Z$-grading,
and $|\al|$ stands above for the total degree of a polyvector field $\al$ (so that $|\al|=2$ does not necessarily
imply that $\al$ is a {\em bi}-vector field).
Therefore, in the context of Poisson geometry, one can skip distinguishing the two notions,
 $\caL eib_\infty$ and
$\caL ie_\infty$, 
and talk simply about exotic automorphisms of finite-dimensional Poisson structures, or, even better,
about exotic automorphisms,
$$
F:(\Poly, d) \lon (\Poly, d),
$$
of a certain very simple dg free  wheeled prop\footnote{In the same vein  the Kontsevich formality map \cite{Ko}
 can be understood as a morphism
 of dg wheeled props $\caD efQ \lon \Poly$ (see \cite{Me-lec}).}, $\Poly$, controlling finite-dimensional
 Poisson geometry
 (see, e.g. \cite{ Me-lec, Me-Perm}
 for an elementary introduction into the language of wheeled operads and props in the context of Poisson
 geometry).

\sip

It is worth emphasizing that our formulae for exotic automorphisms of Poisson structures
 depend on the choice of
a Kontsevich type propagator, $\om(z_1,z_2)$.
Propagators introduced by Kontsevich in his theory of formality
maps \cite{Ko,Ko3} give  suitable propagators for our model (though the weights we produce from these propagators are in general different from Kontsevich ones). In fact, the original Kontsevich propagator gives via our formula a highly non-trivial quasi-isomorphism
from the original Schouten algebra to its $\caL ie_\infty$-extension constructed by Shoikhet
in \cite{Sh2}. The symmetrized singular $\frac{1}{2}$-propagator,
$$
\om_{\frac{1}{2}K}(z_1,z_2)=
\frac{1}{2i}\left(d\log
\frac{z_1-z_2}{\overline{z}_1 - z_2} + d\log
\frac{z_2-z_1}{\overline{z}_2 - z_1}
\right)
$$
introduced by Kontsevich in \cite{Ko2} gives an exotic (that is, homotopy non-trivial)
universal $\caL ie_\infty$ automorphism of the Schouten algebra of polyvector fields.
This propagator is also used  in \S 4.9 to give a de Rham field theory interpretation  of Kontsevich's generalization \cite{Ko2} of the famous  Duflo's strange automorphism which involves an infinite sequence of zeta values,
$
\left\{ \frac{\zeta(n)}{n(2\pi \sqrt{-1})^n}\right\}_{n\in \Z}.
$
The propagator $\om_{\frac{1}{2}K}$ is singular at the strata of collapsing points so that one has yet
to give a rigorous explanation of why it works.

\mip

{\bf 1.2. A motivation}. Let $\caD^\bu_{poly}(\R^d)$ be the Hochschild dg Lie
algebra of polydifferential operators on smooth (formal) functions on $\R^d$. Tamarkin proved
\cite{Ta1} existence of a family of $\caL ie_\infty$-quasi-isomorphisms,
$$
\left\{
F_a:  \caD_{poly}^\bu(\R^d) \lon \wedge^\bullet \cT_{poly}(\R^d)\right\}_{a\in \cM},
$$
parameterized by the set, $\cM$, of all possible
Drinfeld's Lie associators (see the original paper \cite{Dr} or the book \cite{EK} for a definition of $\cM$).
The Grothendieck-Teichmueller group, $GRT$, acts on $\cM$ \cite{Dr}
and hence on the above family, $\{F_a\}$, of formality maps. This in turn defines a map,
$$
\Ba{rccc}
\rho: & GRT & \lon & Aut(\cT_{poly}(\R^d))\\
& G &\lon &  F_{G(a)}\circ F^{-1}_a,
\Ea
$$
where $F^{-1}_a:  \cT_{poly}(\R^d) \rar  \caD_{poly}(\R^d)$ is a $\caL ie_\infty$-morphism which is homotopy
inverse to $F_a$ (it exists but, in general, is not uniquely defined).

\mip

{\sf Conjecture.} {\em There exists a non-trivial representation,
$GT \rar Aut(\cT_{poly}(\R^d))/\hspace{-1mm}\sim$, where $\sim$ stands for the homotopy equivalence relation.}

\mip

Another motivation --- which might lead to a new kind of $GT$ twisted  differential  geometry ---
is outlined in \S 6. That twisted geometry
 might be a useful gadget in the future study of quantum GT invariants.
\mip

{\bf 1.3. Some notation}. The set $\{1,2, \ldots, n\}$ is abbreviated to $[n]$;  its group of automorphisms is
denoted by $\bS_n$. The
cardinality of a finite set
$A$ is denoted by $\# A$. If $V=\oplus_{i\in \Z} V^i$ is a graded vector space, then
$V[k]$ stands for the graded vector space with $V[k]^i:=V^{i+k}$; for $v\in V^i$ we set $|v|:=i$. If $\om_1$ and $\om_2$
are differential forms on manifolds $X_1$ and, respectively, $X_2$, then the form $p_1^*(\om_1)\wedge p_2^*(\om_2)$
on $X_1\times X_2$,
where $p_1: X_1\times X_2\rar X_1$ and  $p_2: X_1\times X_2\rar X_2$ are natural projections,
is often abbreviated to $\om_1\wedge \om_2$.
\sip

The algebra, $\cT_{poly}(\R^d)$, of
 smooth polyvector fields on a finite-dimensional $\Z$-graded vector space $V\simeq \R^d$ is understood in this paper as a $\Z$-graded commutative algebra of smooth functions on the $\Z$-graded
manifold,
$\cT_{\R^d}[1]$, which is isomorphic to the tangent bundle on $\R^d$ with degrees of the fibers  shifted by $1$.
If $x^a$ are homogeneous coordinates on $\R^d$ and $\psi_a:=\p/\p x^a[1]$, then a polyvector field
$\ga\in \cT_{poly}(\R^d)$ is just a smooth function, $\ga(x,\psi)$, of these coordinates, and the Schouten
 brackets are given by,
$$
[\ga_1\bullet \ga_2]:= \Delta(\ga_1\ga_2) - \Delta(\ga_1)\ga_2 - (-1)^{\ga_1}\ga_1\Delta(\ga_2),
$$
where $\Delta=\sum_{a=1}^d (-1)^{|x^a|}\frac{\p^2}{\p x^a\p\psi_a}$. As $|\psi_a|=1-|x^a|$, the operator $\Delta$
and, therefore, the Schouten brackets  have degree $-1$.

\sip

We work in the category of semialgebraic
  sets and use in applications
  $PA$ differential forms on such sets, where $PA$ stands
for ``piecewise semi-algebraic" as defined in \cite{KoSo} and further developed in \cite{HLTV}. However, all our compactified
configuration spaces have also natural structures of smooth manifolds with corners which one can describe
explicitly in terms of metric graphs.

\mip
\section{Configuration space $C_n$}

{\bf 2.1. A Fulton-MacPherson type compactification of $C_n$ \cite{Ko}.}
Let
$$
\Conf_n:=\{z_1, \ldots, z_n\in \C\ |\, z_i\neq z_j\ \mbox{for}\ i\neq j\}
$$
be the configuration space of $n$ pairwise distinct points in the complex plane $\C$.
The space $C_{n}$ is a smooth $(2n-3)$-dimensional real manifold (or, if one prefers, a semialgebraic manifold) defined as the orbit space \cite{Ko},
 $$
C_{n}:=\Conf_n/{G},
$$
with respect to the following action of a real 3-dimensional Lie group,
$$
G= \{z\rar az+b \ |\ a\in \R^+, b\in \C\}.
$$
Its compactification, $\overline{C}_n$,
 was defined in \cite{Ko} (see also \cite{Ga}) as the closure of an embedding,
$$
\Ba{ccc}
C_n & \lon & (\R/2\pi\Z )^{n(n-1)}\times (\R\P^2)^{n(n-1)(n-2)}\\
(z_1, \ldots, z_n) & \lon & \prod_{i\neq j}\exp(Arg(z_i - z_j))\times \prod_{i\neq j \neq k\neq i}[|z_i-z_j: |z_i-z_k|: |z_j-z_k|].
\Ea.
$$
The space $\overline{C}_n$ is a smooth naturally oriented manifold with corners (it also has a natural structure of compact
oriented semialgebraic manifold).
Its codimension 1 strata is given by
$$
\p \overline{C}_n = \bigsqcup_{A\subset [n]\atop \# A\geq 2} C_{n - \# A + 1}\times
 C_{\# A}
$$
where the summation runs over all possible  proper subsets of $[n]$ with cardinality of at least two.
Geometrically, each such a strata corresponds to the $A$-labeled elements of the set $\{z_1, \ldots, z_n\}$ moving
very close
to each other.

\sip

{\bf 2.1.1. Configurations of ordered (or coloured) points.}
The natural action of the permutation groups on the standard {\em face complex}\, (which is a shorthand for the fundamental chain complex)
  of $\overline{C}_\bu$ is trivial
as permutations preserve all the big cells together with their natural orientations; it was observed in \cite{GJ} that
this face complex has a natural structure of an operad of $\caL ie_\infty$-algebras. However,
in the applications below the points of all configuration spaces  considered in this paper are always {\em coloured}, more precisely, they always come decorated
with vertices of certain graphs
(and such decorations extend naturally  to the compactifications).
The natural action of the permutation groups on the face complexes of  compactified spaces of such
{\em coloured}\,
configurations  is non-trivial   and hence can induce in principle  a
non-trivial $\bS_n$-action on the associated de Rham field theories (see \S 4 below). To keep this subtlety
under control
we assume from now  and until \S 5 that all our configuration spaces consist of not only distinct but also
{\em distinctly coloured}\,  points $(z_1,\ldots,z_n)$; equivalently, one may think of a choice of a total order
on the set  $(z_1,\ldots,z_n)$ (which may not coincide with the natural order induced by the integer labels) because such  a structure on $C_n$   extends naturally to its compactification
$\overline{C}_n$; our final formulae involve a summation over all possible graphs and hence over all possible decorations (in particular,
over all possible orderings) so that eventually nothing depends on such a choice. The face complex of
 {\em coloured}\, or {\em totally ordered} configuration spaces $\overline{C}_\bu$ has again a natural structure of a dg operad  which is different
from the operad $\caL ie_\infty$ and which we describe below.

\sip

{\bf 2.2. The face complex of $\{C_n\}$ as an operad of  Leibniz$_\infty$ algebras}.
The faces of $\overline{C}_n$ are isomorphic to the products of the form $C_{k_1}\times\ldots
\times C_{k_m}$.
The stratification of the (totally ordered) configuration space $\overline{C}_n$  is best coded by its face complex, $(\cF\cC hains(\overline{C}_\bu), \p)$,
 which is a dg free
operad, $\cF ree\langle E_\circ\rangle$,   generated by an $\bS$-module $E_{\circ}=\{E_{\circ}(n)\}$ with
\Beq\label{Ch2: white S-module}
E_\circ(n)=\left\{
\Ba{cr}
\C[\Sigma_n][2n-3]= \mbox{span}
\left(
\xy
(1,-6)*{\ldots},
(-13,-7)*{_{\sigma(1)}},
(-6.7,-7)*{_{\sigma(2)}},
(13,-7)*{_{\sigma(n)}},
 (0,0)*{\circ}="a",
(0,5)*{}="0",
(-12,-5)*{}="b_1",
(-8,-5)*{}="b_2",
(-3,-5)*{}="b_3",
(8,-5)*{}="b_4",
(12,-5)*{}="b_5",
\ar @{-} "a";"0" <0pt>
\ar @{-} "a";"b_2" <0pt>
\ar @{-} "a";"b_3" <0pt>
\ar @{-} "a";"b_1" <0pt>
\ar @{-} "a";"b_4" <0pt>
\ar @{-} "a";"b_5" <0pt>
\endxy
\ \
\right)_{\sigma\in \bS_n}
 & \mbox{for}\ n\geq 2 \vspace{6mm}\\
0 & \mbox{otherwise}.
\Ea
\right.
\Eeq
Each plain corolla with $n$ legs corresponds to $\overline{C}_n$. As we prefer working with {\em co}chain complexes, we assign to this corolla
degree $-(2n-3)=3-2n$. The boundary differential is given on the generators by
\Beqr
\p
\xy
(1,-5)*{\ldots},
(-13,-7)*{_1},
(-8,-7)*{_2},
(-3,-7)*{_3},
(7,-7)*{_{n-1}},
(13,-7)*{_n},
 (0,0)*{\circ}="a",
(0,5)*{}="0",
(-12,-5)*{}="b_1",
(-8,-5)*{}="b_2",
(-3,-5)*{}="b_3",
(8,-5)*{}="b_4",
(12,-5)*{}="b_5",
\ar @{-} "a";"0" <0pt>
\ar @{-} "a";"b_2" <0pt>
\ar @{-} "a";"b_3" <0pt>
\ar @{-} "a";"b_1" <0pt>
\ar @{-} "a";"b_4" <0pt>
\ar @{-} "a";"b_5" <0pt>
\endxy
&=&\sum_{k=0}^{n-2}\sum_{[n]\setminus [k+1]=I_1\sqcup I_2\atop
 \# I_1\geq 1}\ \ \
\begin{xy}
<0mm,0mm>*{\circ},
<0mm,0.8mm>*{};<0mm,5mm>*{}**@{-},
<-9mm,-5mm>*{\ldots},
<-9mm,-7mm>*{_{1\ \,  \ldots\  \, k}},
<0mm,-10mm>*{...},
<14mm,-5mm>*{\ldots},
<15mm,-7mm>*{\underbrace{\ \ \ \ \ \ \ \ \ \ \ }},
<15mm,-10mm>*{_{I_2}};
<-0.7mm,-0.3mm>*{};<-13mm,-5mm>*{}**@{-},
<-0.6mm,-0.5mm>*{};<-6mm,-5mm>*{}**@{-},
<0.6mm,-0.3mm>*{};<20mm,-5mm>*{}**@{-},
<0.3mm,-0.5mm>*{};<8mm,-5mm>*{}**@{-},
<0mm,-0.5mm>*{};<0mm,-4.3mm>*{}**@{-},
<0mm,-5mm>*{\circ};
<-5mm,-10mm>*{}**@{-},
<-2.7mm,-10mm>*{}**@{-},
<2.7mm,-10mm>*{}**@{-},
<5mm,-10mm>*{}**@{-},
<-6mm,-12mm>*{_{k+1}},
<1.2mm,-12mm>*{\underbrace{\ \ \ \ \ \ \ }},
<2.2mm,-15mm>*{_{I_1}},
\end{xy}
\nonumber\\
&=&
\sum_{A\varsubsetneq [n]\atop
\# A\geq 2} \ \
\begin{xy}
<0mm,0mm>*{\circ},
<0mm,0.8mm>*{};<0mm,5mm>*{}**@{-},
<-9mm,-5mm>*{\ldots},
<-9mm,-7mm>*{_{1\, ...\,  \inf\hspace{-0.7mm} A\hspace{-0.5mm}-\hspace{-0.5mm}1}},
<0mm,-10mm>*{...},
<14mm,-5mm>*{\ldots},
<15mm,-7mm>*{\underbrace{\ \ \ \ \ \ \ \ \ \ \ }},
<19mm,-10mm>*{_{[n-\inf A +1]\setminus A}};
<-0.7mm,-0.3mm>*{};<-13mm,-5mm>*{}**@{-},
<-0.6mm,-0.5mm>*{};<-6mm,-5mm>*{}**@{-},
<0.6mm,-0.3mm>*{};<20mm,-5mm>*{}**@{-},
<0.3mm,-0.5mm>*{};<8mm,-5mm>*{}**@{-},
<0mm,-0.5mm>*{};<0mm,-4.3mm>*{}**@{-},
<0mm,-5mm>*{\circ};
<-5mm,-10mm>*{}**@{-},
<-2.7mm,-10mm>*{}**@{-},
<2.7mm,-10mm>*{}**@{-},
<5mm,-10mm>*{}**@{-},
<0mm,-12mm>*{\underbrace{\ \ \ \ \ \ \ \ \ \ }},
<0mm,-15mm>*{_{A}},
\end{xy}
\label{Ch2: d on white corollas}.
\Eeqr

A semialgebraic (or smooth with corners) local coordinate chart at the face in $\overline{C}_{n}$ corresponding to
a graph $\Ga\in \cF ree\langle E_\circ\rangle$ can be described, as it is explained in detail
 in \S 5.2 of \cite{Ko}, by an associated {\em metric}\, graph, $\Ga^{met}$,
in which every internal edge, $e$, is assigned a small positive number $\var_e$.
For example, the face in $\overline{C}_7$ corresponding to a graph
$$
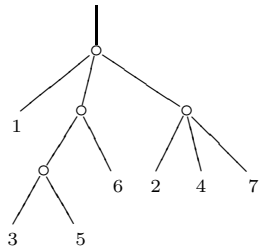
\xy
(-10.5,-2)*{_1},
(-11,-17)*{_3},
(-2,-17)*{_5},
(3,-10)*{_6},
(8,-10)*{_2},
(14,-10)*{_4},
(21,-10)*{_7},
(0,14)*{}="0",
 (0,8)*{\circ}="a",
(-10,0)*{}="b_1",
(-2,0)*{\circ}="b_2",
(12,0)*{\circ}="b_3",
(2,-8)*{}="c_1",
(-7,-8)*{\circ}="c_2",
(8,-8)*{}="c_3",
(14,-8)*{}="c_4",
(20,-8)*{}="c_5",
(-11,-15)*{}="d_1",
(-3,-15)*{}="d_2",
\ar @{-} "a";"0" <0pt>
\ar @{-} "a";"b_1" <0pt>
\ar @{-} "a";"b_2" <0pt>
\ar @{-} "a";"b_3" <0pt>
\ar @{-} "b_2";"c_1" <0pt>
\ar @{-} "b_2";"c_2" <0pt>
\ar @{-} "b_3";"c_3" <0pt>
\ar @{-} "b_3";"c_4" <0pt>
\ar @{-} "b_3";"c_5" <0pt>
\ar @{-} "c_2";"d_1" <0pt>
\ar @{-} "c_2";"d_2" <0pt>
\endxy
$$
has associated the following metric graph,
$$
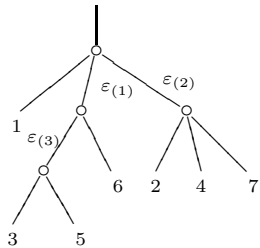
\xy
(3,2.7)*{_{\var_{(1)}}},
(11.2,3.7)*{_{\var_{(2)}}},
(-6.8,-4)*{_{\var_{(3)}}},
(-10.5,-2)*{_1},
(-11,-17)*{_3},
(-2,-17)*{_5},
(3,-10)*{_6},
(8,-10)*{_2},
(14,-10)*{_4},
(21,-10)*{_7},
(0,14)*{}="0",
 (0,8)*{\circ}="a",
(-10,0)*{}="b_1",
(-2,0)*{\circ}="b_2",
(12,0)*{\circ}="b_3",
(2,-8)*{}="c_1",
(-7,-8)*{\circ}="c_2",
(8,-8)*{}="c_3",
(14,-8)*{}="c_4",
(20,-8)*{}="c_5",
(-11,-15)*{}="d_1",
(-3,-15)*{}="d_2",
\ar @{-} "a";"0" <0pt>
\ar @{-} "a";"b_1" <0pt>
\ar @{-} "a";"b_2" <0pt>
\ar @{-} "a";"b_3" <0pt>
\ar @{-} "b_2";"c_1" <0pt>
\ar @{-} "b_2";"c_2" <0pt>
\ar @{-} "b_3";"c_3" <0pt>
\ar @{-} "b_3";"c_4" <0pt>
\ar @{-} "b_3";"c_5" <0pt>
\ar @{-} "c_2";"d_1" <0pt>
\ar @{-} "c_2";"d_2" <0pt>
\endxy
$$
describing an open subset of $\overline{C}_7$ consisting of all possible configurations
of $7$ points obtained in the following way:
\Bi
\item[(i)] take first a standardly positioned\footnote{The projection $\Conf_n\rar C_n$
has a natural section $s: C_n\rar \Conf_n$ representing each point
 $p\in C_n$ as a collection of $n$ pairwise distinct points $(z_1,\ldots,z_n)\in \Conf_n$
 such that the minimal Euclidean circle enclosing $(z_1,\ldots,z_n)$ has radius $1$ and center at $0\in \C$.
The point $s(p)$ is called the {\em standard position}\, of $p$ \cite{Ko}. If $\var$ is a positive
 real number, then the configuration $\var\cdot s(p):=\{\var z_1, \ldots,\var z_n\}\in \Conf_n$ is said to
 be the $\var$-{\em magnified} standard configuration.} configuration of 3 points labeled by $1$ and,
 say, $a$ and $b$,
\item[(ii)]
replace point $a$ (respectively, $b$) by an $\var_{(1)}$-magnified standardly positioned configuration of two points labeled by $c$ and $6$ (respectively, by  an $\var_{(2)}$-magnified
standard configuration of three points labeled by $2$, $4$ and $7$)
\item[(iii)] finally, replace the point $c$ by an $\var_{(3)}$-magnified  standard configuration of two points labeled
by $3$ and $5$, and project the resulting configuration in $\Conf_7$ into $C_7$.
\Ei
The embedding of the boundary faces into this smooth coordinate neighborhood is given by the
equation $\var_{(1)}\var_{(2)}\var_{(3)}=0$.
\mip

We conclude this subsection with a curious observation.
The values of differential (\ref{Ch2: d on white corollas})  on $2$- and and $3$-corollas are given by
$$
\p
\xy
(-3,-6)*{_1},
(3,-6)*{_2},
 (0,0)*{\circ}="a",
(0,4)*{}="0",
(-3,-4)*{}="b_1",
(3,-4)*{}="b_3",
\ar @{-} "a";"0" <0pt>
\ar @{-} "a";"b_3" <0pt>
\ar @{-} "a";"b_1" <0pt>
\endxy=0,
\ \ \
\p
\xy
(-5,-6)*{_1},
(-0,-6)*{_2},
(5,-6)*{_3},
 (0,0)*{\circ}="a",
(0,4)*{}="0",
(-4,-4)*{}="b_1",
(0,-4)*{}="b_2",
(4,-4)*{}="b_3",
\ar @{-} "a";"0" <0pt>
\ar @{-} "a";"b_2" <0pt>
\ar @{-} "a";"b_3" <0pt>
\ar @{-} "a";"b_1" <0pt>
\endxy
=
\Ba{c}
\xy
(-6.2,-10)*{_1},
(0.2,-10)*{_2},
(3.8,-5.8)*{_3},
 (0,0)*{\circ}="a",
(0,4)*{}="0",
(-3,-4)*{\circ}="b_1",
(3,-4)*{}="b_2",
(-6,-8)*{}="c_1",
(0,-8)*{}="c_2",
\ar @{-} "a";"0" <0pt>
\ar @{-} "a";"b_2" <0pt>
\ar @{-} "a";"b_1" <0pt>
\ar @{-} "b_1";"c_1" <0pt>
\ar @{-} "b_1";"c_2" <0pt>
\endxy
\Ea
+
\Ba{c}
\xy
(-6.2,-10)*{_1},
(0.2,-10)*{_3},
(3.8,-5.8)*{_2},
 (0,0)*{\circ}="a",
(0,4)*{}="0",
(-3,-4)*{\circ}="b_1",
(3,-4)*{}="b_2",
(-6,-8)*{}="c_1",
(0,-8)*{}="c_2",
\ar @{-} "a";"0" <0pt>
\ar @{-} "a";"b_2" <0pt>
\ar @{-} "a";"b_1" <0pt>
\ar @{-} "b_1";"c_1" <0pt>
\ar @{-} "b_1";"c_2" <0pt>
\endxy
\Ea
+
\Ba{c}
\xy
(6.2,-10)*{_3},
(-0.2,-10)*{_2},
(-3.8,-5.8)*{_1},
 (0,0)*{\circ}="a",
(0,4)*{}="0",
(3,-4)*{\circ}="b_1",
(-3,-4)*{}="b_2",
(6,-8)*{}="c_2",
(0,-8)*{}="c_1",
\ar @{-} "a";"0" <0pt>
\ar @{-} "a";"b_2" <0pt>
\ar @{-} "a";"b_1" <0pt>
\ar @{-} "b_1";"c_1" <0pt>
\ar @{-} "b_1";"c_2" <0pt>
\endxy
\Ea.
$$
Hence the map
$$
\al:\ \
(\cF ree\langle E_\circ\rangle, \p)\lon
\frac{\cF ree\left\langle \xy
(-3,-6)*{_1},
(3,-6)*{_2},
 (0,0)*{\circ}="a",
(0,4)*{}="0",
(-3,-4)*{}="b_1",
(3,-4)*{}="b_3",
\ar @{-} "a";"0" <0pt>
\ar @{-} "a";"b_3" <0pt>
\ar @{-} "a";"b_1" <0pt>
\endxy\right\rangle}
{
\left\langle
\Ba{c}
\xy
(-6.2,-10)*{_1},
(0.2,-10)*{_2},
(3.8,-5.8)*{_3},
 (0,0)*{\circ}="a",
(0,4)*{}="0",
(-3,-4)*{\circ}="b_1",
(3,-4)*{}="b_2",
(-6,-8)*{}="c_1",
(0,-8)*{}="c_2",
\ar @{-} "a";"0" <0pt>
\ar @{-} "a";"b_2" <0pt>
\ar @{-} "a";"b_1" <0pt>
\ar @{-} "b_1";"c_1" <0pt>
\ar @{-} "b_1";"c_2" <0pt>
\endxy
\Ea
+
\Ba{c}
\xy
(-6.2,-10)*{_1},
(0.2,-10)*{_3},
(3.8,-5.8)*{_2},
 (0,0)*{\circ}="a",
(0,4)*{}="0",
(-3,-4)*{\circ}="b_1",
(3,-4)*{}="b_2",
(-6,-8)*{}="c_1",
(0,-8)*{}="c_2",
\ar @{-} "a";"0" <0pt>
\ar @{-} "a";"b_2" <0pt>
\ar @{-} "a";"b_1" <0pt>
\ar @{-} "b_1";"c_1" <0pt>
\ar @{-} "b_1";"c_2" <0pt>
\endxy
\Ea
+
\Ba{c}
\xy
(6.2,-10)*{_3},
(-0.2,-10)*{_2},
(-3.8,-5.8)*{_1},
 (0,0)*{\circ}="a",
(0,4)*{}="0",
(3,-4)*{\circ}="b_1",
(-3,-4)*{}="b_2",
(6,-8)*{}="c_2",
(0,-8)*{}="c_1",
\ar @{-} "a";"0" <0pt>
\ar @{-} "a";"b_2" <0pt>
\ar @{-} "a";"b_1" <0pt>
\ar @{-} "b_1";"c_1" <0pt>
\ar @{-} "b_1";"c_2" <0pt>
\endxy
\Ea
\right\rangle
}
$$
which sends to zero all generating $n$-corollas except the ones with $n=2$ is a morphism of dg operads.
The operad
on the r.h.s.\ is a quadratic operad generated by $\C[\bS_2]$ modulo an ideal of relations
which is shown explicitly in the formula. This quadratic operad is in fact a
 well known Koszul operad \cite{Lo}, $\caL eib$, of Leibniz algebras\footnote{The author is grateful to
Michel van den Bergh for this observation.}
(with degree shifted by $-1$). Moreover, the minimal resolution of this operad was computed
in \cite{Li} and coincides precisely with $(\cF ree\langle E_\circ\rangle, d)$.\ Hence the map $\al$ is a quasi-isomorphism
giving us the following nice interpretation of the face complex of $\overline{C}_\bu$.

\mip

{\bf 2.2.1. Proposition}. {\em As a dg operad, the face complex, $\cF\cC hains(\overline{C}_\bu)$,
 of the family of
compactified configurations spaces,
$\{\overline{C}_n\}_{n\geq 2}$, is canonically  isomorphic to the minimal resolution, $\caL eib_\infty$,
 of the operad, $\caL eib$, of Leibniz algebras}.

\mip

Thus a structure of $\caL eib_\infty$-{\em algebra}\, on a dg vector space $(\fg, d)$ is given by a collection of
linear maps
$$
\left\{
\Ba{rccc}
\mu_n: & \ot^n\fg &\lon &  \fg[3-2n]\\
& \ga_1\ot \ldots\ot \ga_n     & \lon & \mu_n(\ga_1, \ldots, \ga_n),
\Ea
\right\}_{n\geq 2}
$$
satisfying the equations,
\Beq\label{Ch2:Leib_infty algebra}
(d\mu_n)(\ga_1, \ldots, \ga_n)=
\sum_{A\subsetneqq [n]\atop
\# A\geq 2}(-1)^{\sum_{k=1}^{\inf A-1}|\ga_k|}
\mu_{n-\# A +1}(\ga_1, \ldots, \ga_{\inf A-1}, \mu_{\# A}(\ga_A),\ga_{[n-\inf A]\setminus A}).
\Eeq
Here and elsewhere we use a notation,
$$
 \ga_S:= \ga_{i_1}\ot \ldots \ot \ga_{i_l},
$$
for a naturally ordered subset $S=\{i_i,\ldots,i_l\}$ of $[n]$. If the tensors $\mu_n$
happen to be graded symmetric, $\mu_n: \odot^n\fg\rar \fg[3-2n]$,
 then the above equation is precisely the equation for a $\caL ie_\infty$
structure on $\fg$, i.e. there exists a canonical morphism of dg operads,
$$
\caL eib_\infty \lon \caL ie_\infty.
$$
Thus any Lie algebra, $\fg$, is also a Leibniz
algebra. However, the groups of automorphisms of $\fg$ in the categories of $\caL ie_\infty$-algebras and of
$\caL eib_\infty$-algebras may be different.
\sip

It is worth noting that the
 symmetrization of a generic $\caL eib_\infty$-algebra structure, $\{\mu_n \}$,
does {\em not}\, give, in general, a $\caL ie_\infty$-structure, i.e.\ there is no simple morphism of dg operads
of the form,   $\caL ie_\infty \rar \caL eib_\infty$.

\mip

{\bf 2.3. A  semialgebraic structure on $\overline{C}_n$.}
 For future reference we shall make now the topological compactification $\overline{C}_n$ into a semialgebraic set by representing each point $p\in C_n$ by a  configuration of points in the upper half plane.
  Let $\Conf_A(\C)$ stand for the
 space of immersions, $A\hook \C$, of a finite non-empty set $A$ into the complex plane and
$\widetilde{\Conf}_A(\C)$ for the space of all possible maps, $A\rar \C$. For $\# A\geq 2$ we define the quotient space,
$C_A:=\frac{\Conf_A(\C)}{G}$,
so that $C_{[n]}=C_n$.

\sip

For any given a point $z_0\in \bbH:=\{x+\ii y\in \C\, |\, y>0\}$  and any  positive real number $\var$ there is an
associated  section,
$$
\Ba{rccc}
s_{(z_0,\var)}: & C_A & \lon &  \Conf_A(\C) \\
   & p=\{z_i=x_i+\ii y_i\}_{i\in A} & \lon & p_{(z_0,\var)}:= \var \frac{p - z_{min}(p)}{|p - z_{min}(p)|}+ z_0
   \Ea
$$
of the natural projection $\Conf_A(\C)\rar C_A$,
where
$$
z_{min}(p):=  \frac{1}{\# A}\sum_{i\in A} x_i + \ii \inf_{i\in A} y_i
$$
The image, $C_A(z_0, \var):=s_{(z_0,\var)}(C_A)$, consists of configurations
$p=\{z_i\}_{i\in A}$ in $\Conf_A(\C)$ satisfying two conditions
\Bi
\item[$(i)$]  $z_{min}(p)=z_0$ (so that $z_i\in \bbH$ for all $i\in A$);
\item[$(ii)$] $\sqrt{\sum_{i\in A}|z_i-z_0|^2}=\var$.
\Ei

If $z_0=i$ and $\var=1$ we use a simpler notation,
$C_n^h:= C_n(\ii,1)$, e.g.
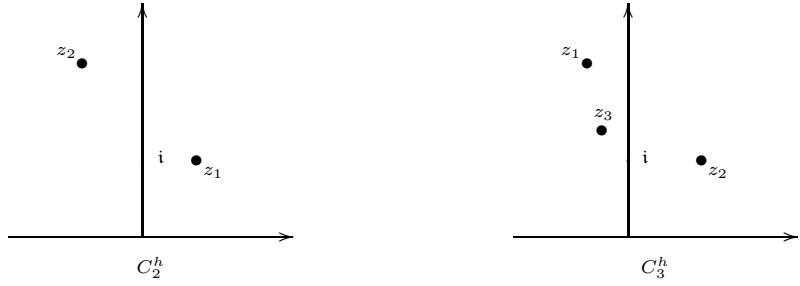
$$
\underset{C_2^h}
{
\Ba{c}
{
\xy
(-8,23)*{\bullet}="e",
(-10,24.5)*{_{z_2}},
(7.2,10)*{\bullet}="f",
(9.5,8.5)*{_{z_1}},
(0,10)*{\cdot},
(0.2,10)*{\ \ \ \ ^{\ii}},
 (-17.8,0)*{}="a",
(20,0)*{}="b",
(0,0)*{}="c",
(0,31)*{}="d",
\ar @{->} "a";"b" <0pt>
\ar @{->} "c";"d" <0pt>
\endxy}
\Ea
}
\ \ \ \ \ \ \ \ \ \ \ \ \ \ \ \ \ \
\ \ \ \
\underset{C_3^h}
{
\Ba{c}
{
\xy
(-8,23)*{\bullet}="e",
(7.2,10)*{\bullet}="f",
(-10,24.5)*{_{z_1}},
(9.5,8.5)*{_{z_2}},
(-6,14)*{\bullet},
(-5.7,15.8)*{^{z_3}},
(-2.5,10)*{\cdot},
(-2.5,10)*{\ \ \ \ ^{\ii}},
 (-17.8,0)*{}="a",
(20,0)*{}="b",
(-2.5,0)*{}="c",
(-2.5,31)*{}="d",
\ar @{->} "a";"b" <0pt>
\ar @{->} "c";"d" <0pt>
\endxy}
\Ea
}
$$
The upper half space representation, $p^h:=                                           p_{(\ii,1)}\in C_n^h$, of a point, $p\in C_n$,
is called its {\em standard
hyperbolic position}\, (or, shortly, {\em hyperposition}).
Note that $C_n^{h}$ is a union, $\bigcup_{i=1}^n\{p^{h}\in C_n^h\, |\, y^i=1\}$,
of semialgebraic sets and hence is itself a semialgebraic set (of dimension $2n-3$).
 For example, $C_2^h\simeq S^1$ is, as a semialgebraic set, the union of two  semialgebraic intervals,
$$
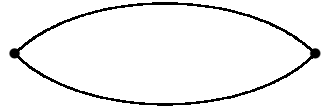
{
\xy
(-20,0)*{\bullet},
(20,0)*{\bullet},
   \ar@{-}@(ur,lu) (-20.0,0.0)*{};(20.0,0.0)*{},
   \ar@{-}@(dr,ld) (-20.0,0.0)*{};(20.0,0.0)*{},
\endxy}
$$
with the lower interval corresponding to the configurations
given in the left picture above and the upper interval
to the similar configurations but with $z_1$ and $z_2$ swapped; the left
and right $\bullet$-points represent, respectively,  the following two configurations,
$$
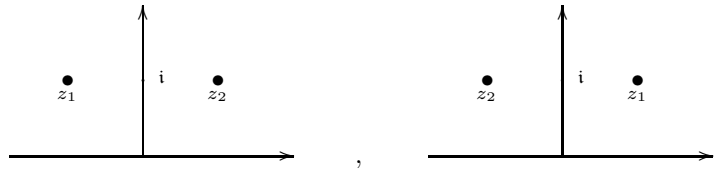
\xy
(-10,10)*{\bullet}="e",
(-10,8)*{_{z_1}},
(10,10)*{\bullet}="f",
(10,8)*{_{z_2}},
(0,10)*{\cdot},
(0.2,10)*{\ \ \ \ ^{\ii}},
 (-17.8,0)*{}="a",
(20,0)*{}="b",
(0,0)*{}="c",
(0,20)*{}="d",
\ar @{->} "a";"b" <0pt>
\ar @{->} "c";"d" <0pt>
\endxy
\ \ \ \ \ \ \
,
\ \ \ \ \ \ \
\xy
(-10,10)*{\bullet}="e",
(-10,8)*{_{z_2}},
(10,10)*{\bullet}="f",
(10,8)*{_{z_1}},
(0,10)*{\cdot},
(0.2,10)*{\ \ \ \ ^{\ii}},
 (-17.8,0)*{}="a",
(20,0)*{}="b",
(0,0)*{}="c",
(0,20)*{}="d",
\ar @{->} "a";"b" <0pt>
\ar @{->} "c";"d" <0pt>
\endxy
$$

\sip

It is clear from the construction that the configuration spaces $C_n(z_0,\var)\subset Conf_n$ for different values of the parameters $z_0$ and
  $\var$ are canonically isomorphic to each other,
$$
\Phi_{(z_1,\var_1)}^{(z_2,\var_2)}: C_n(z_1,\var_1)\lon C_n(z_2,\var_2),
$$
for some uniquely determined element $\Phi_{(z_1,\var_1)}^{(z_2,\var_2)}\in G$.

\sip
 The map,
\Beq\label{Ch2: magnific-1}
\Ba{rccc}
\var\cdot:=\Phi_{(\ii,1)}^{(\ii,\var)}: & C_n^h & \lon & C_n(\ii,\var)\\
&  p_{(\ii,1)} & \lon & p_{(\ii,\var)},
\Ea
\Eeq
is called $\var$-{\em magnification}, e.g., for $\var<1$,
\bip
$$
\underset{\mbox{\sc Figure 1}}{
\Ba{c}
\underset{C_3^h}
{
\xy
(0,17)*{
\xycircle(10,10){.}},
(5,22.1)*{\bullet}="g",
(-2,26.5)*{\bullet}="e",
(2.2,7.1)*{\bullet}="f",
(1.1,9.3)*{\ \ \ \ ^{\ii}},
 (-17.8,0)*{}="a",
(20,0)*{}="b",
\ar @{->} "a";"b" <0pt>
\ar @{.} "e";"f" <0pt>
\ar @{.} "e";"g" <0pt>
\ar @{.} "g";"f" <0pt>
\endxy
}
\Ea
\ \ \ \ \
\stackrel{\var\cdot}{\lon}
\ \ \ \ \
\Ba{c}
\underset{C_{3}}
{
\xy
(0,17)*{
\xycircle(10,10){.}},
(1,14)*{
\xycircle(7,7){}},
(5,22.1)*{}="gg",
(-2,26.5)*{}="ee",
(3.9,17.5)*{\bullet}="g",
(-0.8,20.5)*{\bullet}="e",
(2.2,7.1)*{\bullet}="f",
(1.1,9.3)*{\ \ \ \ \ \ ^{\ii}},
 (-17.8,0)*{}="a",
(20,0)*{}="b",
\ar @{->} "a";"b" <0pt>
\ar @{.} "e";"g" <0pt>
\ar @{.} "ee";"f" <0pt>
\ar @{.} "gg";"f" <0pt>
\ar @{.} "gg";"ee" <0pt>
\endxy
}
\Ea
}
$$
It is given explicitly by the affine map $z\rar \var(z-\ii) + \ii$.

\mip

Finally, let
$p\in C_{k}(\ii,\var)$ and  $p'\in C_{n}(\ii,\var')$ be magnified standard configurations with the ratio
 $\var/\var' \gg 0$. With
any point $z$ in the configuration  $p$                                                                                                        we associate a translation map,
\Beq\label{2: map T_z}
\Ba{rccc}
T_{z}: & C_n(\ii,\var')& \lon &  C_n(z,\var')\\
& p' &\lon & T_z(p')= p' + z - \ii
\Ea
\Eeq
which preserves all the relative Euclidean angles of the points in $p'$ and moves $z_{min}(p')$
from $\ii$ to $z$.
The image of $T_z(p')$ under the projection $Conf_n\rar C_n$ is called
 {\em the configuration  $p'$ placed at the point $z$
of the configuration $p$}.

\sip

%

With these new notions of a {\em standard position}, a {\em magnification}\,
and {\em placing of a magnified standard configuration at a given point of another
magnified standard configuration}\,
one can apply the same idea of metric graphs as in the previous subsection to define a
 semialgebraic atlas
 on $\overline{C}_n$.  It is worth noting that
 this atlas makes sense not only for  {\em small} values of $\var$-parameters  but also for {\em large}\, ones. Indeed, the boundary strata are given, in the above notations, by equations of the form $\var'/\var=0$, and this can be achieved,
for example,  by keeping $\var$ finite and $\var'\rar 0$ and also by keeping $\var'$  finite and $\var\rar +\infty$ (or, better, when both $\var$ and $\var'$ tend to $+\infty$ in such a way that
$\var/\var'\gg 0$).

\sip

In fact, the most suitable for our purposes semialgebraic structure on $\overline{C}_n$
can be described directly, without any reference to metric graphs.  Let $\widetilde{C}^{h}_A$ be a subset of $\widetilde{\Conf}_A(\C)$  consisting of all
possible configurations  $p=\{z_i\}_{i\in A}$  satisfying the  following two conditions,
\Bi
\item[$(i)$]  $z_{min}(p)=\ii$;
\item[$(ii)$] ${\sum_{i\in A}|z_i-\ii|^2}=1$.
\Ei
 Note that
 $\widetilde{C}^{h}_A$ is a {\em compact}\, semialgebraic set in $\R^{2\# A}$ containing $C_A^{h}$ as a semialgebraic subset.
Note also that
$C_2^{h}= \widetilde{C}^{h}_2$.
 It is not hard to check that
the topological  compactification, $\overline{C}_n$, of $C_n$ can be defined as the  closure
of the following embedding (cf.\ \cite{AT})
$$
C_n \stackrel{\prod \pi_A}{\lon} \prod_{A\subseteq [n]\atop \# A\geq 2} C_A
 \stackrel{\simeq }{\lon} \prod_{A\subseteq [n]\atop \# A\geq 2} C_A^{h}
\hook  \prod_{A\subseteq [n]\atop \# A\geq 2} \widetilde{C}^{h}_A.
$$
where the product runs over all possible subsets $A$ of $[n]$ with $\# A\geq 2$, and
$$
\Ba{rccc}
\pi_A : & C_n & \lon & C_A\\
        & p=\{z_i\}_{i\in [n]} & \lon & p_A:=\{z_i\}_{i\in A}
\Ea
$$
is the natural forgetful (semialgebraic) map. The r.h.s.\ above is a compact semialgebraic set
so that  $\overline{C}_n$ comes equipped naturally with an induced semialgebraic structure which we
declare from now a default one. Thus for us $\overline{C}=\{\overline{C}_n\}_{n\geq 2}$ is a non-unital operad in the category of  semialgebraic sets. If configurations in $C_n$ are totally ordered, then
the dg suboperad, $\cF \cC hains(\overline{C})\subset \cC hains(\overline{C})$, of semialgebraic fundamental chains is isomorphic to the dg operad $\caL eib_\infty$, if not, then to $\caL ie_\infty$.

\bip

\section{Configuration space $C_{n,0}$ and its new compactifications}

\sip

{\bf 3.1. ${C}_{n,0}$ as a magnified ${C}_n$}.
  Let $\Conf_{A,0}(\bbH)$ stand for the set of
all injections, $A\hook \bbH$, of a finite set $A$ into the upper half-plane. In this section we study several  compactifications, ${\widehat{C}}_{A,0}$,  ${\widehat{C}}_{A,0}'$ and
${\widehat{C}}_{A,0}''$, of the quotient configuration space,
$$
C_{A,0}:=\frac{\Conf_{A,0}(\bbH)}{\{z\rar \R^+ z+ \R\}}, \ \ \ \ \# A\geq 1,
$$
which
are different from each other as {\em semialgebraic manifolds}\, but which  give us three  geometric models  for {\em one and the same}\, 2-coloured operad,
 $\cM or(\caL ie_\infty)$, of strong homotopy morphisms of  $\caL ie_\infty$-algebras. If the sets
 $A$ are totally ordered, then one gets instead geometric models for the 2-coloured operad
  $\cM or(\caL eib_\infty)$, of strong homotopy morphisms of  $\caL eib_\infty$-algebras.
These  compactifications have little in common with Kontsevich's compactification, $\overline{C}_{n,0}$, of $C_{n,0}$ studied in \cite{Ko2}.
 \sip

Define a section,
$$
\Ba{rccc}
h: &C_{A,0} & \lon & \Conf_{A,0}(\bbH)\\
\displaystyle
&\displaystyle p={\{z_i=x_i+\ii y_i\in \bbH\}_{i\in A}} & \lon &
\displaystyle p^{h}:=  \frac{p - x_c(p)}{\inf_{i\in A}y_i}.
\Ea
$$
where $x_c(p):=\frac{1}{\# A}\sum_{i=1}^{\# A}x_{i}$, and set $C_{A,0}^{h}:=\Img h$. For an arbitrary  $p\in C_{A,0}$ we denote
$$
||p||:=|p^{h}-\ii|=\frac{|p-z_{min}(p)|}{\inf_{i\in A}y_i},
$$
where, as before, $z_{min}(p)= x_c(p) + \ii \inf_{i\in A}y_i$.
 Note that every point in the configuration $p^{h}$ lies
in the subspace $\Im z \geq 1 \subset \bbH$ and at least one point lies on the line $\Im z=1$.
Thus
$$
C_{A,0}^{h}=\left\{p=\{z_i\}_{i\in A} \in  \Conf_{A,0}(\bbH)\ |\ z_{min}(p)=\ii
\right\}.
$$
and $C_A^{h}$ is a subspace of $C_{A,0}^{h}$ consisting of configurations $p^{h}$ with $|p^{h}-\ii|=1$.
Note that $C_{A,0}^{h}$ is a union, $\bigcup_{i\in A}\{p^{h}\in C_{A,0}^h\, |\, y^i=1\}$,
of semialgebraic sets and hence is itself a semialgebraic set (of dimension $2\# A-2$), and we make
$C_{A,0}$ into a semialgebraic set by identifying it with $C_{A,0}^{h}$.

\sip

 There is a semialgebraic homeomorphism,
\Beq\label{5': X_n iso of C(H)}
\Ba{rccccc}
\Psi_n: & C_{n,0}& \lon &  C_{n}^{h} & \times & (0,1)  \\
& p & \lon & \frac{p- z_{min}(p)}{|p-z_{min}(p)|} +\ii &\times & \frac{||p||}{||p||+1}.
\Ea
\Eeq
It is worthwhile pointing out that in geometric terms the inverse isomorphism,
$$
\Ba{ccc}
 C_{n}^h \times \R^+ &\lon& C^h_{n,0} \\
 (p_0, \la) & \lon & p:=\la\cdot p_0,
\Ea
$$
is given by ``exploding"\footnote{We write here ``exploding" instead of ``expanding" in the anticipation of the main trick of the new
compactification which formally allows $\la$ to be equal not only to zero (as in \S 2)
but also to $+\infty$.} a diameter $1$ configuration $p_0\in C_n^h$ by the factor $\la$.
Here we used a  magnification map (cf.\ (\ref{Ch2: magnific-1})),
$$
\Ba{ccc}
C_{n,0} \times \R^+&\lon & C_{n,0}\\
(p, \var)  &\lon& \var\cdot p
\Ea
$$
defined  as the composition,
\Beq
\Ba{rccccccc}\label{Ch3: magnification map on C_n,0}
\var\cdot: & C_{n,0} \simeq C_{n,0}^h& \stackrel{\Psi_n}{\lon} & C_{n}^h\times \R^+  & \lon & C_{n}^h\times \R^+ &
\stackrel{\Psi_n^{-1}}{\lon} & C_{n,0}\\
          &  p    &       \lon    & (p_0, \la) & \lon   & (p_0, \var \la) & \lon& \var\cdot p.
\Ea
\Eeq
It is worth noting that this magnification map preserves all the  Euclidean angles
$Arg(z_i-z_j)$ and all the relative Euclidean  distances, $\frac{|z_i-z_j|}{|z_k-z_l|}$, of points
in a configuration $p\in C_{n,0}$ as it can be represented by an action of an element of the group $G$
on $Conf_{n,0}$.

\mip
{\bf 3.2. New compactifications of $C_{n,0}$.}
With a subset $A\subset [n]$ one associates a forgetful map,
$$
\Ba{rccc}
\pi_A: &  C_{n,0}(\bbH) & \lon & C_{A,0}(\bbH)\\
       &  p=\{z_i\}_{i\in [n]} & \lon &  p_A=\{z_i\}_{i\in A}.
\Ea
$$
A semialgebraic {\em compactification}, $\widehat{C}_{n,0}'$, of $C_{n,0}$
is defined  as the closure of the semialgebraic injection
\Beq\label{2: first compactifn of fC(H)}
C_{n,0}(\bbH)\stackrel{\prod \pi_A}{\lon} \prod_{A\subseteq [n]\atop A\neq \emptyset} C_{A,0}
\stackrel{\prod \Psi_A}{\lon} \prod_{A\subseteq [n]\atop A\neq \emptyset}
 C_{A}^{h}\times (0, 1) \lon \prod_{A\subseteq [n]\atop A\neq \emptyset}
 \widetilde{C}^{h}_{A}\times [0,1].
\Eeq
equipped with the induced structure of a semialgebraic set.

\sip

Next, for a pair of subspaces $B\subsetneq A\subseteq [n]$  we consider two maps
$$
\Ba{rcccccc}
\Psi_{A,B}: & C_{n,0} & \lon &  C_B^{h}&\times & (0,1)\\
& p & \lon &  \frac{p_B-x_c(p_B)}{||p_B||}
&& \frac{||p_{A,B}||}{||p_{A,B}||+1}
\Ea
$$
and
$$
\Ba{rcccccc}
\Psi_{A,B}'': & C_{n,0} & \lon &  C_B^{h}&\times & (0,1)\\
& p & \lon &  \frac{p_B-x_c(p_B)}{||p_B||}
&& \frac{||p_{A,B}||''}{||p_{A,B}||''+1}
\Ea
$$

where
$$
||p_{A,B}||:=\frac{y_{min}(p_B)}{y_{min}(p_A)}
||p_B||=  \frac{|p_B-z_{min}(p_B)|}{y_{min}(p_A)}
\ \ \ \mbox{and}\ \ \ \ ||p_{A,B}||'':=||p_A||\cdot ||p_B||.
$$
 A semialgebraic {\em compactification}, $\widehat{C}_{n,0}$, of $C_{n,0}$
is defined as the closure of the following semialgebraic map
\Beq\label{2: second compactifn of fC(H)}
C_{n,0}\stackrel{\Psi_n\times \prod \Psi_{A,B}}{\lon} C_n^{st}\times (0,1)\times
\hspace{-1mm}  \prod_{B\subsetneq A\subseteq [n]\atop \# B\geq 2} \hspace{-2mm} \left(C_B^{h} \times (0,1)\right)
 \hook \widetilde{C}_n^{h}\times [0,1]\times \hspace{-3mm}\prod_{B\subsetneq A\subseteq [n]\atop \# B\geq 2}
 \hspace{-2mm}  \left(\widetilde{C}_B^{h} \times [0,1]\right).
\Eeq
and  a {\em compactification}, $\widehat{C}_{n,0}''$, of $C_{n,0}$
is defined similarly by replacing the  maps $\Psi_{A,B}$ with $\Psi''_{A,B}$.
 The induced semialgebraic structures in
$\widehat{C}_{n,0}$,  $\widehat{C}_{n,0}'$ and $\widehat{C}_{n,0}''$ can also be  described explicitly with the help of metric graphs (see below); in fact,  one could use metric graphs to make these spaces into smooth manifolds with corners but  we never use such smooth structures on these spaces as all  propagators we work with in this paper give rise to PA differential forms which are only piecewise smooth.

\sip

The boundary strata in all three compactifications are given   by the limit values $0$ or $+\infty$
of the parameters $||p_A||$, $||p_{A,B}||$ and, respectively,  $||p_{A,B}||''$. It is an easy exercise to check that in all  three cases
the codimension 1 boundary strata are given by\footnote{We refer to \cite{Me-Config} for a detailed discussion and explicit examples.}

\Beq\label{Ch3:codimension 1 boundary}
\displaystyle
\p \left\{\Ba{c} \widehat{C}_{n,0}\\ \widehat{C}_{n,0}'\\
 \widehat{C}_{n,0}''
\Ea\right\} = \bigsqcup_{A\subseteq [n]\atop \# A\geq 2} \left(C_{n - \# A + 1,0}\times
 C_{\# A}(\C)\right)\
 \bigsqcup_{[n]=B_1\ \sqcup \ldots \sqcup B_k\atop{
 2\leq k\leq n \atop \#B_1,\ldots, \#B_k\geq 1}}\left( C_{k}\times C_{\# B_1,0}\times \ldots\times
 C_{\# B_k,0}
 \right)
\Eeq
where the first summation runs over all possible  subsets, $A$, of $[n]$
 with cardinality  at least two, and the
second summation runs over all possible decompositions of $[n]$ into (at least two)  disjoint
non-empty subsets $B_1, \ldots, B_k$.
Geometrically,  a stratum in the first group of summands  corresponds to  $A$-labeled elements of the set $\{z_1, \ldots, z_n\}$ moving {\em close}\,
to each other, while a stratum in the second group of summands corresponds to $k$ clusters of points (labeled, respectively,
 by disjoint  subsets $B_1, \ldots B_k$ of $[n]$) moving {\em far}\,  from each other as in the picture below
$$
 \Ba{c}
{
\xy
 (0,2.5)*\ellipse(40,40),=:a(180){.};
 (-15,2.5)*\ellipse(7,7),=:a(180){.};
(-35,5)*{\bullet},
(-27,8)*{\bullet},
(-32,10)*{\bullet},
(-30,7)*{\bullet},
(-28,13)*{^{B_i}},
 (18,33)*{
\xycircle(5.5,5.5){.}};
(15,33)*{\bullet},
(19,37)*{\bullet},
(21,35)*{\bullet},
(18,29)*{\bullet},
(13,28)*{^{B_l}},
 (-17,35)*{
\xycircle(5,5){.}};
(-15,33)*{\bullet},
(-19,37)*{\bullet},
(-16,35)*{\bullet},
(-15,28)*{^{B_p}},
 (32,13)*{
\xycircle(5.5,5.5){.}};
(35,13)*{\bullet},
(30,17)*{\bullet},
(31,8)*{\bullet},
(32,14)*{\bullet},
(25,8)*{^{B_q}},
 (-50,0)*{}="a",
(50,0)*{}="b",
(0,0)*{}="c",
(0,50)*{}="d",
(-40,5)*{}="1",
(40,5)*{}="2",
(-29.7,32.7)*{}="3",
(-34,37)*{}="4",
(29.7,32.7)*{}="3'",
(34,37)*{}="4'",
\ar @{->} "a";"b" <0pt>
\ar @{->} "c";"d" <0pt>
\ar @{.} "1";"2" <0pt>
\ar @{.>} "3";"4" <0pt>
\ar @{.>} "3'";"4'" <0pt>
\endxy}
\Ea
$$
The only difference between the three compactifications $\widehat{C}_{n,0}$, $\widehat{C}_{n,0}'$ and $\widehat{C}_{n,0}''$ is at the infinity strata:
\Bi
\item
in the case of $\widehat{C}_{n,0}'$ the codimension 1 infinity strata  are given by the equations $||p||=+\infty$, $||p_{B_i}||$ a finite number for all  $i\in[k]$ (here $||p||$ is the ``size" of the whole configuration $p=\cup_{i=1}^k{B_i}$ and $||p_{B_i}||$
is the ``size" of the subgroup $B_i$), i.e.\ in this case the groups of points $\{B_i\}_{i\in [k]}$ tend to infinite Euclidean distance from each other
while keeping their sizes $||p_{B_i}||$ finite,
\item in the case of $\widehat{C}_{n,0}''$ the codimension 1 infinity strata  is given by the equations $||p||=+\infty$, $||p||\cdot ||p_{B_i}||$ a finite number for all $i\in [k]$, i.e.\ in this case the groups of points $\{B_i\}_{i\in [k]}$ tend to infinite Euclidean distance from each other
while having their sizes $||p_{B_i}||$ decreasing with the speed $\sim ||p||^{-1}$.
\item in the case of $\widehat{C}_{n,0}$ the codimension 1 infinity strata  are given by the equations $||p||=+\infty$, $\frac{|p_{B_i}-z_{min}(p_{B_i})|}{y_{min}(p)}$ a finite number; as $y_{min}(p)$ is always
    finite and, moreover,
    is always normalized  in our compactification embedding to be $1$, the latter condition simply says that
     the groups of points $\{B_i\}_{i\in [k]}$ tend to infinite Euclidean distance from each other
     while keeping their Euclidean sizes   $|p_{B_i}-z_{min}(p_{B_i})|$ finite.
\Ei
The disjoint unions,
$$
\widehat{C}:=\overline{C}_\bu \sqcup \widehat{C}_{\bu,0}\sqcup \overline{C}_\bu\  , \ \ \ \ \widehat{C}':=\overline{C}_\bu \sqcup \widehat{C}_{\bu,0}'\sqcup \overline{C}_\bu \ \ \ \ \mbox{and}\ \ \ \
\widehat{C}'':=\overline{C}_\bu \sqcup \widehat{C}_{\bu,0}\\\sqcup \overline{C}_\bu
$$
have natural structures of 2-coloured operads in the category of compact semialgebraic sets; if we forget the semialgebraic structures and view them as operads of sets, then all the three 3 operads
are identical to each other and to the free 2-coloured operad generated by the set
${C}_\bu \sqcup {C}_{\bu,0}\sqcup {C}_\bu$.

\sip

We shall be most interested in this paper in the semialgebraic compactification $\widehat{C}_{n,0}$.
By analogy to $\overline{C}_{n}$,
the face stratification of  $\widehat{C}_{n,0}$ can be nicely described in terms of graphs (in fact, in terms of a
dg operad describing strongly homotopic morphisms of $\caL eib_\infty$-algebras) while the  semialgebraic structure
on $\widehat{C}_{n,0}$ is best described in terms of the associated {\em metric}\, graphs. These are the main
themes of the rest of this section.

\mip

{\bf 3.3. The face complex of $\widehat{C}_{n,0}$ as an operad $\cM or(\caL eib_\infty)$.} Note that all the boundary faces of the compactification  $\widehat{C}_{n,0}$ are products of the form
$$
\overline{C}_{p_1}\times\ldots \overline{C}_{p_k}\times \widehat{C}_{n_1,0}\times \ldots \times \widehat{C}_{n_k,0}\times \overline{C}_{q_1}\times
\times \ldots \times \overline{C}_{q_l}
$$
where factors $\overline{C}_{q_\bullet}$ correspond to the first group of boundary terms in
(\ref{Ch3:codimension 1 boundary}), that is, to the collapsing strata, and the factors $\overline{C}_{p_\bullet}$ to the second one, i.e.\ to infinity strata.
There are no boundary strata which would contain mixed products of  $\overline{C}_{q_\bullet}$ and
$\overline{C}_{p_\bullet}$ only. This fact forces us to interpret these two types of boundary strata as  operads in {\em two}\, different colours and represent the associated generators, $\overline{C}_p$ and, respectively, $\overline{C}_q$, by {\em different}\, corollas,   say, by {\em dashed}\, white corollas
and, respectively,  {\em solid}\, white corollas,
$$
\overline{C}_p \simeq
\xy
(1,-5)*{\ldots},
(-13,-7)*{_1},
(-8,-7)*{_2},
(-3,-7)*{_3},
(7,-7)*{_{p-1}},
(13,-7)*{_p},
 (0,0)*{\circ}="a",
(0,5)*{}="0",
(-12,-5)*{}="b_1",
(-8,-5)*{}="b_2",
(-3,-5)*{}="b_3",
(8,-5)*{}="b_4",
(12,-5)*{}="b_5",
\ar @{.} "a";"0" <0pt>
\ar @{.} "a";"b_2" <0pt>
\ar @{.} "a";"b_3" <0pt>
\ar @{.} "a";"b_1" <0pt>
\ar @{.} "a";"b_4" <0pt>
\ar @{.} "a";"b_5" <0pt>
\endxy,
\ \ \ \ \
\overline{C}_q \simeq
\xy
(1,-5)*{\ldots},
(-13,-7)*{_1},
(-8,-7)*{_2},
(-3,-7)*{_3},
(7,-7)*{_{q-1}},
(13,-7)*{_q},
 (0,0)*{\circ}="a",
(0,5)*{}="0",
(-12,-5)*{}="b_1",
(-8,-5)*{}="b_2",
(-3,-5)*{}="b_3",
(8,-5)*{}="b_4",
(12,-5)*{}="b_5",
\ar @{-} "a";"0" <0pt>
\ar @{-} "a";"b_2" <0pt>
\ar @{-} "a";"b_3" <0pt>
\ar @{-} "a";"b_1" <0pt>
\ar @{-} "a";"b_4" <0pt>
\ar @{-} "a";"b_5" <0pt>
\endxy,
$$
Next we represent a boundary factor of the type $\widehat{C}_{n,0}$  as a degree $2-2n$ black vertex corolla,
$$
\widehat{C}_{n,0}\simeq
\xy
(1,-5)*{\ldots},
(-13,-7)*{_1},
(-8,-7)*{_2},
(-3,-7)*{_3},
(7,-7)*{_{n-1}},
(13,-7)*{_n},
 (0,0)*{\bullet}="a",
(0,5)*{}="0",
(-12,-5)*{}="b_1",
(-8,-5)*{}="b_2",
(-3,-5)*{}="b_3",
(8,-5)*{}="b_4",
(12,-5)*{}="b_5",
\ar @{.} "a";"0" <0pt>
\ar @{-} "a";"b_2" <0pt>
\ar @{-} "a";"b_3" <0pt>
\ar @{-} "a";"b_1" <0pt>
\ar @{-} "a";"b_4" <0pt>
\ar @{-} "a";"b_5" <0pt>
\endxy
$$
After these notational preparations one can describe the face complex of the compactification $\widehat{C}_{\bullet,0}$ as follows.

\mip

{\bf 3.3.1. Proposition.} {\em The face complex of the disjoint union $\overline{C}_\bu\sqcup \widehat{C}_{\bu,0}\sqcup
\overline{C}_\bu$ has naturally a structure of a free 2-coloured operad,
$$
\cM or(\caL eib_\infty):= \cF ree
\left\langle
\xy
(1,-5)*{\ldots},
(-13,-7)*{_1},
(-8,-7)*{_2},
(-3,-7)*{_3},
(7,-7)*{_{p-1}},
(13,-7)*{_p},
 (0,0)*{\circ}="a",
(0,5)*{}="0",
(-12,-5)*{}="b_1",
(-8,-5)*{}="b_2",
(-3,-5)*{}="b_3",
(8,-5)*{}="b_4",
(12,-5)*{}="b_5",
\ar @{-} "a";"0" <0pt>
\ar @{-} "a";"b_2" <0pt>
\ar @{-} "a";"b_3" <0pt>
\ar @{-} "a";"b_1" <0pt>
\ar @{-} "a";"b_4" <0pt>
\ar @{-} "a";"b_5" <0pt>
\endxy,
\ \ \
\xy
(1,-5)*{\ldots},
(-13,-7)*{_1},
(-8,-7)*{_2},
(-3,-7)*{_3},
(7,-7)*{_{n-1}},
(13,-7)*{_n},
 (0,0)*{\bullet}="a",
(0,5)*{}="0",
(-12,-5)*{}="b_1",
(-8,-5)*{}="b_2",
(-3,-5)*{}="b_3",
(8,-5)*{}="b_4",
(12,-5)*{}="b_5",
\ar @{.} "a";"0" <0pt>
\ar @{-} "a";"b_2" <0pt>
\ar @{-} "a";"b_3" <0pt>
\ar @{-} "a";"b_1" <0pt>
\ar @{-} "a";"b_4" <0pt>
\ar @{-} "a";"b_5" <0pt>
\endxy
\ \ \ ,
\xy
(1,-5)*{\ldots},
(-13,-7)*{_1},
(-8,-7)*{_2},
(-3,-7)*{_3},
(7,-7)*{_{q-1}},
(13,-7)*{_q},
 (0,0)*{\circ}="a",
(0,5)*{}="0",
(-12,-5)*{}="b_1",
(-8,-5)*{}="b_2",
(-3,-5)*{}="b_3",
(8,-5)*{}="b_4",
(12,-5)*{}="b_5",
\ar @{.} "a";"0" <0pt>
\ar @{.} "a";"b_2" <0pt>
\ar @{.} "a";"b_3" <0pt>
\ar @{.} "a";"b_1" <0pt>
\ar @{.} "a";"b_4" <0pt>
\ar @{.} "a";"b_5" <0pt>
\endxy
\right\rangle_{p,q\geq 2, n\geq 1}
$$
equipped with a differential which is given on white corollas of both colours by
formula (\ref{Ch2: d on white corollas}) and on  black corollas by the following formula
  \Beqr\label{Ch3: d on black corollas}
\p
\xy
(1,-5)*{\ldots},
(-13,-7)*{_1},
(-8,-7)*{_2},
(-3,-7)*{_3},
(7,-7)*{_{n-1}},
(13,-7)*{_n},
 (0,0)*{\bullet}="a",
(0,5)*{}="0",
(-12,-5)*{}="b_1",
(-8,-5)*{}="b_2",
(-3,-5)*{}="b_3",
(8,-5)*{}="b_4",
(12,-5)*{}="b_5",
\ar @{.} "a";"0" <0pt>
\ar @{-} "a";"b_2" <0pt>
\ar @{-} "a";"b_3" <0pt>
\ar @{-} "a";"b_1" <0pt>
\ar @{-} "a";"b_4" <0pt>
\ar @{-} "a";"b_5" <0pt>
\endxy
&=& -\sum_{k=0}^{n-2}\sum_{[n]\setminus [k+1]=I_1\sqcup I_2\atop
 \# I_1\geq 1}\ \ \
\Ba{c}
\begin{xy}
<0mm,0mm>*{\bullet},
<0mm,0.8mm>*{};<0mm,5mm>*{}**@{.},
<-9mm,-5mm>*{\ldots},
<-9mm,-7mm>*{_{1\ \,  \ldots\  \, k}},
<0mm,-10mm>*{...},
<14mm,-5mm>*{\ldots},
<15mm,-7mm>*{\underbrace{\ \ \ \ \ \ \ \ \ \ \ }},
<15mm,-10mm>*{_{I_2}};
<-0.7mm,-0.3mm>*{};<-13mm,-5mm>*{}**@{-},
<-0.6mm,-0.5mm>*{};<-6mm,-5mm>*{}**@{-},
<0.6mm,-0.3mm>*{};<20mm,-5mm>*{}**@{-},
<0.3mm,-0.5mm>*{};<8mm,-5mm>*{}**@{-},
<0mm,-0.5mm>*{};<0mm,-4.3mm>*{}**@{-},
<0mm,-5mm>*{\circ};
<-5mm,-10mm>*{}**@{-},
<-2.7mm,-10mm>*{}**@{-},
<2.7mm,-10mm>*{}**@{-},
<5mm,-10mm>*{}**@{-},
<-6mm,-12mm>*{_{k+1}},
<1.2mm,-12mm>*{\underbrace{\ \ \ \ \ \ \ }},
<2.2mm,-15mm>*{_{I_1}},
\end{xy}
\Ea
\nonumber\\
&& +\ \, \sum_{k=2}^n \sum_{[n]=B_1\sqcup\ldots\sqcup B_k
\atop \inf B_1<\ldots< \inf B_k}
\Ba{c}
\xy
(-15.5,-7)*{...},
(19,-7)*{...},
(7.5,0)*{\ldots},
(-17.8,-12)*{_{B_1}},
(-3.2,-12)*{_{B_2}},
(17.8,-12)*{_{B_k}},
(-1.8,-7)*{...},
%
(-3.2,-9)*{\underbrace{\ \ \ \ \ \  \ \ \ \   }},
%
(-17.8,-9)*{\underbrace{\ \ \ \ \ \  \ \ \ \   }},
%
(16.8,-9)*{\underbrace{\ \ \ \ \ \  \ \ \ \ \  }},
%
 (0,7)*{\circ}="a",
(-14,0)*{\bullet}="b_0",
(-4.5,0)*{\bullet}="b_2",
(14,0)*{\bullet}="b_3",
(0,13)*{}="0",
(1,-7)*{}="c_1",
(-8,-7)*{}="c_2",
(-5,-7)*{}="c_3",
(-22,-7)*{}="d_1",
(-19,-7)*{}="d_2",
(-13,-7)*{}="d_3",
(12,-7)*{}="e_1",
(15,-7)*{}="e_2",
(22,-7)*{}="e_3",
\ar @{.} "a";"0" <0pt>
\ar @{.} "a";"b_0" <0pt>
\ar @{.} "a";"b_2" <0pt>
\ar @{.} "a";"b_3" <0pt>
\ar @{-} "b_2";"c_1" <0pt>
\ar @{-} "b_2";"c_2" <0pt>
\ar @{-} "b_2";"c_3" <0pt>
\ar @{-} "b_0";"d_1" <0pt>
\ar @{-} "b_0";"d_2" <0pt>
\ar @{-} "b_0";"d_3" <0pt>
\ar @{-} "b_3";"e_1" <0pt>
\ar @{-} "b_3";"e_2" <0pt>
\ar @{-} "b_3";"e_3" <0pt>
\endxy
\Ea.
\Eeqr
Representations of this operad in a pair of dg vector spaces, $V_{in}$ and $V_{out}$, is the same
as a triple, $(\mu_{in}, \mu_{out}, F)$, consisting of a $\caL eib_\infty$ structure, $\mu_{in}$, on $V_{in}$,
a $\caL eib_\infty$ structure, $\mu_{out}$, on $V_{out}$, and  of a morphism, $F:(V_{in},\mu_{in})\rar (V_{out},\mu_{out})$,
of $\caL eib_\infty$ algebras.
}

\bip

This first half of the above claim follows from formula (\ref{Ch3:codimension 1 boundary}). The second half
(more precisely, the last three rows) follows from Proposition 2.2.1 and some standard manipulations
 with the bar-cobar constructions
 of the operad of Leibniz algebras. We omit these manipulations so that the reader can interpret the
 sentence in the last
 row
  as a {\em definition}\, of the notion of $\caL ieb_\infty$ {\em morphism}.
The only
    property of such a morphism which we use in this paper is stated and  proven in Appendix 3.

\mip

{\bf 3.4. A semialgebraic atlas on $\widehat{C}_{n,0}$}. Let $z_0=x_0+\ii y_0$ be any point in $\bbH$.
 A {\em placed at $z_0$ configuration}\,
$p\in C_{n,0}$ means a configuration in
$\Conf_{n,0}(\bbH)$ obtained from $z_0$ and $p$ in the following two steps
\Bi
\item[(i)] put $p$ into  its  hyperposition $p^h$ in $C_{n,0}^h$;
\item[(ii)] apply the  translation map (\ref{2: map T_z}),
$$
T_{z_0}: p^h \lon p^h - \ii + z_0,
$$
which moves $z_{min}(p^h)$  from $\ii$ to $z_0$.
\Ei
Note that this operation of {\em placing}\,  a configuration $p\in C_{n,0}$
{\em at a position}\, $z_0\in \bbH$ preserves its Euclidean size, $|p^h-z_{min}(p^h)|$, the property, which is most important in the context of the compactification $\widehat{C}_{n,0}$.\footnote{A suitable analogue of the operation of {\em placing $p\in C_{n,0}$ into position $z_0$}\, in the context of, say, compactification $\widehat{C}_{n,0}'$ will have to use a ``hyperbolic" translation map instead of $T_{z_0}$,
$$
p^h \lon y_0 p^h + x_0,
$$
as it again moves $z_{min}(p^h)$  from $\ii$ to $z_0$ but preserves now the invariant $||p||$ rather than the invariant  $|p^h-z_{min}(p^h)|$.}
\mip

The boundary strata of  $\widehat{C}_{n,0}$ are given by graphs $G\in\cM or(\caL eib_\infty)$ containing
at least one black corolla. A structure of a semialgebraic set in the neighborhood of a
boundary face corresponding to a graph $G$ is best defined in terms of the associated {\em metric}\,
graph, $G_{metric}$,
by the following procedure:
\Bi
\item[(a)] every internal edge of the form
$\xy (0,4)*{\circ}="a", (0,-4)*{\circ}="b",
           \ar @{-} "b";"a" <0pt>\endxy$,
$
\xy (0,4)*{\bullet}="a", (0,-4)*{\circ}="b"
           \ar @{-} "b";"a" <0pt>\endxy$
or
$\xy (0,4)*{\circ}="a", (0,-4)*{\circ}="b",
           \ar @{.} "b";"a" <0pt>\endxy$,
           is assigned a {\em small}\, positive real number $\var \ll +\infty$,

\item[(b)] every white vertex of a dashed corolla  is assigned a {\em large}\, positive real number $\tau \gg 0$,
$$
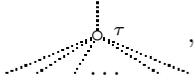
\xy
(3,0)*{^\tau},
(1,-5)*{\ldots},
 (0,0)*{\circ}="a",
(0,5)*{}="0",
(-12,-5)*{}="b_1",
(-8,-5)*{}="b_2",
(-3,-5)*{}="b_3",
(8,-5)*{}="b_4",
(12,-5)*{}="b_5",
\ar @{.} "a";"0" <0pt>
\ar @{.} "a";"b_2" <0pt>
\ar @{.} "a";"b_3" <0pt>
\ar @{.} "a";"b_1" <0pt>
\ar @{.} "a";"b_4" <0pt>
\ar @{.} "a";"b_5" <0pt>
\endxy,
$$

\item[(c)] for every (if any) two vertex subgraph of $G_{metric}$ of the form
$\xy
(3,4)*{^{\tau_1}};
(3,-4)*{^{\tau_2}};
(3,0)*{^{\var}};
(0,4)*{\circ}="a",
(0,-4)*{\circ}="b",
\ar @{.} "b";"a" <0pt>
\endxy$
there is associated a relation, $\tau_2=\var \tau_1$, between the parameters
(which essentially says that $\tau_1\gg \tau_2\gg 0$).

\Ei
Such a metric graph defines a semialgebraic local coordinate chart, $\cU_G$, on $\widehat{C}_{n,0}$ in which
the face $G$ is given by the equations: all  $\var=0$ and all $\tau=+\infty$ (or, better,
all $\tau':=\tanh \tau=1$). The construction of $\cU_G$ should be clear from
a pair of explicit examples one of which we show here and another (illustrating a relation of the type
$\tau_2=\var \tau_1$) in Appendix 2.
\sip

Let $G$ be the face of $\widehat{C}_{n,0}$ corresponding to a graph
$$
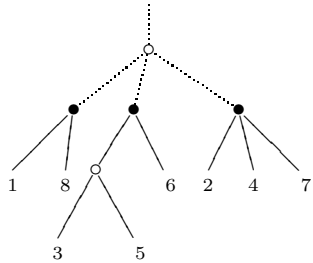
\xy
(-18,-10)*{_1},
(-11,-10)*{_8},
(-12,-19)*{_3},
(-1,-19)*{_5},
(3,-10)*{_6},
(8,-10)*{_2},
(14,-10)*{_4},
(21,-10)*{_7},
(0,14)*{}="0",
 (0,8)*{\circ}="a",
(-10,0)*{\bullet}="b_1",
(-2,0)*{\bullet}="b_2",
(12,0)*{\bullet}="b_3",
(-18,-8)*{}="c_0",
(2,-8)*{}="c_1",
(-7,-8)*{\circ}="c_2",
(8,-8)*{}="c_3",
(14,-8)*{}="c_4",
(20,-8)*{}="c_5",
(-11,-8)*{}="c_8",
(-12,-17)*{}="d_1",
(-2,-17)*{}="d_2",
\ar @{.} "a";"0" <0pt>
\ar @{.} "a";"b_1" <0pt>
\ar @{.} "a";"b_2" <0pt>
\ar @{.} "a";"b_3" <0pt>
\ar @{-} "b_1";"c_0" <0pt>
\ar @{-} "b_1";"c_8" <0pt>
\ar @{-} "b_2";"c_1" <0pt>
\ar @{-} "b_2";"c_2" <0pt>
\ar @{-} "b_3";"c_3" <0pt>
\ar @{-} "b_3";"c_4" <0pt>
\ar @{-} "b_3";"c_5" <0pt>
\ar @{-} "c_2";"d_1" <0pt>
\ar @{-} "c_2";"d_2" <0pt>
\endxy
$$
The associated metric graph is given by
$$
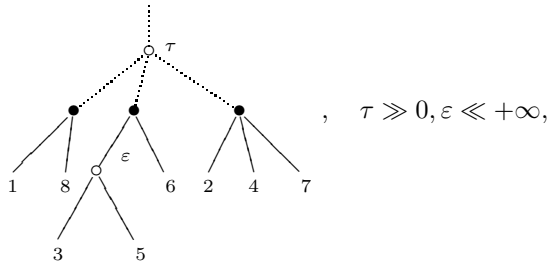
\xy
(-18,-10)*{_1},
(-11,-10)*{_8},
(-12,-19)*{_3},
(-1,-19)*{_5},
(3,-10)*{_6},
(8,-10)*{_2},
(14,-10)*{_4},
(21,-10)*{_7},
%
(-3,-6)*{_{\var}},
(3,8)*{^\tau},
(0,14)*{}="0",
 (0,8)*{\circ}="a",
(-10,0)*{\bullet}="b_1",
(-2,0)*{\bullet}="b_2",
(12,0)*{\bullet}="b_3",
(-18,-8)*{}="c_0",
(2,-8)*{}="c_1",
(-7,-8)*{\circ}="c_2",
(8,-8)*{}="c_3",
(14,-8)*{}="c_4",
(20,-8)*{}="c_5",
(-11,-8)*{}="c_8",
(-12,-17)*{}="d_1",
(-2,-17)*{}="d_2",
\ar @{.} "a";"0" <0pt>
\ar @{.} "a";"b_1" <0pt>
\ar @{.} "a";"b_2" <0pt>
\ar @{.} "a";"b_3" <0pt>
\ar @{-} "b_1";"c_8" <0pt>
\ar @{-} "b_1";"c_0" <0pt>
\ar @{-} "b_2";"c_1" <0pt>
\ar @{-} "b_2";"c_2" <0pt>
\ar @{-} "b_3";"c_3" <0pt>
\ar @{-} "b_3";"c_4" <0pt>
\ar @{-} "b_3";"c_5" <0pt>
\ar @{-} "c_2";"d_1" <0pt>
\ar @{-} "c_2";"d_2" <0pt>
\endxy
\ ,\ \ \  \tau\gg 0, \var\ll +\infty,
$$
and the semialgebraic coordinate chart $U(G)\cap C_{8,0}$ is, by definition,
 an open subset of $C_{8,0}$ consisting of all those configurations, $p$,
  of 8 points in $\bbH$
 which result from the following four step construction:
 \Bi
\item[{\em Step 1}:] take an arbitrary hyperpositioned configuration, $p^{(1)}=(z',z'',z''')\in C_3^h$,  and magnify it (see \S 2.3),
    $$
    p^{(1)}\rar \tau\cdot p^{(1)}=\left(z'_\tau,,z''_\tau, ,z'''_\tau\right):=\left(\tau(z'-\ii) +\ii, \tau(z''-\ii) +\ii ,
    \tau(z'''-\ii) +\ii\right);
    $$
\item[{\em Step 2:}] take arbitrary hyperpositioned configurations of points, $p^{(2)}_1\in C_{2,0}$,  $p^{(2)}_2\in C_{2,0}$ and
$p^{(2)}_3\in C_{3,0}$, labelled, respectively,
 by sets, $\{1,8\}$, $\{z'''', 6\}$ and $\{2,4,7\}$,
and place them at the positions $z'_\tau$, $z''_\tau$ and, respectively, $z'''_\tau$;
\item[{\em Step 3}:] take an arbitrary hyperpositioned configuration, $p^{(4)}\in C_2$, of two points
labelled by $3$ and $5$, $\var$-shrink it as explained in \S 2.3, and
finally place the result at the point $z''''$.
\Ei
The final result is a hyperpositioned point $p=(z_1, z_3,z_5, z_6,z_2,z_4,z_7,z_8)$
in $C_{8,0}$ of the form
$$
 \Ba{c}
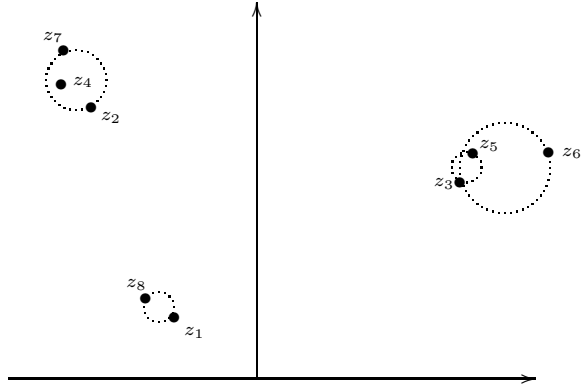
{
\xy
(-23,39)*{\bullet},
(-20,39)*{^{z_4}},
(-22.7,43.5)*{\bullet},
(-24,45)*{^{z_7}},
(30.0,26)*{\bullet}="z_3",
(28,26)*{_{z_3}},
(31.7,29.8)*{\bullet},
(34,31)*{_{z_5}},
(41.8,30)*{\bullet},
(45,30)*{_{z_6}},
(-8,8)*{{\bullet}}="z_1",
(-5.2,6.1)*{_{z_1}},
%
 (30.9,28.1)*{
\xycircle(2,2){.}};
(36.0,28)*{
\xycircle(6,6){.}};
 (-10,9.5)*{
\xycircle(2,2){.}};
(-11.8,10.6)*{\bullet},
(-12.9,12.5)*{_{z_8}},
 (-21,39.7)*{
\xycircle(4,4){.}};
(-19,36)*{\bullet}="z_2",
(-16.3,34.7)*{_{z_2}},
 (-30,0)*{}="a",
(40,0)*{}="b",
(3,0)*{}="c",
(3,50)*{}="d",
\ar @{->} "a";"b" <0pt>
\ar @{->} "c";"d" <0pt>
\endxy}
\Ea
$$
Thus $\cU_G\cap C_{8,0}\simeq  C_3\times C_{2,0}\times C_{2,0}\times C_{3,0}\times C_2\times \R \times \R$.
The boundary strata are given by setting formally $\tau=\infty$ and/or $\var=0$.

\sip
Semialgebraic atlases on the compactifications $\widehat{C}_{n,0}'$ and $\widehat{C}_{n,0}''$
can be described analogously; one has only to take care about suitable analogues of the notion of
{\em placing a configuration at point in $\bbH$}.

\mip

{\bf 3.5. Angle functions
 on $\widehat{C}_{2,0}$.}
 The spaces $\widehat{C}_{2,0}$,
$\widehat{C}_{2,0}'$  and
$\widehat{C}_{2,0}''$ are all identical to each other as semialgebraic sets. They all are given as the closure of the semialgebraic
embedding,
$$
\Ba{rcccc}
C_{2,0} & \lon & C_2^h & \times & [0, 1]\\
p=(z_1,z_2) & \lon & \left\{\Ba{cc}
\left(\ii - \frac{\frac{1}{2}(x_2-x_1)}{\sqrt{\frac{1}{2}(x_2-x_1)^2 + (y_2-y_1)^2}}, \ii+ \frac{\frac{1}{2}(x_2-x_1) +\ii(y_2-y_1)}{\sqrt{\frac{1}{2}(x_2-x_1)^2 + (y_2-y_1)^2}}
\right) & \mbox{if}\  \ y_1\leq y_2 \vspace{2mm}\\
\left(\ii+ \frac{\frac{1}{2}(x_1-x_2) +\ii(y_1-y_2)}{\sqrt{\frac{1}{2}(x_2-x_1)^2 + (y_2-y_1)^2}}, \ii
 - \frac{\frac{1}{2}(x_1-x_2)}{\sqrt{\frac{1}{2}(x_2-x_1)^2 + (y_2-y_1)^2}}
\right) & \mbox{if}\  \ y_1\geq y_2\\
\Ea
\right. && \frac{||p||}{||p||+1}
\Ea
$$
where
$$
||p||=\left\{\Ba{cc}
\frac{\sqrt{\frac{1}{2}(x_2-x_1)^2+(y_2-y_1)^2}}{y_1} & \mbox{if} \ y_1\leq y_2\vspace{2mm}\\
\frac{\sqrt{\frac{1}{2}(x_2-x_1)^2+(y_2-y_1)^2}}{y_2} & \mbox{if} \ y_1\geq y_2\\
\Ea
\right.
$$
 and hence are the unions of two semialgebraic sets,
$$
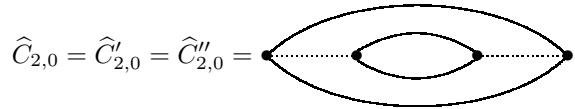
\widehat{C}_{2,0}=\widehat{C}_{2,0}'=\widehat{C}_{2,0}''=
{
\xy
<-20mm,0mm>*{};<-8mm,0mm>*{}**@{.},
<20mm,0mm>*{};<8mm,0mm>*{}**@{.},
(-20,0)*{\bullet},
(20,0)*{\bullet},
(-8,0)*{\bullet},
(8,0)*{\bullet},
(-9,0)*+{};(8.5,0)*+{};
**\crv{(0,6)}
**\crv{(0,-6)}
   \ar@{-}@(ur,lu) (-20.0,0.0)*{};(20.0,0.0)*{},
   \ar@{-}@(dr,ld) (-20.0,0.0)*{};(20.0,0.0)*{},
\endxy}
$$
The inner topological circle, $S^1_{in}$, describes the first boundary component, $C_{1,0}\times C_2^h$
(two points moving very close to each other), in the face decomposition
$$
\p\widehat{C}_{2,0}= C_{1,0}\times C_2\ \bigsqcup \ C_2\times C_{1,0}\times C_{1,0},
$$
while the outer topological circle, $S^1_{out}$,
describes the second boundary component (two points moving very far ---
in the Euclidean or Poincar\'{e} metric --- from each other).

 \mip

{\bf 3.5.1. Definition.} Let $S^1\subset \C$ be given by the equation $z\overline{z}=1$, and let $vol_{S^1}$ stand for
the homogeneous volume form on $S^1$,
$$
vol_{S^1}:= \frac{1}{2i}\left(\overline{z}dz - zd\overline{z}\right)|_{S_1}.
$$
This is a minimal differential 1-form on the semialgebraic manifold $S^1$.
\sip

An {\em angle function}\, on $\widehat{C}_{2,0}$ is a semialgebraic map,
$$
\phi: \widehat{C}_{2,0} \lon S^1\subset \C
$$
such that $\phi^*(vol_{S^1})\mid_{outer\ circle}$ and  $\phi^*(vol_{S^1})\mid_{inner\ circle}$ represent normalized cohomology classes in $H^1(S^1)$, i.e.\
$$
\int_{outer\ circle} \phi^*(vol_{S^1})=2\pi \ \ \ \mbox{and}\ \ \ \ \int_{inner\ circle} \phi^*(vol_{S^1})=2\pi.
$$
The associated  differential form, $\om=\phi^*(vol_{S^1})$, is called a {\em propagator}\, on $\widehat{C}_{2,0}$. It is a
minimal 1-form on the semialgebraic set  $\widehat{C}_{2,0}$.

\mip

{\bf 3.5.2. Example: Kontsevich's propagator.} It is not hard  to check that the function
$$
\Ba{rccc}
\phi_h: & {C}_{2,0} &\lon& S^1\\
& p=(z_1,z_2) & \lon& \phi_h(p):=\sqrt{\frac{(z_1-z_2)(z_1-\overline{z_2})}{(\overline{z_1} -z_2)(\overline{z_1}-\overline{z_2})}}
\Ea
$$
extends to an angle function on $\widehat{C}_{2,0}$
which maps the lower lid of the outer topological circle into a single point.
The associated propagator is denoted by
$$
\om_K:= \phi_h^*(vol_{S^1}).
$$
The geometric meaning of $Arg(\phi_h)$ is explained by the picture,
$$
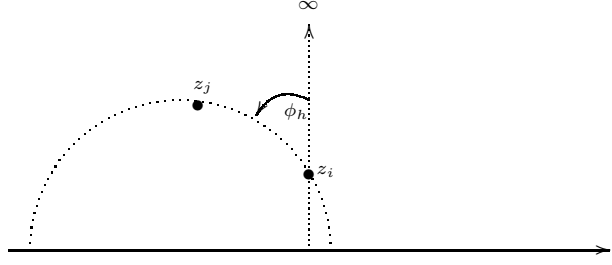
{
\xy
 (-8.56,0)*\ellipse(20,20),=:a(-180){.};
(-14.8,19.24)*{\bullet},
(-14.0,21.4)*{^{z_j}},
(0.0,32.0)*{^{\infty}},
(0,10)*{\bullet},
(0,10)*{\ \ \ \ ^{z_i}},
 (-40,0)*{}="a",
(40,0)*{}="b",
(0,0)*{}="c",
(0,30)*{}="d",
\ar @{->} "a";"b" <0pt>
\ar @{.>} "c";"d" <0pt>
{\ar@/_0.4pc/(0,20)*{};(-7,18)*{}^{{\ \ \phi_h}}}
\endxy}
$$
%
\sip
It measures (in the anticlockwise direction)  the angle between two hyperbolic geodesics passing through $z_i$,
the first one in the vertical direction (towards the point $\infty$) and the second one towards $z_j$.
This propagator was used by Kontsevich in \cite{Ko2} to construct his celebrated formality map.

\mip

{\bf 3.5.3. Example: Kontsevich's antipropagator.} The function
$$
\Ba{rccc}
\hat{\phi}_h: & {C}_{2,0} &\lon& S^1\\
& p=(z_1,z_2) & \lon& \phi_h(p):=\sqrt{\frac{(z_2-z_1)(z_2-\overline{z}_1)}{(\overline{z_2} -z_1)(\overline{z_2}-\overline{z_1})}}
\Ea
$$
extends to an angle function on $\widehat{C}_{2,0}$
which maps the upper lid of the outer topological circle into a single point.
The associated propagator is denoted by
$$
\om_{\overline{K}}:= \hat{\phi}_h^*(vol_{S^1}).
$$

\sip

{\bf 3.5.4. Example: symmetrized Kontsevich's propagator.} The function
$$
\Ba{ccc}
 {C}_{2,0} &\lon& S^1\\
 (z_1,z_2) & \lon& \frac{z_1-z_2}{|z_1-z_2|}=\exp(i \mbox{Arg} (z_1-z_2))
\Ea
$$
extends to an angle function on $\widehat{C}_{2,0}$. The associated propagator is denoted by $\om^{sym}$
 because of the equality
$$
 \om^{sym}= \frac{1}{2}\left(\om_K(z_1,z_2) + \om_K(z_2,z_1)\right)=d  \mbox{Arg} (z_1-z_2).
$$

\mip

{\bf 3.5.5. Example: a one-parameter family of propagators.} For any $t\in [0,1]$ we set
 $$
 \om^{t}_{{K}}= t\om_K + (1-t)\om_{\overline{K}}.
 $$
It i easy to check that
$$
\om_K^t(z_1,z_2)|_{inner\ circle}=dArg(z_1-z_2), \ \ \ \forall t\in [0,1].
$$

\mip

{\bf 3.6. Renormalized forgetful map.} Note that for any pair of integers $i,j\in [n]$, $i\neq j$,
there is an associated forgetting map,
\Beq\label{Ch3: p-forgetting map}
\Ba{rccc}
p_{ij}: & C_{n,0}  &\lon & C_{2,0} \\
        & (z_1, \ldots, z_n) & \lon & (z_i, z_j),
\Ea
\Eeq
which extends to a semialgebraic map of their compactifications,
$$
{p}_{ij}:\widehat{C}_{n,0}\rar
\widehat{C}_{2,0}.
$$
 Hence, for any propagator $\om$ on $\widehat{C}_{2,0}$ the pull-back ${p}_{ij}^*(\om)$ is a well-defined
 minimal one-form on  $\widehat{C}_{n,0}$. 

\sip

We shall need below a ``renormalized"  version of the forgetting map,
\Beq\label{Ch3: fp-forgetting map}
\Ba{rccc}
\fp_{ij}: & \widehat{C}_{n,0} &\lon & \widehat{C}_{2,0}\\
& (z_1, \ldots, z_i, \ldots, z_j, \ldots, z_n) & \lon & \displaystyle\left\{  \Ba{cc}
\left(z_i- z_j+ z_{min}(p), z_{min}(p)\right) & \mbox{if} \ y_i\geq y_j\\
\left(z_{min}(p),z_j-z_i+  z_{min}(p)\right)  & \mbox{if} \ y_i\leq y_j
\Ea
\right.
\Ea
\Eeq
Note that for $n=2$ we have $\fp_{12}=\Id$.

\bip

\section{De Rham field theories on  configuration spaces}
{\bf 4.1. Families of graphs.}
Let ${\cG}_{n,l}$ stand for a family of graphs, $\{\Ga\}$, with $n$ vertices
and $l$ edges such that
\Bi
\item the edges of $\Ga$ are directed, beginning and ending at {\em different}\, vertices;
\item the set of  vertices, $V(\Ga)$, is labeled by the set $[n]$;
\item the set of edges, $E(\Ga)$, is totally ordered.
\Ei
We identify two total orderings on the set $E(\Ga)$ (that is, isomorphisms $E(\Ga)\simeq [\# E(\Ga)]$),
if they differ by an even permutation of $[\# E(\Ga)]$. Thus there are precisely two possible orderings\footnote{It is useful sometimes  to identify an orientation of $\Ga\in  {\cG}_{n,l}$
with a vector $\displaystyle \f_\Ga:=\wedge_{e\in E(\Ga)}e$ in the real one dimensional vector space $\wedge^{l}\R[E(\Ga)]$,
 where $\R[E(\Ga)]$ is the $l$-dimensional vector space spanned over $\R$ by the set $E(\Ga)$.}
on the set $E(\Ga)$
 and the group $\Z_2$ acts freely  on ${\cG}_{n,l}$ by ordering changes; its orbit
space, $(\Ga, \Ga_{opp})$, is denoted by $\fG_{n,l}$.

\sip
With every graph $\Ga\in {\cG}_{n,l}$ one can associate a linear map,
$$
\Ba{rccc}
\Phi_\Ga: &\ot^n \cT_{poly}(\R^d) &\lon &\cT_{poly}(\R^d)[-l]\\
& \ga_1\ot\ldots \ot\ga_n & \lon& \Phi_\Ga(\ga_1,\ldots,\ga_n)
\Ea
$$
where
\Beq\label{Ch4: Phi_Ga}
\Phi_\Ga(\ga_1,...,\ga_n)=
\left[\left(\prod_{e\in E(\Ga)}
\Delta_{e}\right)
\ga_1(\psi_{(1)},x_{(1)})\cdots \ga_n(\psi_{(n)},x_{(n)})
\right]_{\hspace{-1mm}x_{(1)}=\ldots=x_{(n)}\atop
\psi_{(1)}=\ldots=\psi_{(n)}}
\Eeq
and, for an edge $e$ beginning at a vertex labelled by $i$ and ending at a vertex labelled by $j$,
$$
\Delta_e:=
\sum_{a=1}^d 
\frac{\p^2}{\p x_{(j)}^a\p \psi_{(i)\,a}}.
$$

{\bf 4.1.1. Complete subgraphs}. For any subset $A\subset [n]$ and any graph $\Ga$ in $\cG_{n,l}$
 (or in $\fG_{n,l}$), there is an associated complete subgraph $\Ga_A$ of $\Ga$ whose vertices
are, by definition, those vertices of $\Ga$ which are labelled by elements of $A$, and whose edges
are all the edges of $\Ga$ which connect these $A$-labelled vertices. If we shrink all the
$A$-labelled vertices of $\Ga$ (together with all the edges connecting these $A$-labelled vertices) into a single vertex,
then we obtain from $\Ga$ a new graph which we denote by $\Ga/\Ga_A$.

\sip

Similarly, for any family of disjoint subsets $A_1,\ldots, A_k$ of $[n]$
and any graph  $\Ga$ in $\cG_{n,l}$ (or in  $\fG_{n,l}$) one can define complete subgraphs
$\Ga_{A_1},\ldots,\Ga_{A_k}\subset \Ga$ as well as the quotient graph $\Ga/\{\Ga_{A_1},\ldots,\Ga_{A_k}\}$.

\sip


{\bf 4.1.2. Lemma}. {\em Let $A$ be a (naturally ordered) proper subset of $[n]$,
$\Ga_1\in {\cG}_{n-\#A +1, l_1}$ and $\Ga_2\in {\cG}_{\#A ,l_2}$. Then
$$
\Phi_{\Ga_1}(\ga_1, \ldots, \ga_{\inf A-1}, \Phi_{\Ga_2}(\ga_A),
\ga_{[n-\inf A+1]\setminus A})=\sum_{\Ga\in {\cG}_{n,l_1+l_2}(A,\Ga_{1,2})} (-1)^{\sigma_\Ga}
\Phi_\Ga(\ga_1, \ldots, \ga_n),
$$
where $(-1)^{\sigma_\Ga}$ is the standard Koszul sign and ${\cG}_{n,l_1+l_2}(A,\Ga_{1,2})$ is a subset of
${\cG}_{n,l_1+l_2}$ consisting of all those graphs $\Ga$
whose complete subgraph $\Ga_A$ is isomorphic to $\Ga_2$ and the quotient graph $\Ga/\Ga_A$ is isomorphic
to $\Ga_1$ (and which are equipped with the natural ordering of edges induced from the ones
in $\Ga_1$ and $\Ga_2$).
}

\mip

Proof uses only definition (\ref{Ch4: Phi_Ga})  and the Leibniz rule for partial derivatives.

\mip

{\bf 4.1.3. Definitions}. (i) Given a graph  $\Ga\in \fG_{n,2n-4}$. A subset
$A\subset Vert(\Ga)\simeq[n]$ is called {\em admissible}\, if $2\leq \# A\leq n-1$ and  the associated subgraph $\Ga_A$ belongs to
$\fG_{\# A, 2\# A -3}$. Note that in this case
$\Ga/\Ga_A\in \fG_{n-\# A+1, 2(n-\# A +1) -3}$.

\sip

 (ii) Given a graph  $\Ga\in \fG_{n,2n-3}$. A subset
$A\subset Vert(\Ga)\simeq[n]$ is called {\em admissible}\, if the associated subgraph $\Ga_A$ belongs to
$\fG_{\# A, 2\# A -3}$. Note that in this case
$\Ga/\Ga_A\in \fG_{n-\# A+1, 2(n-\# A +1) -2}$.

\sip

(iii)
Given a graph  $\Ga\in \fG_{n,2n-3}$. A decomposition
$Vert(\Ga)=A_1\sqcup\ldots \sqcup A_k$, $k\geq 2$,  is called admissible if
$\Ga_{A_i}\in \fG_{\# A_i, 2\# A_i -2}$, $i=1,\ldots,k$.
In this case $\Ga/\{\Ga_{A_1},...,\Ga_{A_k}\}\in \fG_{k, 2k-3}$.

\sip

Analogous definitions can be made for graphs from the sets  ${\cG}_{n,2n-4}$ and
${\cG}_{n,2n-3}$.

\mip

{\bf 4.2. De Rham field theory on $\overline{C}$}. For any proper subset $A\subset [n]$ of cardinality
at least two there is an associated embedding,
$$
i_A: \overline{C}_{n-\# A +1}\times  \overline{C}_{\# A} \lon \overline{C}_n,
$$
 of the  corresponding  boundary component into $\overline{C}_n$. For example, for $A=\{3,5,6,7\} \subset [7]$ the image of $i_A$ is a face
of $\overline{C}_8$ represented by the following graph
$$
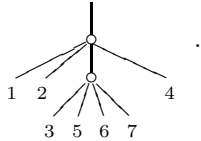
\begin{xy}
<0mm,0mm>*{\circ},
<0mm,0.8mm>*{};<0mm,5mm>*{}**@{-},
<-10.5mm,-7mm>*{_{1}},
<-6.5mm,-7mm>*{_{2}},
<10.5mm,-7mm>*{_{4}};
<-0.7mm,-0.3mm>*{};<-10mm,-5mm>*{}**@{-},
<-0.6mm,-0.5mm>*{};<-6mm,-5mm>*{}**@{-},
<0.3mm,-0.5mm>*{};<10mm,-5mm>*{}**@{-},
<0mm,-0.5mm>*{};<0mm,-4.3mm>*{}**@{-},
<0mm,-5mm>*{\circ};
<-5mm,-10mm>*{}**@{-},
<-1.7mm,-10mm>*{}**@{-},
<1.7mm,-10mm>*{}**@{-},
<5mm,-10mm>*{}**@{-},
<-5.5mm,-12mm>*{_{3}},
<-1.8mm,-12mm>*{_{5}},
<1.8mm,-12mm>*{_{6}},
<5.5mm,-12mm>*{_{7}},
\end{xy}\ \ .
$$
Let $\Omega^\bullet(\overline{C}_n)=\oplus_{p\geq 0} \Omega^p(\overline{C}_n)$ stand for the de Rham algebra of $PA$ differential forms on the space  $\overline{C}_n$.
 A {\em de Rham field theory}\, on the family of compactified
 configuration spaces
 $\{\overline{C}_n\}_{n\geq 2}$ is, by definition, a collections of maps,
$$
\left\{
\Ba{rccc}
\Omega: & {\cG}_{n, l}& \lon & \Omega^{l}(\overline{C}_n)\\
& \Ga & \lon & \Omega_\Ga
\Ea
\right\}_{n\geq2, l\geq 1},
$$
such that $d\Omega_\Ga=0$,  $\Omega_{\Ga_{opp}}=-\Omega_\Ga$,
and, for any $\Ga\in  {\cG}_{n, 2n-4}$
and any proper subset $A\subset Vert(\Ga)$ with $\# A\geq 2$,  one has
\Beq\label{Ch4: de  Rham cond on C_n}
i_A^*(\Omega_\Ga)\simeq (-1)^{\sigma_A}\Omega_{\Ga/\Ga_A} \wedge \Omega_{\Ga_A}.
\Eeq
where the sign $(-1)^{\sigma_A}$ is defined by the equality
$
\f_\Ga=(-1)^{\sigma_A} \f_{\Ga/\Ga_A}\wedge \f_{\Ga_A}$
(see footnote 7). The symbol $\simeq $ means here and below
equality modulo differential forms
whose integral over the corresponding boundary component is zero.

\mip
{\bf 4.2.1. Theorem.}
{\em Given a de Rham field theory on $\overline{C}$, then, for any $d\in \N$,
there is an associated $\caL eib_\infty$-algebra structure,
$$
\Ba{rccc}
\mu_n: & \ot^n\cT_{poly}(\R^d) &\lon &  \cT_{poly}(\R^d)[3-2n]\\
& \ga_1\ot \ldots\ot \ga_n     & \lon & \mu_n(\ga_1, \ldots, \ga_n)
\Ea,
$$
 on $\cT_{poly}(\R^d)$ given by}
\Beq\label{Ch4: induced_Leib_infty}
\mu_n(\ga_1, \ldots, \ga_n):=
\left\{
\Ba{cl}
0 & \mbox{for}\ n=1,\\
 \sum_{\Ga\in \fG_{n,2n-3}} c_{\Ga}\Phi_\Ga(\ga_1, \ldots, \ga_n)
 & \mbox{for}\ n\geq 2
 \Ea
 \right.
\Eeq
with\footnote{$c_\Ga$ and $\Phi_\Ga$ are computed for an arbitrary lift of $\Ga\in \fG_{n,2n-2}$ to an element
of ${\cG}_{n,2n-2}$; the product $c_\Ga\Phi_\Ga$ is independent of the choice of such a lift.}
\Beq\label{Ch4: c_gamma_on_C_n}
c_\Ga:= \int_{\overline{C}_n} \Omega_\Ga.
\Eeq

\mip

\Proof By the Stokes theorem, for any $\Ga\in \fG_{n, 2n-4}$,
$$
0=\int_{\overline{C}_n}d\Omega_\Ga =\int_{\p\overline{C}_n}\Omega_\Ga=
\sum_{A\varsubsetneq [n]\atop \# A\geq 2}(-1)^{\sigma_A}
\int_{\overline{C}_{n-\# A+1}}\Omega_{\Ga/\Ga_A}
\int_{\overline{C}_{\# A}}\Omega_{\Ga_A}  =
\sum_{A\subset V(\Ga)\atop
A\, \mathit i \mathit s\, \mathit a \mathit d \mathit m \mathit i \mathit s \mathit s \mathit i \mathit b \mathit l \mathit e}(-1)^{\sigma_A}c_{\Ga_A}c_{\Ga/\Ga_A}.
$$
Then, using Lemma~{4.1.2}, one obtains  
$$
\sum_{A\subsetneqq [n]\atop
\# A\geq 2}(-1)^{\sum_{k=1}^{\inf A-1}|\ga_k|}
\mu_{n-\# A +1}(\ga_1, \ldots, \ga_{\inf A-1}, \mu_{\# A}(\ga_A),\ga_{[n-\inf A]\setminus A})=
 \hspace{70mm}
$$
\Beqrn
 &=&
\sum_{A\subsetneqq [n]\atop
\# A\geq 2}
\sum_{\Ga_1\in \fG_{N, 2N-3}\atop N:=n-\# A +1}\sum_{\Ga_2\in \fG_{\# A, 2\# A-3}}\hspace{-3mm}
(-1)^{\sum_{k=1}^{\inf A-1}|\ga_k|}
c_{\Ga_1}c_{\Ga_2}\Phi_{\Ga_1}(\ga_1, \ldots, \ga_{\inf A-1}, \Phi_{\Ga_2}(\ga_A),
\ga_{[n-\inf A]\setminus A})\\
 &=&
\sum_{\Ga\in \fG_{n,2n-4}} 
\left(
\sum_{A\subset Vert(\Ga)\atop
A\, \mathit i \mathit s\, \mathit a \mathit d \mathit m \mathit i \mathit s \mathit s \mathit i \mathit b \mathit l \mathit e}(-1)^{\sigma_A}c_{\Ga_A}c_{\Ga/\Ga_A}\right)\Phi_{\Ga}(\ga_1,\ldots,\ga_n)\\
  &=& 0,
\Eeqrn
which proves the claim. \hfill $\Box$

\bip

For any pair of different integers, $i, j\in[n]$,
there is an associated  map,
$$
\Ba{rccc}
\pi_{ij}: & C_n & \lon & C_2=S^1\\
& (z_1,\ldots, z_i, \ldots, z_j,\ldots, z_n) & \lon & \frac{z_i - z_j}{|z_i-z_j|},
\Ea
$$
which extends to the compactifications, ${\pi}_{ij}:\overline{C}_n \rar
\overline{C}_2$.

\bip

{\bf 4.2.2. Proposition.} {\em For any PA 1-form, $\om$, on $C_2\simeq S^1$ satisfying the condition
\Beq\label{Ch4: integral om_0 over C_2}
\int_{S^1}\om=2\pi,
\Eeq
the associated map
$$
\Ba{rccc}
\Omega: & {\cG}_{n, l}& \lon & \Omega^{l}(\overline{C}_n)\\
& \Ga & \lon & \displaystyle \Omega_\Ga:=\bigwedge_{e\in E(\Ga)}\frac{{\pi}_e^*(\om)}{2\pi}
\Ea
$$
with
$
{\pi}_e:={\pi}_{ij}$ for an edge
$e$
beginning at a vertex labelled by $i\in [n]$ and ending at a vertex labelled by $j\in [n]$,
 defines a non-trivial
de Rham field theory on $\overline{C}$ (and hence an associated non-trivial $\caL ie_\infty$-structure
 on $\cT_{poly}(\R^d)$).}

\begin{proof}
Basic condition (\ref{Ch4: de  Rham cond on C_n}) can be easily checked in the coordinate chart near the stratum
$\Img i_A$ (see \S 2.2). The weights, $c_\Ga=\int_{C_n}\Omega_\Ga$, are
independent of the labeling maps, $V(\Ga)\rar [n]$, so that the resulting
$\caL eib_\infty$ operations $\mu_n$
are graded symmetric and give, therefore, a $\caL ie_\infty$-structure on
 $\cT_{poly}(\R^d)$. Its non-triviality follows from a particular example studied
 in the next subsection.
\end{proof}

\mip

{\bf 4.2.3. Example: Schouten brackets from the homogeneous volume form on $S^1$.} The first natural choice for a
differential 1-form on
$C_2$ satisfying  condition (\ref{Ch4: integral om_0 over C_2}) is, of course, the following one
\Beq\label{Ch4: standard volume form on S^1}
\om_0(z_1,z_2)=dArg(z_1-z_2).
\Eeq
 By Kontsevich's ``vanishing" Lemma~6.4 in \cite{Ko}, the associated weights $c_\Ga$ are zero
 for all graphs $\Ga\in \sqcup_{n\geq 2}\cG_{n,2n-3}$ excepts
 for
$
\Ga_1=\xy
<0mm,4.5mm>*{^{1}},
<0mm,-5mm>*{_{2}},
(0,3)*{\bullet}="a",
(0,-3)*{\bullet}="b",
\ar @{->} "a";"b" <0pt>
\endxy
$ and
$
\Ga_2=\xy
<0mm,4.5mm>*{^{1}},
<0mm,-5mm>*{_{2}},
(0,3)*{\bullet}="a",
(0,-3)*{\bullet}="b",
\ar @{<-} "a";"b" <0pt>
\endxy
$\
 which both have weight $1$. Hence all  $\caL ie_\infty$-operations
(\ref{Ch4: induced_Leib_infty}) are zero for $n\geq 3$, and
$$
\mu_2(\ga_1, \ga_2)=c_{\Ga_1}\Phi_{\Ga_1}(\ga_1,\ga_2) + c_{\Ga_2}\Phi_{\Ga_2}(\ga_1,\ga_2)=
(-1)^{|\ga_1|}[\ga_1\bullet \ga_2].
$$
Thus the propagator  $dArg(z_1-z_2)$ on $C_2$ is responsible for the
existence in nature of the Schouten bracket on polyvector fields.

\mip

{\bf 4.2.4. Example: a $\caL ie_\infty$-structure on $\cT_{poly}(\R^d)$ from the Kontsevich propagator}. Consider the following minimal 1-form,
$$
\om:= \om_K|_{S^1_{out}}
$$
on $C_2$. It satisfies the normalization condition (\ref{Ch4: integral om_0 over C_2}) and hence defines by Proposition 4.2.2 a non-trivial $\caL ie_\infty$-structure on $\cT_{poly}(\R^d)$ for any $d$. It is clear that the weights
$$
c_\Ga= \int_{\overline{C}_n} \bigwedge_{e\in E(\Ga)}\frac{{\pi}_e^*(\om)}{2\pi}
$$
can be  non-zero only for graphs $\Ga$ with no closed paths of oriented edges, and if $\Ga\in \fG_{n,2n-3}$ is such a graph then
$$
c_\Ga= \int_{\overline{C_n(\Ga)}} \bigwedge_{e\in E(\Ga)}\frac{{\pi}_e^*\left(d Arg (z_1-z_2)\right)}{2\pi}
$$
where $\overline{C}_n(\Ga)$ is is the closure in $\overline{C}_n$ of a
semialgebraic subset ${C}_n(\Ga) \subset C_n$ consisting of all configurations, $(z_1,\ldots, z_n)$,
 satisfying the following family of inequalities (one for each directed edge in $\Ga$),
 $$
 \Img z_i \leq \Img z_j \ \mbox{if there is an edge $e$ beginning at $j$ and ending at $i$}.
 $$
Thus the $\caL ie_\infty$ structure induced on $\cT_{poly}(\R^d)$ by Proposition~4.2.2 is precisely the one introduced by Shoikhet in \cite{Sh2}. It is clear that $\mu_2$ operation is equal to
the Schouten bracket again. We prove in Appendix 4 that $\mu_4= \sum_{\Ga\in \fG_{4,5}}c_\Ga \Phi_\Ga$ is spanned by following three graphs,
\Beq\label{Ch4: three graphs}
\Ga_1=  \underset{_{31\wedge 32\wedge 41\wedge 42\wedge21}}{\Ba{c}\xy
<0mm,-2.0mm>*{_{1}},
<0mm,10.0mm>*{^{2}},
<-7mm,17.5mm>*{^{3}},
<7mm,17.5mm>*{^{4}},
(0,0)*{\bullet}="1",
(-7,16)*{\bullet}="2",
(7,16)*{\bullet}="3",
(0,8)*{\bullet}="4",
\ar @{->} "2";"4" <0pt>
\ar @{->} "3";"4" <0pt>
\ar @{->} "4";"1" <0pt>
\ar @{->} "2";"1" <0pt>
\ar @{->} "3";"1" <0pt>
\endxy\Ea}, \ \ \
\Ga_2=  \underset{_{42\wedge 43\wedge 31\wedge 21\wedge 32}}{\Ba{c}\xy
<0mm,-2.0mm>*{_{1}},
<-6.5mm,7.5mm>*{^{2}},
<6.5mm,11.5mm>*{^{3}},
<0mm,17.5mm>*{^{4}},
(0,0)*{\bullet}="1",
(-6,6)*{\bullet}="2",
(6,10)*{\bullet}="3",
(0,16)*{\bullet}="4",
\ar @{->} "4";"3" <0pt>
\ar @{->} "4";"2" <0pt>
\ar @{->} "3";"2" <0pt>
\ar @{->} "2";"1" <0pt>
\ar @{->} "3";"1" <0pt>
\endxy\Ea}, \ \ \
\Ga_3=   \underset{_{41\wedge 31\wedge 32\wedge 32\wedge 43}}{ \Ba{c}\xy
<-7mm,-2.0mm>*{_{1}},
<7mm,-2.0mm>*{^{2}},
<0mm,4.0mm>*{^{3}},
<0mm,17.5mm>*{^{4}},
(0,16)*{\bullet}="1",
(-7,0)*{\bullet}="2",
(7,0)*{\bullet}="3",
(0,8)*{\bullet}="4",
\ar @{<-} "2";"4" <0pt>
\ar @{<-} "3";"4" <0pt>
\ar @{<-} "4";"1" <0pt>
\ar @{<-} "2";"1" <0pt>
\ar @{<-} "3";"1" <0pt>
\endxy\Ea}
\Eeq
which contribute into $\mu_4$ --- for any fixed  numbering of vertices, say the one shown in the pictures above, and for the ordering of edges
shown under the pictures ---  with the one and the same weight
$$
c_{\Ga_1}=c_{\Ga_2}=c_{\Ga_3}=\frac{1}{12}.
$$
 One can show that Maurer-Cartan elements of this $\caL ie_\infty$ can be quantized by standard iteration, i.e.\ their deformation quantization does not require sophisticated mathematics, and all the subtleties are, therefore, hidden in the $\caL ie_\infty$ quasi-isomorphism between the standard Schouten
algebra and  this exotic $\caL ie_\infty$ algebra. We construct explicitly such a quasi-isomorphism below
using the Kontsevich propagator and the renormalized forgetful map (\ref{Ch3: fp-forgetting map}).
\sip

{\bf 4.2.5. On homotopy equivalence of $\caL ie_\infty$ algebras}. Let $\R[t,dt]$ stand for the polynomial de Rham algebra
 on $\R$ and $d$ for the de Rham differential. The tensor product $\cT_{poly}(\R^d)[t,dt]:=\cT_{poly}(\R^d)  \ot_\R \R[[t,dt]]$ is naturally a dg module
over  $\R[t,dt]$. A $\caL ie_\infty$-structure, $\mu_\bu(t,dt)$, on  $\cT_{poly}(\R^d)[t,dt]$ such that
$$
\mu_1(t,dt)=d
$$
and all the higher operations, $\mu_n(t,dt)$, are morphisms of  $\R[t,dt]$-modules
is uniquely determined by
two families of operations,
$$
\left\{\Ba{ccc}
\mu'_n(t): \odot^n\cT_{poly}(\R^d) &\lon&  \cT_{poly}(\R^d)[3-2n]\ot \R[t], \\
\mu''_n(t): \odot^n\cT_{poly}(\R^d) &\lon &  \cT_{poly}(\R^d)[2-2n]\ot \R[t]
\Ea
\right\}_{n\geq 2}
$$
such that $\mu_n(t,dt)= \mu'_n(t) + dt \mu''_n(t)$ for $n\geq 2$.
We call such a  $\caL ie_\infty$-structure on  $\cT_{poly}(\R^d)[t,dt]$ a {\em path}\, one.
Two minimal $\caL ie_\infty$-structures, say $\mu_\bu$ and $\hat{\mu}_\bu$, on $\cT_{poly}(\R^d)$ are  called
{\em gauge or homotopy equivalent}\,  \cite{Fu} if there exists a path  $\caL ie_\infty$-structure $\mu_\bu(t,dt)$
such that $\mu'_\bu(t)|_{t=0}= \mu_\bu$ and $\mu'_\bu(t)|_{t=1}= \hat{\mu}_\bu$.

\sip
Any PA 1-form on the semialgebraic set $C_2$ can be written in the form,
$\om=\om_{0} + df$,
where $\om_{0}=dArg(z_1-z_2)$ and $f$ is some $PA$ function on $C_2$.

\sip

{\bf 4.2.6. Proposition}. {\em For any  PA function $f$ on $C_2$, the operators
$\mu'_n(t):=\sum_{\Ga\in \fG_{n,2n-3}} c'_{\Ga}(t)\Phi_\Ga$ and  $\mu''_n(t):=\sum_{\Ga\in \fG_{n,2n-2}}  c''_{\Ga}(t)\Phi_\Ga$
with
\Beq\label{Ch4: c' c''_gamma_on_C_n}
c'_\Ga(t):= \int_{\overline{C}_n}
\bigwedge_{e\in E(\Ga)}\frac{{\pi}_{e'}^*(\om_{00} + tdf)}{2\pi}, \ \ \
\mbox{and}\ \ \ dt c''_\Ga(t):= \int_{\overline{C}_n}
\bigwedge_{e\in E(\Ga)}\frac{{\pi}_{e}^*(\om_0 +tdf+ fdt)}{2\pi}.
\Eeq
define a path  $\caL ie_\infty$-structure,
\Beq\label{Ch4: mu(t,dt)^f}
\mu^f_n(t,dt) =  \mu'_n(t) + dt \mu''_n(t), \ \ \ n\geq 2,
\Eeq
  on   $\cT_{poly}(\R^d)[t,dt]$
such that $\mu'_\bu(t)|_{t=0}$ equals the standard Schouten algebra structure on  $\cT_{poly}(\R^d)$
and  $\mu'_\bu(t)|_{t=1}$ equals the $\caL ie_\infty$-structure on $\cT_{poly}(\R^d)$ associated by Proposition
4.2.2 to the generic propagator $\om=\om_{0} + df$}.

\begin{proof} First we note that the definition of the weight $c''_\Ga(t)$ makes sense. In fact, one has,
$$
c''_\Ga(t):= \sum_{e\in E(\Ga)}(-1)^{|e|}\int_{\overline{C}_n} \frac{\pi_e^*(f)}{2\pi}
\bigwedge_{e'\in E(\Ga)\atop e'\neq e}\frac{{\pi}_{e'}^*(\om_0 +tdf)}{2\pi}
$$
where $|e|$ counts the number of edges of $\Gamma$ staying before the edge $e$ in the chosen total ordering,
$o:   E(\Ga)\rar [\# E(\Ga)]$, of edges, i.e.\ $|e|:= o(e)-1$.
\sip

It is obvious that the required conditions  on the boundary values $\mu'_\bu(t)|_{t=0}$ and  $\mu_\bu(t)_{t=1}$
are satisfied.

\sip

Next it is a straightforward calculation (which is fully analogous to the one made on the proof of Theorem~4.2.1)
to check that  maps (\ref{Ch4: mu(t,dt)^f}) define a $\caL ie_\infty$-algebra structure on  $\cT_{poly}(\R^d)\ot \K[[t,dt]]$
if and only if one has, for any  $\Ga\in \cG_{n,2n-4}$,
$$
\sum_{A\subset V(\Ga)\atop
A\, \mathit i \mathit s\, \mathit a \mathit d \mathit m \mathit i \mathit s \mathit s \mathit i \mathit b \mathit l \mathit e}(-1)^{\sigma_A}c'_{\Ga_A}(t)c'_{\Ga/\Ga_A}(t)=0
$$
 and, for any $\Ga\in \cG_{n,2n-3}$,
$$
\frac{dc'_\Ga(t)}{dt}= \sum_{A\subset V(\Ga)\atop
{such\ that \atop \Ga_A\in \cG_{\# A, 2\# A-3}  }}(-1)^{\sigma_A}c'_{\Ga_A}(t)c''_{\Ga/\Ga_A}(t) +
 \sum_{A\subset V(\Ga)\atop
{such\ that \atop \Ga_A\in \cG_{\# A, 2\# A-2}  }}(-1)^{\sigma_A}c''_{\Ga_A}(t)c'_{\Ga/\Ga_A}(t).
$$
The first condition is obvious (see again the proof of Theorem 4.2.1 above) while the second one follows from the following
calculation,
\Beqrn
\frac{dc'_\Ga(t)}{dt}&=& \sum_{e\in E(\Ga)}(-1)^{|e|}\int_{\overline{C}_n} \frac{\pi_e^*(df)}{2\pi}
\bigwedge_{e'\in E(\Ga)\atop e \neq e}\frac{\overline{\pi}_{e'}^*(\om_0 +tdf)}{2\pi}\\
&=&
\sum_{e\in E(\Ga)}(-1)^{|e|}\int_{\p\overline{C}_n} \frac{\pi_e^*(f)}{2\pi}
\bigwedge_{e'\in E(\Ga)\atop e \neq e}\frac{\overline{\pi}_{e'}^*(\om_0 +tdf)}{2\pi}\\
&=&\hspace{-2mm}
\sum_{e\in E(\Ga)}\hspace{-1mm}(-1)^{|e|}\hspace{-1mm}
\left(\sum_{A\subset V(\Ga)\atop e\neq E(\Ga_A)}(-1)^{\sigma_A}\hspace{-1mm}
\underset{\overline{C}_A}{\int} \hspace{-1mm}
\bigwedge_{e'\in E(\Ga_A)}\hspace{-4mm}\frac{\overline{\pi}_{e'}^*(\om_0 +tdf)}{2\pi}\underset{{C}_{n-\# A+ 1}}{\int}
\hspace{-3mm}
\frac{\pi_e^*(f)}{2\pi}\hspace{-3mm}
\bigwedge_{e'\in E(\Ga/\Ga_A)\atop e \neq e}\hspace{-4mm} \frac{\overline{\pi}_{e'}^*(\om_0 +tdf)}{2\pi}
\right.\vspace{0mm}\\
&&\hspace{-2mm}
+\ \left.
\sum_{A\subset V(\Ga)\atop e\in E(\Ga_A)}(-1)^{\sigma_A}\hspace{-1mm}
\underset{\overline{C}_A}{\int} \hspace{-1mm}\frac{\pi_e^*(f)}{2\pi}\hspace{-3mm}
\bigwedge_{e'\in E(\Ga_A)\atop e \neq e}\hspace{-4mm} \frac{\overline{\pi}_{e'}^*(\om_0 +tdf)}{2\pi}
\underset{{C}_{n-\# A+ 1}}{\int}
\hspace{0mm}
\bigwedge_{e'\in E(\Ga/\Ga_A)}\hspace{-4mm}\frac{\overline{\pi}_{e'}^*(\om_0 +tdf)}{2\pi}
\right)\\
&=&\sum_{A\subset V(\Ga)\atop
{such\ that \atop \Ga_A\in \cG_{\# A, 2\# A-3}  }}(-1)^{\sigma_A}c'_{\Ga_A}(t)c''_{\Ga/\Ga_A}(t) +
 \sum_{A\subset V(\Ga)\atop
{such\ that \atop \Ga_A\in \cG_{\# A, 2\# A-2}  }}(-1)^{\sigma_A}c''_{\Ga_A}(t)c'_{\Ga/\Ga_A}(t).
\Eeqrn
\end{proof}

{\bf 4.2.7. De Rham field theory of path  $\caL ie_\infty$-structures}. One can generalize the notion of a
de Rham field theory on $\overline{C}$ to a map
$$
\left\{
\Omega(t,dt):  {\cG}_{n, l} \lon  \Omega^{l}(\overline{C}_n\times \R)
\right\}_{n\geq2, l\geq 1},
$$
taking values in $\Omega^\bullet(\overline{C}_n\times \R):= \Omega^\bullet(\overline{C}_n)[t,dt]$
and satisfying the same closeness  and factorization condition (\ref{Ch4: de  Rham cond on C_n}) as the map $\Omega$
in the beginning of \S 4.2. An analogue of Theorem~4.2.1 claiming that every such a de Rham field theory
$\Omega(t,dt)$
 defines a {\em path}\, $\caL eib_\infty$-algebra structure on $\cT_{poly}(\R^d)[t,dt]$ holds true with the only
difference that the summation in formula (\ref{Ch4: induced_Leib_infty}) goes over graphs
$\Ga\in \fG_{n,2n-3}\sqcup \fG_{n,2n-2}$;  weights (\ref{Ch4: c_gamma_on_C_n}) take now values in $\R[t,dt]$ rather than in $\R$. Proposition
 4.2.6 gives us an explicit example of such a generalized de Rham field theory and of the corresponding {\em path}\ $\caL ie_\infty$-algebra associated with
the propagator
\Beq\label{Ch3: gauge equiv propagator}
\om^{f}(z_1,z_2):= d Arg(z_1-z_2) + t\, d f\hspace{-1mm} \left(\frac{z_i - z_j}{|z_i-z_j|}\right)
+ f\hspace{-1mm}\left(\frac{z_i - z_j}{|z_i-z_j|}\right)dt \in \Omega^1(C_2\times \R).
\Eeq

\mip

{\bf 4.2.8. Homotopy equivalence theorem}. {\em  Let $f$ be an arbitrary smooth
(or $PA$) function on $S^1$ and let
$\mu^f_\bu(t,dt)$ be the associated path $\caL ie_\infty$ structure (\ref{Ch4: mu(t,dt)^f}) on   $\cT_{poly}(\R^d)$.
There exists a morphism of $\caL ie_\infty$-algebras,
$$
H(t,dt): \left(\cT_{poly}(\R^d), [\ \bu\ ]\right) \lon  \left(\cT_{poly}(\R^d)[t,dt], \mu^f_\bu(t,dt)\right)
$$
whose composition with the evaluation map at $t=0$,
$$
\left(\cT_{poly}(\R^d)[t,dt], \mu^f_\bu(t,dt)\right) \stackrel{ev_{t=0}}\lon
\left(\cT_{poly}(\R^d), [\ \bullet\ ]\right),
$$
 equals the identity map.
}

\bip

This theorem can be proven by a direct but very tedious inspection of the polynomial dependence of the weights $c'_\Ga(t)$ on $t$
along the lines of the proof of Proposition 4.2.6. We shall show below in \S 4.5 a short and more elegant  proof (as well as explicit
formulae for $H(t,dt)$)  using the compactified configuration space $\widehat{C}_{n,0}$.

\mip

{\bf 4.2.9. Remark}.
Let $Aut(\cT_{poly}(\R^d))$ be the group of $\caL ie_\infty$ automorphisms of the Schouten algebra.
Two elements, $F_0, F_1\in Aut(\cT_{poly}(\R^d))$ are called {\em homotopy equivalent}\, \cite{Fu} if there exists
a $\caL ie_\infty$ morphism, $F(t,dt)$, from the Lie algebra $(\cT_{poly}(\R^d), [\ \bu\ ])$ to the dg Lie algebra
$(\cT_{poly}(\R^d)[t,dt], [\ \bu\ ], d)$ such that the compositions of   $F(t,dt)$ with the evaluations
at $t=0$ and, respectively,
at $t=1$ give $F_0$ and, respectively, $F_1$.
An element $F=\{F_n\}_{n\geq 1}\in Aut(\cT_{poly}(\R^d)$ is called {\em homotopy trivial}\,
if it is homotopy equivalent to the identity map, and {\em exotic}\, if it is homotopy non-trivial.
If $F_1=\Id$, then the first non-zero higher composition, $F_{min\geq 2}$, of an automorphism $F$ gives a
 cycle in the Chevalley-Eilenberg complex $C_{CE}^\bu(\cT_{poly}(\R^d),\cT_{poly}(\R^d))$   of the Schouten algebra.
 It is easy to check that if the automorphism $F$ is homotopy trivial, then $F_{min}$ is a coboundary.
 Thus if $F_{min}$ gives a non-trivial cohomology class in $H^\bu(\cT_{poly}(\R^d),\cT_{poly}(\R^d))$, then the automorphism is exotic.

\mip


{\bf 4.3. De Rham field theory on $\overline{C}\sqcup \widehat{C}\sqcup \overline{C}$}. For any proper subset $A\subset [n]$ of cardinality
at least two and for any decomposition, $[n]=A_1\sqcup \ldots \sqcup A_k$, of $[n]$ into disjoint
non-empty subsets, there are associated embeddings,
\Beq\label{Ch4: boundary embeddings into C_n,0}
j_A: \widehat{C}_{n-\# A +1,0}\times  \overline{C}_{\# A} \hook \widehat{C}_{n,0},
\ \ \
j_{A_1,...,A_k}: \overline{C}_{k}\times \widehat{C}_{\# A_1,0}\times\ldots\times
\widehat{C}_{\# A_k,0} \hook  \widehat{C}_{n,0},
\Eeq
of the corresponding boundary components into $\widehat{C}_{n,0}$ (see (\ref{Ch3: d on black corollas})).

\mip

 A {\em de Rham field theory on}\, $\overline{C}\sqcup \widehat{C}\sqcup \overline{C}$ is, by definition, a
 pair, $\Omega^{in}$ and $\Omega^{out}$, of de Rham field theories  on
 $\overline{C}$ together with a  family  of maps,
$$
 \left\{
\Ba{rccc}
\Xi: & {\cG}_{n, l}& \lon & \Omega^{l}(\widehat{C}_{n,0})\\
& \Ga & \lon & \Xi_\Ga
\Ea
\right\}_{n\geq 2}
$$
such that  $d\Xi_\Ga=0$, $\Xi_{\Ga_{opp}}=-\Xi_\Ga$ and, for any $\Ga\in  \cG_{n, 2n-3}$,
and any boundary embedding (\ref{Ch4: boundary embeddings into C_n,0})
 one has
\Beq\label{Ch4: boundary in}
j_A^*(\Xi_\Ga) \simeq (-1)^{\sigma_A}\Xi_{\Ga/\Ga_A}  \wedge \Omega^{in}_{\Ga_A},
\Eeq
\Beq\label{Ch4: boundary out}
j_{A_1,...,A_k}^*(\Xi_\Ga)\simeq (-1)^{\sigma_{A_1,...,A_k}}
\Omega^{out}_{\Ga/\{\Ga_{A_1},...,\Ga_{A_k}\}}\wedge \Xi_{\Ga_{A_1}}\wedge\ldots
\wedge \Xi_{\Ga_{A_k}},
\Eeq
where the sign $(-1)^{\sigma_{A_1...A_k}}$ is defined by the equality
$$
\f_\Ga=(-1)^{\sigma_{A_1...A_k}} \f_{\Ga/\{\Ga_{A_1},...,\Ga_{A_k}\}}\wedge \f_{\Ga_{A_1}}
\wedge\ldots\wedge  \f_{\Ga_{A_k}},
$$
i.e.\ it is given just by a rearrangement of the wedge product of edges of $\Ga$.
\mip

{\bf 4.3.1. Theorem.}
{\em Given a de Rham field theory, $(\Omega^{in},\Omega^{out},\Xi)$,
 on $\overline{C}\sqcup \widehat{C}\sqcup \overline{C}$, then, for any $d\in\N$, there are associated
\Bi
\item[(i)]
two  $\caL eib_\infty$-algebra structures,
$\mu^{in}$ and $\mu^{out}$,
on $\cT_{poly}(\R^d)$ given by formulae
(\ref{Ch4: induced_Leib_infty})-(\ref{Ch4: c_gamma_on_C_n})
for $\Omega=\Omega^{in}$ and, respectively, $\Omega=\Omega^{out}$, and
\item[(ii)] a $\caL ieb_\infty$ morphism,
$$
F^{Leib}=
\left\{
\Ba{rccc}
F^{Leib}_n: & \ot^n\cT_{poly}(\R^d) &\lon &  \cT_{poly}(\R^d)[2-2n]\\
& \ga_1\ot \ldots\ot \ga_n     & \lon & F^{Leib}_n(\ga_1, \ldots, \ga_n)
\Ea
\right\}_{n\geq 1},
$$
from $\mu^{in}$-structure to  $\mu^{out}$-structure
given by the formulae,
\Beq\label{Ch4: induced_Aut_Leib_infty}
F^{Leib}_n(\ga_1, \ldots, \ga_n):=
\left\{
\Ba{cl}
\Id & \mbox{for}\ n=1,\\
 \sum_{\Ga\in \fG_{n,2n-2}} C_{\Ga}\Phi_\Ga(\ga_1, \ldots, \ga_n)
 & \mbox{for}\ n\geq 2
 \Ea
 \right.
\Eeq
with
\Beq\label{Ch4: weight C_Ga}
C_\Ga:= \int_{\widehat{C}_{n,0}} \Xi_\Ga.
\Eeq
\Ei
}
\begin{proof} Proof is completely analogous to
the proof of Theorem~4.2.1 above. The required claim follows immediately from the definitions
and the Stokes theorem,
\Beqr \label{Ch4: main}
0&=&\int_{\widehat{C}_{n,0}}d\Xi_\Ga =\int_{\p\widehat{C}_{n,0}}\Xi_\Ga=\nonumber \\
&=&\hspace{-2mm}
-\hspace{-1mm}\sum_{A\varsubsetneq [n]\atop \# A\geq 2}(-1)^{\sigma_A}\hspace{-2mm}
\underset{\overline{C}_{\# A}}{\int}\hspace{-1mm}\Omega^{in}_{\Ga_A}\hspace{-2mm}
\underset{\widehat{C}_{n-\# A+1,0}}
{\int}\hspace{-5mm}\Xi_{\Ga/\Ga_A}
 +\hspace{-2mm}\sum_{V(\Ga)=A_1\sqcup ...\sqcup A_k\atop 2\leq k\leq n}\hspace{-4mm}
 (-1)^{\sigma_{A_1...A_k}}
 \underset{\overline{C}_{k}}{\int}\Omega^{out}_{\Ga/\{\Ga_{A_1}...\Ga_{A_k}\}}
  \hspace{-3mm}\underset{\widehat{C}_{\# A_1,0}}{\int}\hspace{-3mm}\Xi_{\Ga_{A_1}}...\hspace{-2mm}
  \underset{\widehat{C}_{\# A_k,0}}{\int}\hspace{-3mm}\Xi_{\Ga_{A_k}} \nonumber
\\
&=&\hspace{-2mm}
-\hspace{-4mm}
\sum_{A\subset V(\Ga)\atop
A\, \mathit i \mathit s\, \mathit a \mathit d \mathit m \mathit i \mathit s \mathit s \mathit i \mathit b \mathit l \mathit e}
\hspace{-2mm}
(-1)^{\sigma_{A}}c^{in}_{\Ga_A}C_{\Ga/\Ga_A}  +
\sum_{k=2}^n\ \sum_{V(\Ga)=\underbrace{A_1\sqcup ...\sqcup A_k}_{
\mathit a \mathit d \mathit m \mathit i \mathit s \mathit s \mathit i \mathit b \mathit l \mathit e
}}\hspace{-3mm}
(-1)^{\sigma_{A_1...A_k}}
c^{out}_{\Ga/\{\Ga_{A_1}...\Ga_{A_k}\}}C_{\Ga_{A_1}}\ldots C_{\Ga_{A_k}}.
\Eeqr
\end{proof}

{\bf 4.4. De Rham field theories from angular functions on $\widehat{C}_{2,0}$}.
Let $\phi$ be an angular function on
 $\widehat{C}_{2,0}$, $\om=\phi^*(vol_{S^1})$ the
 associated propagator on $\widehat{C}_{2,0}$, and let
$$
\om_{in}:= \om|_{inner\ circle} \ \ \ \mbox{and}\ \ \  \om_{out}:= \om|_{outer\ circle}
$$
be the 1-forms on $S^1$ obtained by restricting $\om$ to the inner and, respectively, outer circles
of $\widehat{C}_{2,0}$.
Define a series of maps,
\Beq\label{Ch4: Omega-in maps}
\Ba{rccc}
\Omega^{in}: & \cG_{n,l}\hspace{-2mm} & \rar &\hspace{-2mm} \Omega^l(\overline{C}_n)\vspace{2mm} \\
& \Ga \hspace{-2mm}& \rar & \hspace{-2mm}\Omega^{in}_\Ga:=\hspace{-3mm}\displaystyle
\bigwedge_{e\in E(\Ga)}\hspace{-2mm}
\frac{{\pi}^*_e\left(\om_{in}\right)}{2\pi}
\Ea
\ \ \ \ \ \ \ \ \ \ \
\Ba{rccc}
\Omega^{out}: & \cG_{n,l}\hspace{-2mm} & \rar &\hspace{-2mm} \Omega^l(\overline{C}_n)\vspace{2mm}  \\
& \Ga \hspace{-2mm}& \rar & \hspace{-2mm}\Omega^{in}_\Ga:=\hspace{-3mm}\displaystyle
\bigwedge_{e\in E(\Ga)}\hspace{-2mm}
\frac{{\pi}^*_e\left(\om_{out}\right)}{2\pi}
\Ea
\Eeq
and
\Beq\label{Ch4: Xi maps}
\Ba{rccc}
\Xi: &  \cG_{n,l} & \lon & \Omega^l(\widehat{C}_{n,0})\vspace{2mm} \\
& \Ga & \lon & \Xi_\Ga:=\displaystyle\hspace{-3mm} \bigwedge_{e\in E(\Ga)}\hspace{-1mm}
\frac{{\fp}^*_e\left(\om\right)}{2\pi}
\Ea
\Eeq
where, as before,
${\pi}_e: {C}_n\rar  {C}_2$ is the  map which  forgets all the points in the configurations except the two ones which
are the boundary vertices of the edge $e$, and
${\fp}_e: {C}_{n,0}\rar  {C}_{2,0}$ is a ``renormalized" version (see (\ref{Ch3: fp-forgetting map}))
of the  forgetful map $p_e: {C}_{n,0}\rar  {C}_{2,0}$. As Mayer-Vietoris sequence for minimal differential forms on locally closed semialgebraic
sets trivializes (see \S 8.3 in \cite{KoSo}), ${\Xi}_\Ga$ is a well-defined minimal differential
form on $\widehat{C}^{}_{n,0}$ for any graph
$\Ga$ with $n$ vertices.

\mip
{\bf 4.4.1. Theorem}. {\em For any angular function on $\widehat{C}_{2,0}$ the
associated  data (\ref{Ch4: Omega-in maps})-(\ref{Ch4: Xi maps})
define a de Rham field theory on  $\overline{C}\sqcup \widehat{C}\sqcup \overline{C}$.
}

\begin{proof} The proof of the factorization (\ref{Ch4: boundary in})  is standard (see \cite{Ko}):
this equation is equivalent to the following one,
$$
\int_{ \widehat{C}_{n-\# A +1,0}\times  \overline{C}_{\# A}}
j_A^*(\Xi_\Ga)= (-1)^{\sigma_A}\int_{{C}_{n-\# A +1,0}} \Xi_{\Ga/\Ga_A} \int_{C_{\# A}}  \Omega^{in}_{\Ga_A}.
$$
By definition of the boundary strata $\widehat{C}_{n-\# A +1,0}\times  \overline{C}_{\# A} \hook \widehat{C}_{n,0}$, both sides of the above equation are zero unless $\Ga_A$ is an admissible subgraph
of $\Ga$ in which case the equality is obvious when one uses local coordinates defined in \S 3.

\sip

Consider next, for a partition $[n]=A_1\sqcup \ldots\sqcup A_k$, the associated boundary strata at ``infinity",
$$
j_{A_1,...,A_k}: {C}_{k}\times {C}_{\# A_1,0}\times\ldots\times
{C}_{\# A_k,0} \hook  \widehat{C}_{n,0}.
$$
By definition, the boundary stratum $j_{A_1,...,A_k}({C}_{k}\times {C}_{\# A_1,0}\times\ldots\times
{C}_{\# A_k,0})$ is a subset of $\widehat{C}_{n,0}$, obtained in the limit
 $\var\rar\infty$ from a class of configurations, $p_\var(p_0,p_1,\ldots, p_k)$, in $C_{n,0}$
 defined as follows:
\Bi
\item[(a)]  let $p_0=(z_1, \ldots, z_k)$ be an arbitrary configuration in $C^h_k$ and let
$$
p_0^\var:= \var (p_0- \ii) + \ii= (z^\var_1:=\var(z_1-\ii)+\ii , \ldots, z^\var_k:=\var(z_k-\ii) +\ii )
$$
 be the associated  configuration
in $C^h_{k,0}$ with $||p_0^\var||=\var$;
\item[(b)]  let $p_1\in  {C}^h_{ A_1,0}$,
$\ldots$, $p_k\in   {C}^h_{ A_k,0}$ be arbitrary configurations,
\Ei
then, for sufficiently large $\var$,
$$
p_\var(p_0,p_1,\ldots, p_k):= \bigcup_{i=1}^k \{\underbrace{  p_i + \var(z_i-\ii)}_{p_{A_i}} \}
$$
is a well-defined configuration in $C^h_{n,0}$. Note that $z_{min}(p_{A_i})= z_i^\var$
and $|p_{A_i}-z_{min}(p_{A_i})|$ is equal to $|p_i-\ii|$ for all $i\in [k]$ and hence is a finite number
independent of
$\var$.

\sip

 Therefore, for any $\Ga\in \fG_{n,l}$
we have
\Beq\label{3: Equality for Xi'}
 j_{A_1,...,A_k}^*(\Xi_\Ga)=(-1)^{\sigma_{A_1,,...,A_k}} \lim_{\var\rar +\infty}\left(
\bigwedge_{e\in E(\Ga/\{\Ga_{A_1}\ldots \Ga_{A_k}\} )}\hspace{-1mm}
\frac{{\fp}^*_e\left(\om\right)}{2\pi}\ \
\prod_{i=1}^k
\bigwedge_{e\in E(\Ga_{A_i})}\hspace{-1mm}
\frac{{\fp}^*_e\left(\om\right)}{2\pi}\right).
\Eeq
Let $e$ be  an edge in $E(\Ga/\{\Ga_{A_1}\ldots \Ga_{A_k}\} )$, then $In(e)$ and $Out(e)$ correspond to points, say $w_{p_0}$ and $w_{q_0}$, which belong to two {\em different}\,
groups, say $p_{A_p}$ and $p_{A_q}$.
Thus
$$
w_{p_0}= z^h_{p_0} + \var(z_{p} -\ii), \ \ \ \
w_{q_0}= z^h_{q_0} + \var(z_{q} -\ii),
$$
for some uniquely defined  $z^h_{p_0}\in C^h_{A_p,0}$, $z^h_{q_0}\in C^h_{A_q,0}$,  $z_{p}, z_q\in C_k^h$.
Then, for sufficiently large $\var$, we have,
\Beqrn
{\fp}^*_e\left(\om\right)(w_{p_0}, w_{q_0})&=& \displaystyle\left\{  \Ba{cc}
\om\left(w_{p_0}- w_{q_0}+ \ii, \ii\right) & \mbox{if} \ \Im z_p \geq  \Im z_q\\
\om \left(\ii, w_{q_0}- w_{p_0}+ \ii\right)  & \mbox{if} \ \ \Im z_p \leq  \Im z_q
\Ea
\right.\\
&=& \displaystyle\left\{  \Ba{cc}
\om\left(\var(z_p-z_q +\frac{z^h_{p_0}- z^h_{q_0}}{\var}) + \ii, \ii\right) & \mbox{if} \ \Im z_p \geq  \Im z_q\\
\om \left(\ii, \var(z_q-z_p +\frac{z^h_{q_0}- z^h_{p_0}}{\var}) + \ii \right)  & \mbox{if} \ \ \Im z_p \leq  \Im z_q
\Ea
\right.\\
&\underset{\var\rar +\infty}{\lon} & \om_{out}(z_p,z_q)
\Eeqrn
so that
$$
\bigwedge_{e\in E(\Ga/\{\Ga_{A_1}\ldots \Ga_{A_k}\} )}\hspace{-1mm}
\frac{\underset{\var\rar +\infty}{\lim}{\fp}^*_e\left(\om\right)}{2\pi}=
\Omega^{out}_{\Ga/\{\Ga_{A_1},...,\Ga_{A_k}\}}.
$$

It is also clear that for any $i\in [k]$ and any edge $e$ in $E(\Ga_{A_i})$ the 1-form
${\fp}^*_e\left(\om\right)$ is independent of $\var$ so that, in the above notations,
$$
{\fp}^*_e\left(\om\right)(w_{p_0}, w_{q_0})= {\fp}^*_e\left(\om\right)(z^h_{p_0}, z^h_{q_0})
$$
and hence
$$
\bigwedge_{e\in E(\Ga_{A_i})}\hspace{-1mm}
\frac{\underset{\var\rar +\infty}{\lim}{\fp}^*_e\left(\om(w_{In(v)}, w_{Out(v)}\right)}{2\pi}=
\bigwedge_{e\in E(\Ga_{A_i})}\hspace{-1mm}
\frac{{\fp}^*_e\left(\om(z^h_{In(v)}, z^h_{Out(v)}\right)}{2\pi}= \Xi_{\Ga_{A_i}}.
$$
 Then equality (\ref{3: Equality for Xi'}) implies
equality (\ref{Ch4: boundary out}).
\end{proof}

\mip

{\bf 4.4.2. Remark}. The above proof shows that our maps $\Omega^{in}$, $\Omega^{out}$ and $\Xi$ satisfy factorization equations
(\ref{Ch4: boundary in}) and (\ref{Ch4: boundary out}) strictly, i.e.\ one can replace there the equivalence sign $\simeq$ with the equality sign $=$. This fact means that the class of de Rham field theory on $\overline{C}\sqcup \widehat{C}\sqcup \overline{C}$ associated with angular functions  can be used to construct representations of not only the operad of fundamental chains, but also
 of the {\em full}\, chain operad,
 $Chains(\overline{C}\sqcup \widehat{C}\sqcup \overline{C})$.

\mip

{\bf 4.4.3. Corollary}. {\em For any angular function $\phi$ on $\widehat{C}_{2,0}$ there is an
associated  $\caL ie_\infty$ automorphism of the Schouten algebra of
polyvector fields.}

\begin{proof}
As weights (\ref{Ch4: weight C_Ga}) are $\bS_n$-invariant,
 Proposition~4.2.2 implies that the map $F$ given by formulae (\ref{Ch4: induced_Aut_Leib_infty}) describes a $\caL ie_\infty$ morphism between the $\caL ie_\infty$
 structures $\mu_\bullet^{in}$ and $\mu^{out}_\bullet$ on  $\cT_{poly}(\R^d)$ corresponding to the 1-forms
 $\om^{in}$ and $\om^{out}$ respectively. By Theorem 4.2.8, both these $\caL ie_\infty$ structures
 are homotopy equivalent to the Schouten Lie algebra via certain homotopy equivalence maps, $H_{in}$ and $H_{out}$ (see the next section \S 4.5  for their explicit formulae).  Then the composition
 $$
\cF:
\left(\cT_{poly}(\R^d),[\ \bullet \ ]\right)\stackrel{H_{in}}{\lon}
\left(\cT_{poly}(\R^d),\mu^{in}_\bullet\right)
)\stackrel{F}{\lon} \left(\cT_{poly}(\R^d),\mu^{out}_\bullet\right) \stackrel{H_{out}^{-1}}{\lon}
\left(\cT_{poly}(\R^d),[\ \bullet \ ]\right)
$$
is an automorphism of the Schouten algebra.
\end{proof}
\sip

{\bf 4.4.4. An example}. Let us test Main Theorem~4.4.1 for the Kontsevich propagator $\om_K$.
The statement of that theorem is equivalent to saying that, for any graph $\Ga\in \fG_{n,2n-3}$, the associated weights
$$
c^{in}_\Ga=\int_{\overline{C}_n}\bigwedge_{e\in E(\Ga)}\hspace{-2mm}
\frac{{\pi}^*_e\left(\om_K|_{inner\ circle}\right)}{2\pi},\ \  c^{out}_\Ga= \int_{\overline{C}_n}\bigwedge_{e\in E(\Ga)}\hspace{-2mm}
\frac{{\pi}^*_e\left(\om_K|_{outer\ circle}\right)}{2\pi},\  \  C_\Ga=
\int_{\widehat{C}_{n,0}}\bigwedge_{e\in E(\Ga)}\hspace{-2mm}
\frac{{\fp}^*_e\left(\om_K\right)}{2\pi},\
$$
satisfy equation (\ref{Ch4: main}). For $n=1$ this equation is trivial. For the cases $n=2$ and $n=3$
it is obvious. The first non-obvious numerical identities come from graphs $\Ga$ in the set $\fG_{4,5}$.
 For the three graphs in (\ref{Ch4: three graphs}) the associated equations (\ref{Ch4: main}) reduce to the following ones,
$$
c^{out}_{\Ga_1}= 2C_{\Ga'}, \ \ \ c^{out}_{\Ga_3}= 2C_{\Ga''}, \ \ \ c^{out}_{\Ga_2}= C_{\Ga'} + C_{\Ga''},
$$
where
$$
\Ga'=  \Ba{c}\xy
(0,-2)*{_1},
(-13.5,6)*{^2},
(0,14)*{^3},
(0,0)*{\bullet}="1",
(-12,6)*{\bullet}="2",
(0,12)*{\bullet}="3",
   {\ar@/^0.6pc/(0,12)*{\bu};(0,0)*{\bu}};
 {\ar@/^0.6pc/(0,0)*{\bu};(0,12)*{\bu}};
\ar @{->} "2";"3" <0pt>
\ar @{<-} "1";"2" <0pt>
\endxy\Ea,
\ \ \ \ \ \ \
\Ga''=  \Ba{c}\xy
(0,-2)*{_1},
(-13.5,6)*{^2},
(0,14)*{^3},
(0,0)*{\bullet}="1",
(-12,6)*{\bullet}="2",
(0,12)*{\bullet}="3",
   {\ar@/^0.6pc/(0,12)*{\bu};(0,0)*{\bu}};
 {\ar@/^0.6pc/(0,0)*{\bu};(0,12)*{\bu}};
\ar @{<-} "2";"3" <0pt>
\ar @{->} "1";"2" <0pt>
\endxy\Ea
$$
Thus the Main Theorem 4.4.1 implies the identity,
$$
 c^{out}_{\Ga_2}=\frac{1}{2}\left(  c^{out}_{\Ga_1} +  c^{out}_{\Ga_3}  \right)
$$
which indeed holds true as all these weights are equal to $\frac{1}{12}$ (see \S 4.2.4). Similarly one can check
several other identities between the weights computed in Appendix 4.

\mip

{\bf 4.5. Gauge equivalence propagators and  a proof of the homotopy equivalence theorem.} A semialgebraic function,
$$
\Ba{rccc}
\fl: &  C_{2,0} &\lon &  (0,1)\\
&p= (z_1, z_2)  & \lon & \frac{||p||}{||p||+1}
\Ea
$$
extends to a semialgebraic function on its compactification,  $\fl: \widehat{C}_{2,0} \rar [0,1]$, which takes
values $0$ at the inner circle and value $1$  at the outer one. For an arbitrary semialgebraic function $f$ on $S^1$
we consider a propagator
$$
\om(z_1,z_2,t,dt)= dArg(z_1-z_2) + t d\left(\hspace{-1mm}\fl(z_1,z_2)\,f
\hspace{-1mm}\left(\frac{z_i - z_j}{|z_i-z_j|}\right)\right) +
\fl(z_1,z_2)f\hspace{-1mm}\left(\frac{z_i - z_j}{|z_i-z_j|}\right)dt \in \Omega^1(\widehat{C}_{2,0}\times \R)
$$
which satisfies the boundary conditions,
$$
\om(z_1,z_2,t,dt)|_{inner\ circle}= dArg(z_1-z_2), \ \ \
\mbox{and}\ \ \ \om(z_1,z_2,t,dt)|_{outer\ circle}= (\ref{Ch3: gauge equiv propagator}).
$$
Hence formulae (\ref{Ch4: induced_Aut_Leib_infty}) with summation $\sum_{\Ga\in \fG_{n,2n-2}}$ extended to
 $\sum_{\Ga\in \fG_{n,2n-2}\sqcup \fG_{n,2n-1}}$ give us a $\caL ie_\infty$-morphism
$$
H(t,dt): \left(\cT_{poly}(\R^d), [\ \bu\ ]\right) \lon  \left(\cT_{poly}(\R^d)[t,dt], \mu_\bu^f(t,dt)\right)
$$
which obviously has the property stated at the end of Theorem 4.2.8.
\sip

The map $H_{out}$ used in the proof of Corollary 4.4.3 is equal to  $H(t,dt)|_{t=1}$. The map $H_{in}$ is constructed similarly.

\mip

\mip
{\bf 4.6. Exotic transformations of Poisson structures}. Any $\caL ie_\infty$-automorphism
$F$ of the algebra $\cT_{poly}(\R^d)[[\hbar]]$ acts on its set of Maurer-Cartan elements,
$$
\ga \rar F(\ga)=\sum_{n\geq 1}\frac{\hbar^{n-1}}{n!} F_n(\ga,\ldots,\ga).
$$
If $F$ is a $\caL ie_\infty$-automorphism given by a de Rham field theory
on $\overline{C}_\bu\sqcup \widehat{C}_{\bu,0}\sqcup \overline{C}_\bu$, then
\Beq\label{ch4: action on MC elements}
 F(\ga):= \ga +
 \sum_{n\geq 2}\hbar^{n-1}\sum_{\Ga\in G_{n,2n-2}} \frac{C_\Ga}{\# Aut (\Ga)} \Phi_\Ga(\ga^{\ot n})
\Eeq
where $\# Aut(\Ga)$ is the cardinality of the group of automorphisms of the graph $\Ga$ and $G_{n,2n-2}$ means the family of graphs
$\cG_{n,2n-2}$ with labeling of vertices forgotten.
In particular, $F$ acts on the set of ordinary Poisson structures on $\R^d$, that is,
on the set of bivector fields,
$\ga=\frac{1}{2}\sum_{i,j}\ga^{ij}(x)\psi_i\psi_j$, satisfying the equation $[\ga\bu\ga]=0$. In this case
 only those graphs $\Ga\in G_{n,2n-2}$ can give a non-trivial contribution to
  $F(\ga)$ which have at most two output edges at each vertex. The wheels,
\Beq\label{Ch4:w_n}
w_n=\ \ \xy
(5,0.7)*{^{n+1}},
(-12,13)*{^{1}},
(0,17)*{^{2}},
(12,13)*{^{3}},
(-17,0)*{^{n}},
(-12,-14)*{^{n-1}},
 (0,0)*{\bullet}="1",
(-11,11)*{\bullet}="2",
(0,15)*{\bullet}="3",
(11,11)*{\bullet}="4",
(15,0)*{\bullet}="5",
(11,-11)*{\bullet}="6",
(0,-15)*{\bullet}="7",
(-15,0)*{\bullet}="n",
(-11,-11)*{\bullet}="n-1",
\ar @{<-} "1";"2" <0pt>
\ar @{<-} "1";"3" <0pt>
\ar @{<-} "1";"4" <0pt>
\ar @{<-} "1";"5" <0pt>
\ar @{<-} "1";"6" <0pt>
\ar @{<-} "1";"7" <0pt>
\ar @{<-} "1";"n" <0pt>
\ar @{<-} "1";"n-1" <0pt>
\ar @{<-} "2";"3" <0pt>
\ar @{<-} "3";"4" <0pt>
\ar @{<-} "4";"5" <0pt>
\ar @{<-} "5";"6" <0pt>
\ar @{<-} "6";"7" <0pt>
\ar @{<-} "7";"n-1" <0pt>
\ar @{<-} "n-1";"n" <0pt>
\ar @{<-} "n";"2" <0pt>
\endxy
\ \ , \ \ \
n\geq 2,
\Eeq
and their unions
do have this property. 
An easy calculation based on definition (\ref{Ch4: Phi_Ga}) gives
(with the total ordering
  of $E(w_n)$ chosen  to be $\{(1,2),(2,3),\ldots,(n-1,n), (1,n+1),\ldots, (n,n+1)\}$),
$$
\Phi_{w_n}(\ga^{\ot n+1}) = -(-1)^{n(n-1)/2}\frac{1}{2}\sum \frac{\p^n \ga^{ij}}{\p x^{k_1}\cdots\p x^{k_n}}
\frac{\p \ga^{k_1 l_1 }}{\p x^{l_2}}
\frac{\p \ga^{k_2 l_2 }}{\p x^{l_3}}\cdots
\frac{\p \ga^{k_n l_n }}{\p x^{l_{1}}}
(\psi_i\psi_j).
$$

As  $Aut(w_n)=\Z/n\Z$,  the contribution of
graphs $w_n$ into
into $F$ is given by
\Beqr
F(\ga)\hspace{-3mm}
&=&\ga+\sum_{n\geq 2}\hbar^{n}\frac{C_{w_n}}{\# Aut(w_n)}\Phi_{w_n}(\ga^{\ot n+1})+\ldots\nonumber\\
 &=&\hspace{-1mm}
  \ga-\frac{1}{2}\sum_{n=2}^\infty (-1)^{n(n-1)/2}\frac{C_{w_n}\hbar^{n}}
 {n}
\frac{\p^{n} \ga^{ij}}{\p x^{k_1}\cdots\p x^{k_{2n+1}}}
\frac{\p \ga^{k_1 l_1 }}{\p x^{l_2}}
\frac{\p \ga^{k_2 l_2 }}{\p x^{l_3}}\cdots
\frac{\p \ga^{k_{n} l_{n}}}{\p x^{l_{1}}}
(\psi_i\psi_j)+... \label{Ch4: F for wheels}
\Eeqr

\mip

{\bf 4.7. Examples: symmetrized Kontsevich's  propagators}. The symmetrized
propagator (see \S 3.5.5),
$$
\om_{{K}}^{sym}=\frac{1}{2}\left(\om_K(z_1,z_2) + \om_K(z_2,z_1)\right)=d\mbox{Arg} (z_1-z_2),
$$
on $\widehat{C}_{2,0}$ restricts to its inner and outer boundary circles as $d\mbox{Arg} (z_1-z_2)$
and hence defines, by Example~4.2.3 and Theorem 4.4.1,  an automorphism, $F_{{K}}^{sym}$,
 of
the Schouten algebra $\cT_{poly}(\R^d)$. This is, however,
a  trivial automorphism (i.e.\ the one with  $F_{n\geq 2}=0$) as the differential
$(2n-2)$-forms $\prod_{e\in E(\Ga)} \fp_e^*( \om_{{K}}^{sym})$
are invariant under the action of the semigroup, $z\rar  z + \nu$, $\nu\in \R^+$, and hence
vanish identically on $C_{n,0}$ for dimensional reasons.

\mip

 In \cite{Ko3} Kontsevich introduced a ``$\frac{1}{2}$"-propagator,
 $$
\om_{\frac{1}{2}K}(z_i,z_j):=\frac{1}{i}d\log
\frac{z_i-z_j}{\overline{z}_i - z_j}\ ,
$$
 and claimed that ``all
identities proven in \cite{Ko} remain true". In particular, all the integrals
$\int_{C_{n,0}}\wedge_{e\in E(\Ga)}p_e^*(\om)$,
$\Ga\in \cG_{n,2n-2}$, are finite.
This is by no means an obvious claim as the differential
forms $p_{e}^*(\om_{\frac{1}{2}K})$, $i,j\in[n]$, extend neither to the compactification
$\overline{C}_{n,0}$   nor
to $\widehat{C}_{n,0}$. In fact such forms extend nicely to all boundary components of both
compactifications except to those of the form $C_{n-k+1,0}\times C_{k}$ which describe a group of $k$
points moving too close to each other in $\bbH$ (and which are the only ones which are common to
 both compactifications for all $n$). We refer to \cite{AT}
for a discussion of why the Kontsevich  $\frac{1}{2}$-propagator works.
As the symmetrized version of this propagator,
$$
\om^{sym}_{\frac{1}{2}K}(z_i,z_j):=\frac{1}{2}\left(\om_{\frac{1}{2}K}(z_i,z_j) + \om_{\frac{1}{2}K}(z_j,z_i)\right),
$$
restricts to the outer circle of $\widehat{C}_{2,0}$ as  $d\mbox{Arg} (z_1-z_2)$ and tends towards the inner circle
 as  $d\mbox{Arg} (z_1-z_2) + d\ln \var$, $\var \rar 0$, we infer from Example~4.2.3 and Theorem~{4.3.1}
 that the associated
 universal map
$F^{sym}_{\frac{1}{2}K}$ given by formulae (\ref{Ch4: induced_Aut_Leib_infty})-(\ref{Ch4: weight C_Ga})
gives an exotic automorphism of the Schouten algebra without any homotopy adjustments.
The automorphism $F_{\frac{1}{2}K}^{sym}$ is homotopy non-trivial
 as its lowest in $n$ non-trivial component
is given by the graph $w_3$ whose weight with respect to the propagator $\om^{sym}_{\frac{1}{2}K}$
and the ordinary forgetful map
 is equal to  $\frac{\zeta(3)}{(4\pi)^3\ii}$ \cite{Gr2}; one can use this fact to show that its weight is also non-zero with respect to the renormalized forgetful map $\fp_e$.
It is worth noting that for any  symmetrized propagator the choice of arrows on a graph $\Ga\in \cG_{n,2n-2}$ does not affect its weight $C_\Ga$ (but does affect the associated operator $\Phi_\Ga$). Such theories are better understood as de
Rham field theories on {\em braid}\, configuration spaces, see \S 5 below.
\mip

{\bf 4.8.  Example:  Kontsevich's (anti)propagator}. In the case of the Kontsevich propagator
$$
\om_{{K}}(z_1,z_2)=d\mbox{Arg} \frac{z_1-z_2}{\overline{z_1} -z_2}
$$
or the antipropagator,
$$
\om_{\overline{K}}(z_1,z_2)=d\mbox{Arg} \frac{z_2-z_1}{\overline{z_2} -z_1}
$$
formulae  (\ref{Ch4: induced_Aut_Leib_infty})-(\ref{Ch4: weight C_Ga}) define a $\caL ie_\infty$ morphism from the Schouten algebra to its $\caL ie_\infty$-extension constructed
by Shoikhet in \cite{Sh2} (see \S 4.2.4). Both such morphisms, $F_{K}$ and, respectively,
 $F_{\overline K}$,  must be highly non-trivial
as they encode all the obstructions to existence of universal Kontsevich type
formality morphism for {\em infinite-dimensional}\, Schouten algebras (non-existence of a such a formality morphism,
i.e.\ non-emptiness of the set of obstructions,  was proven in \cite{Me-lec}).

\mip

{\bf 4.9. De Rham field theory of Duflo's strange automorphism}.
Let $\om$ be a propagator on $\widehat{C}_{2,0}$.
Had we defined the map $\Xi$ in
(\ref{Ch4: Xi maps})  with the help
of the ordinary forgetful map $p_e: C_{n,0}\rar C_{2,0}$
(rather than with $\fp_e$),
\Beq\label{Ch4: Xi maps nonrenorm}
\Ba{rccc}
\Xi': &  \cG_{n,l} & \lon & \Omega^l(\widehat{C}_{n,0})\vspace{2mm} \\
& \Ga & \lon & \Xi'_\Ga:=\displaystyle\hspace{-3mm} \bigwedge_{e\in E(\Ga)}\hspace{-1mm}
\frac{{p}^*_e\left(\om\right)}{2\pi}
\Ea
\Eeq
then we would not get in general a de Rham field theory
$(\Omega^{in}, \Xi', \Omega^{out})$ on $\widehat{C}_{n,0}$
as for generic graphs $\Ga$  the factorization (\ref{Ch4: boundary out}) might fail.
For a fixed propagator $\om$, let us denote by $\cG_{sing}(\om)$ that set of graphs $\Ga$
for which it fails indeed.

\sip

The boundary values of the propagator $\om$ determine the associated $\mu_\bu^{in}$
and  $\mu_\bu^{in}$ $\caL ieb_\infty$-structures on $\R^d$ as it is explained in \S 3.
It is clear that if $\ga$ is a Maurer-Cartan
element of the  $\mu_\bu^{in}$ structure
such that for any $\Ga\in \cG_{sing}(\om)$ the value,
$\Phi_\Ga(\al, \ldots,\al)$, of the associated operator $\Phi_\Ga$ vanishes, then the transformation
\Beq\label{ch4: formular for partial de rham auto}
 F(\ga):= \ga +
 \sum_{n\geq 2}\hbar^{n-1}\sum_{\Ga\in G_{n,2n-2}} \frac{C'_\Ga}{\# Aut (\Ga)} \Phi_\Ga(\ga^{\ot n})
\Eeq
$$
C'_\Ga:=\int_{C_{n,0}} \Xi'_\Ga,
$$
defines a Maurer-Cartan element, $F(\al)$, of the
the $\mu^{out}_\bu$-structure on $\R^d$, i.e.\ in such a case our machinery works with non-renormalized weights
$C'_\Ga$.

\sip

As an illustration, let
consider a family, $\{\ga^{P.D.}\}$, of polyvector fields on $\R^d$  of the form $\ga^{P.D.}=\sum_{i\geq 0}\ga^i$, $\ga^i\in
\wedge^i\cT_{\R^d}$, with all $\ga^i$ vanishing
 except for $i= 0$ and $2$, and with $\ga^2=\frac{1}{2}\sum_{i,j}\al^{ij}_k x^k \psi_i\psi_j$ being a linear
Poisson structure. Equation $[\ga,\ga]_{Schouten}=0$  implies then that $\ga^0=\ga^0(x)$ is an invariant
polynomial on $\R^d$, that is, an element of $(\odot^\bullet \fg)^\fg$, where $\fg$ is the space dual
 to $\R^d$ and  equipped  with the Lie algebra structure determined by $\ga^2$.
It is not hard to see that this class of {\em Poisson-Duflo} structures, $\{\ga^{P.D.}\}$, does satisfy
the condition $\Phi_\Ga(\ga^{P.D.}, \ldots,\ga^{P.D.})=0$
for any $\Ga\in \cG_{sing}(\om_{\overline{K}})$. Therefore, the Kontsevich antipropagator defines an automorphism (\ref{ch4: formular for partial de rham auto}) of Poisson structures,
$$
\ga^{P.D.}\lon F^{\overline{K}}(\ga^{P.D.})
$$
which can be computed explicitly as, by Kontsevich vanishing theorems
in \cite{Ko}, only graphs (\ref{Ch4:w_n}) and their unions (with the same center) may contribute to
(\ref{ch4: formular for partial de rham auto}),
\Beqrn
F^{\overline{K}}(\ga^{P.D.})&=&\ga^{P.D.} + \sum_{n\geq 1}\sum_{\Ga\in G_{n+1, 2n-1}}
 \frac{\hbar^nC_\Ga'}{\# Aut(\Ga)}
 \Phi_\Ga({\ot^{ n+1}}\ga^{P.D.})\\
&=&\ga^{P.D.} +  \sum_{m=1}^\infty \frac{1}{m!}\left(\sum_{n\geq 2} \frac{\hbar^n}{n}C_{w_n}'\Phi_{w_n}
\left({\ot^{ n+1}}\ga^{P.D.}\right)\right)^m\\
\\
 &=& \ga^0 + \ga^2 + \sum_{m=1}^\infty \frac{1}{m!}\left( \sum_{n\geq 2} \cB_n\hbar^n
\frac{\p^n \ga^0}{\p x^{k_1}\cdots\p x^{k_n}}
 \al^{k_1 l_1}_{l_2}
\al^{k_2 l_2 }_{l_3}\cdots
\al^{k_n l_n }_{l_{1}}\right)^m
\\
&=& \ga^2 + e^{\sum_{n\geq 2} \cB_n \hbar
{\mathrm T \mathrm r \mathrm a \mathrm c \mathrm e}(\ad^n)}
\ga^0\\
&=& \ga^2 + \det\sqrt{\frac{e^{\frac{\hbar}{2}\ad}- e^{-\frac{\hbar}{2}\ad}}{\ad}}
\ga^0.
\Eeqrn
Here we used the fact that the weight,
 $C_{w_n}'$,  of wheel (\ref{Ch4:w_n}) with respect to  Kontsevich's antipropagator
(and the ordinary projection $p_e:C_{n,0}\rar C_{2,0}$)  is zero for $n$ odd (see \cite{Ko}) and equals to
 $
C_{w_n}'= -(-1)^{n(n-1)/2}n{\cB_n}$
for $n$ even  \cite{VdB}.
The weight of a union of $m$ wheels is equal to the product of their weights giving rise above to the summation
over $m\geq 1$.

\sip

The conclusion is that the exotic transformation  $F^{\overline{K}}$ associated with Kontsevich's antipropagator
preserves the class of Poisson-Duflo structures, and, at $\hbar=1$, coincides  precisely
with the famous {\em strange Duflo automorphism}
 (see, e.g., \cite{Ko,CaRo} and references cited there).

\mip

Analogously, the transformation of Poisson-Duflo structures,
$$
\ga^{P.D.}\lon F^{\frac{1}{2}\overline{K}}(\ga^{P.D.})
$$
 associated  with Kontsevich's $\frac{1}{2}$-antipropagator,
$$
\om_{\frac{1}{2}\overline{K}}(z_i,z_j)=\frac{1}{i}d\log
\frac{z_i-z_j}{z_i -\overline{z}_j}
$$
is given  by
\Beqrn
F^{\frac{1}{2}\overline{K}}(\ga^{P.D.})
&=& \ga^2 + e^{\sum_{n\geq 2} \frac{\zeta(n)}{n
(2\pi \ii)^n}\hbar
{\mathrm T \mathrm r \mathrm a \mathrm c \mathrm e}(\ad^n)}
\ga^0
\Eeqrn
and at $\hbar=1$ coincides  with Kontsevich's modification
(see \S 4.6 in \cite{Ko3})
of Duflo's strange transformation. Here we used the fact that the weight,
 $C_{w_n}'$,  of wheel (\ref{Ch4:w_n}) with respect to  Kontsevich's $\frac{1}{2}$-antipropagator
(and the ordinary projection $p_e:C_{n,0}\rar C_{2,0}$) is equal to
$(-1)^{n(n-1)/2}\frac{\zeta(n)}{(2\pi\ii)^n}$, see Appendix 1.

\bip

\section{Braid configuration spaces}

\sip

{\bf 5.1. Compactified braid configuration spaces as operads $\caL ie_\infty$ and $\cM or(\caL ie_\infty)$.}
In the previous sections we studied configurations of points, $C_n$ and $C_{n,0}$, which were both ordered
and numbered. Let $B_n$ and $B_{n,0}$ be their versions in which the total ordering of points is forgotten. Their compactifications, $\overline{B}_n$ and $\widehat{B}_{n,0}$ can be defined
in a full analogy to the previous case, and can be characterized as follows.

\mip

{\bf 5.1.1. Proposition} \cite{GJ}. {\em The face complex, $C^\bullet(\overline{B})$,
 of the family of
compactified braid configurations spaces,
$\{\overline{B}_n\}_{n\geq 2}$, has a structure of a dg operad canonically isomorphic
to the operad, $\caL ie_\infty$, of strong homotopy Lie algebras.}

\mip

{\bf 5.1.2. Proposition}. {\em The face complex of $\overline{B}\sqcup \widehat{B}\sqcup \overline{B}$
has structure of  a dg 2-coloured operad
which is canonically isomorphic
to the dg 2-coloured  operad $\cM or(\caL ie_\infty)$  describing pairs of
  $\caL ie_\infty$-algebras and
 $\caL ie_\infty$-morphisms between them.}

\mip

{\bf 5.2. De Rham field theories on braid configuration spaces.}
Let $\cB_{n,l}$ stand for a family of graphs, $\{\Ga\}$,
such that (i) $\# V(\Ga)=n$, (ii) $\# E(\Ga)=l$, (iii) $\Ga$  has no loop type edges, and (iv) the set $E(\Ga)$ is totally ordered (up to an even permutation).
Let $\fB_{n,l}$ be the version of $\cB_{n,l}$ with data (iv) forgotten.
 Note that edges of these graphs
are {\em not}\, directed, and their vertices  are {\em not}\, numbered.

\sip

With any graph $\Ga\in \cB_{n,l}$ one can associate a linear map
$$
\Ba{rccc}
\Phi^s_\Ga: &\odot^n \cT_{poly}(\R^d) &\lon &\cT_{poly}(\R^d)[-l]\\
& \ga_1\ot\ldots \ot\ga_n & \lon& \Phi^s_\Ga(\ga_1,\ldots,\ga_n),
\Ea
$$
$$
\Phi^s_\Ga(\ga_1,\ldots,\ga_n)
:=\frac{1}{n!}\sum_{f: Vert(\Ga)\rar [n]}
\left[\left(\prod_{e\in Edges(\Ga)}\Delta_{e}\right)
\prod_{v\in Vert(\Ga)} \ga_{f(v)}(\psi_{f(v)},x_{f(v)})
\right]_{\hspace{-2mm}x_1=\ldots=x_n\atop
\psi_1=\ldots=\psi_n},
$$
where the operator $\Delta_e$ corresponding to an edge $e$ beginning at a vertex
labelled by $i\in [n]$ and ending at a vertex labelled by $j\in [n]$ is given by
$$
\Delta_e:= \sum_{a=1}^d\left(\frac{\p^2}{\p x_{(i)}^a\p \psi_{(j)\,a}}
+  \frac{\p^2}{\p x_{(j)}^a\p \psi_{(i)\,a}}
\right)
.
$$

\sip

{\em De Rham field theories}\, on $\overline{B}$ and $\widehat{B}$  can be defined in
a full analogy with the ones on  $\overline{C}$ and $\widehat{C}$.  Moreover, Theorems 4.2.1 and 4.3.1 hold true
with symbols  $\overline{C}$, $\widehat{C}$  and  $\caL eib_\infty$ replaced, respectively,
 by  $\overline{B}$, $\widehat{B}$  and $\caL ie_\infty$.

\sip

A class of de Rham theories  on $\overline{C}\sqcup\widehat{C}\sqcup \overline{C}$
determined  by a propagator $\om(z_i,z_j)$ on $\widehat{C}_{2,0}$ satisfying the symmetry condition
(cf.\ \S 4.5),
$$
\om(z_i,z_j)=\om(z_j,z_i),
$$
 comes in fact from a class of
 de Rham field theories on $\overline{B}\sqcup\widehat{B}\sqcup \overline{B}$.
It is easy to see that for any choice of such a symmetric propagator
 the weight,
$C_{\mathfrak w_n}(s)$, of the following graph
$$
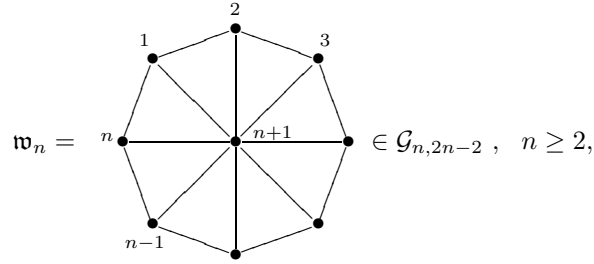
\mathfrak w_n=\ \ \xy
(5,0.7)*{^{n+1}},
(-12,13)*{^{1}},
(0,17)*{^{2}},
(12,13)*{^{3}},
(-17,0)*{^{n}},
(-12,-14)*{^{n-1}},
 (0,0)*{\bullet}="1",
(-11,11)*{\bullet}="2",
(0,15)*{\bullet}="3",
(11,11)*{\bullet}="4",
(15,0)*{\bullet}="5",
(11,-11)*{\bullet}="6",
(0,-15)*{\bullet}="7",
(-15,0)*{\bullet}="n",
(-11,-11)*{\bullet}="n-1",
\ar @{-} "1";"2" <0pt>
\ar @{-} "1";"3" <0pt>
\ar @{-} "1";"4" <0pt>
\ar @{-} "1";"5" <0pt>
\ar @{-} "1";"6" <0pt>
\ar @{-} "1";"7" <0pt>
\ar @{-} "1";"n" <0pt>
\ar @{-} "1";"n-1" <0pt>
\ar @{-} "2";"3" <0pt>
\ar @{-} "3";"4" <0pt>
\ar @{-} "4";"5" <0pt>
\ar @{-} "5";"6" <0pt>
\ar @{-} "6";"7" <0pt>
\ar @{-} "7";"n-1" <0pt>
\ar @{-} "n-1";"n" <0pt>
\ar @{-} "n";"2" <0pt>
\endxy
\ \in \cG_{n,2n-2} \ , \ \
n\geq 2,
$$
is equal to zero for {\em even}\, $n$. For odd $n$ its weight is, in general, non-zero.
 The {\em infinitesimal}\,
part, $\delta\al$,  of the associated exotic transformation of an ordinary Poisson structure,
 $$
 F_{\frac{1}{2}K}^{sym}(\al)=\al + \frac{C_{\mathfrak w_3}}{\# Aut(w_3)} \Phi^s_\Ga(\al,\al,\al,\al) + \mbox{higher order (in $\al$) terms},
 $$
 is controlled by the graph $\mathfrak w_3$ (which is the same as the
 tetrahedron graph
 $ \Ba{c}
\xy
 (0,0)*{\bullet}="a",
(12,0)*{\bullet}="b",
(6,9)*{\bullet}="c",
(14,4)*{\bullet}="d",
\ar @{-} "a";"b" <0pt>
\ar @{-} "a";"c" <0pt>
\ar @{-} "b";"c" <0pt>
\ar @{-} "d";"c" <0pt>
\ar @{-} "b";"d" <0pt>
\ar @{.} "d";"a" <0pt>
\endxy
\Ea
 $) so that we get, up to a non-zero numerical factor, the following infinitesimal change of the Poisson structure,
$$
\delta \al\sim 
\sum_{
i,j,k,l,m\atop k',l',m'}\hspace{-2mm}\left(
\frac{\p^3 \al^{ij}}{\p x^k\p x^l\p x^m}
\frac{\p \al^{k k'}}{\p x^{l'}}
\frac{\p \al^{l l'}}{\p x^{m'}}
\frac{\p \al^{m m'}}{\p x^{k'}}
+
\frac{4}{3}
\frac{\p^3 \al^{im}}{\p x^k\p x^l}
\frac{\p \al^{k k'}}{\p x^{l'}}
\frac{\p \al^{l l'}}{\p x^{m'}}
\frac{\p \al^{j m'}}{\p x^{k'}\p x^m}
\right)(\p_i\wedge \p_j).
$$
According to Kontsevich \cite{Ko2}, the second term  vanishes  identically
for any Poisson structure $\al$  so that the above relation simplifies further,
$$
\delta\al \sim 
\sum_{
i,j,k,l,m,k',l',m'}
\frac{\p^3 \al^{ij}}{\p x^k\p x^l\p x^m}
\frac{\p \al^{k k'}}{\p x^{l'}}
\frac{\p \al^{l l'}}{\p x^{m'}}
\frac{\p \al^{m m'}}{\p x^{k'}}(\p_i\wedge \p_j).
$$
Thus the flow of Poisson structures,
$$
\frac{d\al}{d t}
:=\delta \al,
$$
associated with the infinitesimal part of the exotic $\caL ie_\infty$ morphism $F_{\frac{1}{2}K}^{sym}$
is precisely the one which was found by Kontsevich long ago in  \S 4.6.3 of \cite{Ko}
as an example of an exotic (i.e.\ homotopy non-trivial, see Remark 4.2.9) infinitesimal
$\caL ie_\infty$ automorphism
of the Schouten algebra.

\section{Towards a new differential geometry}

The classical architecture  of geometry and
 theoretical\ physics can be described as follows:
a geometric  structure is a function (a ``field" or an ``observable")
 on a manifold (``space-time") satisfying some differential
equations
$$
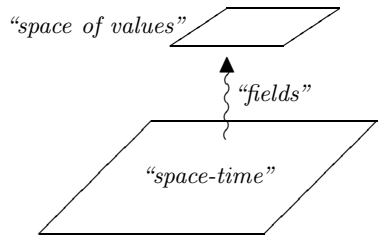
\begin{xy}
 <0mm,0mm>*{};<30mm,0mm>*{}**@{-},
 <0mm,0mm>*{};<15mm,15mm>*{}**@{-},
 <15mm,15mm>*{};<45mm,15mm>*{}**@{-},
 <0mm,0mm>*{};<20mm,0mm>*{}**@{-},
 <30mm,0mm>*{};<45mm,15mm>*{}**@{-},
  <22.5mm,7.5mm>*{\mbox{\small\em ``space-time"}};
 <25mm,12.5mm>*{};<25mm,22.5mm>*{}**@{~},
 <25mm,22.5mm>*{\blacktriangle};
<31mm,18.5mm>*{\mbox{\small\em ``fields"}};
<17.5mm,25mm>*{};<32.5mm,25mm>*{}**@{-},
 <17.5mm,25mm>*{};<25mm,30mm>*{}**@{-},
 <25mm,30mm>*{};<40mm,30mm>*{}**@{-},
 <32.5mm,25mm>*{};<40mm,30mm>*{}**@{-},
<7.5mm,27.5mm>*{\mbox{\small\em ``space of values"}};
\end{xy}
$$
The theory of (wheeled) props offers a different picture \cite{Me2,Me-lec}
in which
``space-time" equipped with a geometric structure is itself a function (a representation)
on a more fundamental object --- a certain dg free prop, a kind of a graph complex,
$$
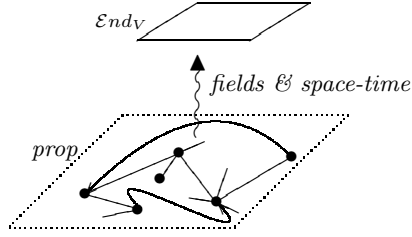
\begin{xy}
 <0mm,0mm>*{};<30mm,0mm>*{}**@{.},
 <0mm,0mm>*{};<15mm,15mm>*{}**@{.},
 <15mm,15mm>*{};<45mm,15mm>*{}**@{.},
 <30mm,0mm>*{};<45mm,15mm>*{}**@{.},
  <6mm,10mm>*{\mbox{\small \em prop}};
 <25mm,12mm>*{};<25mm,22mm>*{}**@{~},
 <25mm,22mm>*{\blacktriangle};
<40mm,19mm>*{\mbox{{\small\em  fields \& space-time}}};
<17mm,25mm>*{};<32mm,25mm>*{}**@{-},
 <17mm,25mm>*{};<25mm,30mm>*{}**@{-},
 <25mm,30mm>*{};<40mm,30mm>*{}**@{-},
 <32mm,25mm>*{};<40mm,30mm>*{}**@{-},
<15mm,27mm>*{^{\cE nd_{V}}};
<17mm,2.5mm>*{\bullet};
<27.5mm,3.5mm>*{\bullet};
<20mm,6.5mm>*{\bullet};
<37.5mm,9.5mm>*{\bullet};
<10mm,4.5mm>*{\bullet};
<22.5mm,10mm>*{\bullet};
<37.5mm,9.5mm>*{};<27.5mm,3.5mm>*{}**@{-},
<30.5mm,2mm>*{};<27.5mm,3.5mm>*{}**@{-},
<22.5mm,10mm>*{};<26mm,11.5mm>*{}**@{-},
<27.5mm,3.5mm>*{};<24mm,4mm>*{}**@{-},
<27.5mm,3.5mm>*{};<29mm,8mm>*{}**@{-},
<27.5mm,3.5mm>*{};<22.5mm,10mm>*{}**@{-},
<17.5mm,2.5mm>*{};<10mm,4.5mm>*{}**@{-},
<17.5mm,2.5mm>*{};<12.5mm,1.5mm>*{}**@{-},
<22.5mm,10mm>*{};<10mm,4.5mm>*{}**@{-},
<20mm,6.5mm>*{};<22.5mm,10mm>*{}**@{-},
 (17.5,2.5)*{}
   \ar@{->}@(ul,ur) (37.5,9.5)*{};
 (10.0,4.5)*{}
   \ar@{->}@(ul,dr) (27.5,3.5)*{}
\end{xy}
$$
Here $V$ stands for a vector space  modelling a local coordinate chart of some (say, real analytic)
manifold $M$. A real analytic gluing map, $V\rar W$, of two coordinate charts on $M$ can be understood as a morphism of props $\cE nd_V\rar \cE nd_W$. Then a
{\em consistent}\, gluing of local
geometric structures, $\{\cP\rar \cE nd_V\}$,
 controlled by some prop $\cP$ into a global one on the manifold $M$ can be understood as commutativity
of diagrams of the form,
$$
\xy
(0,0)*{\cP}="a",
(-10,10)*{\cE nd_V\ }="b_1",
(10,10)*{\ \cE nd_W}="b_2",
\ar @{->} "a";"b_1" <4pt>
\ar @{->} "a";"b_2" <-4pt>
\ar @{->} "b_1";"b_2" <0pt>
\endxy
$$
If, however, the dg prop $\cP$ admits a non-trivial group
of automorphisms, then the above gluing pattern can be replaced by the following one,
$$
\xymatrix{
\cE nd_V\ar[r] & \cE nd_W \\
\cP\ \ar[r]^{f_{VW}\in Aut(\cP)}\ar[u] &\ \cP\ar[u]
}
$$
i.e.\ the group $Aut(\cP)$ is allowed to twist the standard local coordinate gluing mappings.
Note that the group $Aut(\cP)$ is universal and does not depend on a particular manifold $M$,
i.e.\ one modifies in this way the whole category of geometric structures of type $\cP$.

\sip

More concretely, one can think of our explicit formulae for exotic $\caL ie_\infty$
automorphisms of the Schouten algebra in $\R^d$ as of gluing mappings defining a {\em quantum
manifold}\, out of local coordinate charts,
$$
\mathfrak U=(\odot^\bu(\cT_{poly}(\R^d))[\hbar]], \Delta, d_{S}),$$
which, by definition, are dg coalgebras $\odot^\bu(\cT_{poly}(\R^d))[\hbar]]$ equipped with the standard
comultiplication $\Delta$ and with the co-differential $d_{S}$ corresponding to the
Schouten bracket in $\cT_{poly}(\R^d)$. Every ordinary manifold is a special case of such a quantum manifold in which all gluing mapping are  homotopy trivial morphisms
of coordinate charts  $f_{ij}:\mathfrak U_i\rar \mathfrak U_j$ given by the ordinary changes of coordinates. In general, however, such an ordinary manifold is only an $\hbar\rar 0$ limit of a quantum one, and the GT group might twist the gluing mappings
$f_{ij}$
in a homotopy non-trivial way. Such quantum manifolds can be useful in the study
of {\em geometric}\, invariants of the Grothendieck-Teichmueller group; the nature, perhaps, does not depend on a particular choice of the quantization scheme made by an observer so that quantum observables must be GT invariants.

\sip

Poisson structures on any finite-dimensional manifold are controlled by a dg wheeled prop $\cP oly$
(see Appendix 4)
whose automorphism group is very non-trivial. It is worth noting that another important class
of geometric structures --- the so called Nijenhuis structures (coming from the famous Nijenhuis
integrability condition for an almost complex structure) --- are also controlled by a certain dg (wheeled) prop
(see \cite{Me2,St}).

\newpage

\begin{center}
{\sc Appendix 1: Wheels and zeta function}
\end{center}

{\bf Theorem A.} {\em The weight,
$$
C_{W_n}:=\frac{1}{(2\pi)^{2n}}
\int_{C_{n+1,0}} \om_{\frac{1}{2}}(z_{n+1},z_1)\wedge\ldots \wedge\om_{\frac{1}{2}}(z_{n+1},z_n)
\wedge \om_{\frac{1}{2}}(z_1,z_2)\wedge\om_{\frac{1}{2}}(z_2,z_3)\wedge\ldots\wedge
\om_{\frac{1}{2}}(z_n,z_1)
$$
of the following graph with $n+1$ vertices,
$$
W_n=\ \ 
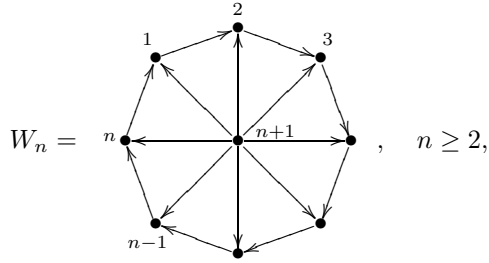
\xy
(5,0.7)*{^{n+1}},
(-12,13)*{^{1}},
(0,17)*{^{2}},
(12,13)*{^{3}},
(-17,0)*{^{n}},
(-12,-14)*{^{n-1}},
 (0,0)*{\bullet}="1",
(-11,11)*{\bullet}="2",
(0,15)*{\bullet}="3",
(11,11)*{\bullet}="4",
(15,0)*{\bullet}="5",
(11,-11)*{\bullet}="6",
(0,-15)*{\bullet}="7",
(-15,0)*{\bullet}="n",
(-11,-11)*{\bullet}="n-1",
\ar @{->} "1";"2" <0pt>
\ar @{->} "1";"3" <0pt>
\ar @{->} "1";"4" <0pt>
\ar @{->} "1";"5" <0pt>
\ar @{->} "1";"6" <0pt>
\ar @{->} "1";"7" <0pt>
\ar @{->} "1";"n" <0pt>
\ar @{->} "1";"n-1" <0pt>
\ar @{->} "2";"3" <0pt>
\ar @{->} "3";"4" <0pt>
\ar @{->} "4";"5" <0pt>
\ar @{->} "5";"6" <0pt>
\ar @{->} "6";"7" <0pt>
\ar @{->} "7";"n-1" <0pt>
\ar @{->} "n-1";"n" <0pt>
\ar @{->} "n";"2" <0pt>
\endxy
\ \ , \ \ \
n\geq 2,
$$
with respect to  Kontsevich's 1/2-propagator
$\om_{\frac{1}{2}}(z_i,z_j):= \frac{1}{ \ii}d\log \frac{z_1-z_2}{\overline{z}_1-z_2}$
is given by
$$
C_{W_n}=(-1)^{n(n-1)/2}\frac{\zeta(n)}{(2\pi\ii)^n}=(-1)^{n(n-1)/2}
\frac{\sum_{p=1}^{\infty}\frac{1}{p^n}}{(2\pi\ii)^n}.
$$
}

\Proof We identify $C_{n+1,0}$ with a subspace of $\Conf_{n+1,0}$
consisting of all configurations, $\{z_1, \ldots, z_n,z_{n+1}\}$, with $z_{n+1}=\ii$, and introduce
in $C_{n+1,0}$ a system of coordinates,
$\{\rho_i, \phi_i\ \mid 0<\rho_i<1, 0<\phi_i\leq 2\pi\}_{1\leq i\leq n}$,
as follows
$$
\frac{z_i - \ii}{z_i+\ii}=: \rho_ie^{\ii\phi_i}.
$$
Thus $\rho_i=\tanh (\frac{1}{2}h(\ii,z_i))$, where $h(\ii,z_i)$ is the hyperbolic distance from $\ii$ to $z_i\in \bbH$
and $\phi_i$ is the angle between the vertical line and the hyperbolic
geodesic from $\ii$ to $z_i$.

\sip
Let $I$ be the ideal in the de Rham algebra on $C_{n,0}$ generated by 1-forms
$d(\rho_ie^{\ii\phi_i})$, $1\leq i\leq n$.
As
$$
z_i= \ii\frac{1+\rho_ie^{\ii\phi_i}}{1-\rho_ie^{\ii\phi_i}}
$$
we have
\Beqrn
\frac{1}{\ii}d \log\frac{z_i-z_j}{\overline{z}_i -z_j} &=&
 \frac{1}{\ii}d \log\frac{\left(1-\rho_ie^{-\ii\phi_i}\right)
 \left(\rho_je^{\ii\phi_j} - \rho_ie^{\ii\phi_i}
 \right)}
 {\left(1-\rho_ie^{\ii\phi_i}\right)
 \left(1-\rho_i\rho_je^{\ii(\phi_j-\phi_i)}
 \right)}\\
 &=&\frac{1}{\ii}\left(\frac{\rho_je^{\ii\phi_j}}{1-\rho_i\rho_je^{\ii(\phi_j-\phi_i)}}
 -\frac{1}{1 - \rho_ie^{-\ii\phi_i}}\right)
 d\left(\rho_ie^{-\ii\phi_i}
 \right)\bmod I.
\Eeqrn
Hence, modulo $I$,
\Beqrn
\frac{1}{\ii}d\log(\rho_ie^{\ii\phi_i})\wedge \frac{1}{\ii}d
\log\frac{z_i-z_j}{\overline{z}_i -z_j}
&=&\frac{2}{\ii} d\phi_i\wedge d\rho_i
\left(\frac{\rho_j e^{\ii(\phi_j-\phi_i)}}{1-\rho_i\rho_je^{\ii(\phi_j-\phi_i)}}
-\frac{e^{-\ii\phi_i}}{1 - \rho_ie^{-\ii\phi_i}}
\right)\\
&=&\frac{2}{\ii} d\phi_i\wedge d\rho_i
\left(\sum_{k_i=0}^\infty\rho_i^{k_i}\rho_j^{k_i+1}e^{\ii(k_i+1)(\phi_j-\phi_i)}
- \sum_{k_i=0}^\infty\rho_i^{k_i}e^{-\ii(k_i+1)\phi_i}
\right).
\Eeqrn
As
$$
\int e^{\ii k\phi}d\phi= \left\{
\Ba{cc}
2\pi & \mbox{if}\ k=0\\
0    & \mbox{otherwise},
\Ea
\right.
$$
we finally get (identifying $k_{-1}$ with $k_{i_n}$)
\Beqrn
C_{W_n} &=&\frac{2^n(-1)^{n(n-1)/2}}{(\ii)^n (2\pi)^{2n}}\int d\phi_1...d\phi_n d\rho_1...d\rho_n
\sum_{k_1,...,k_n=0}^\infty\prod_{i=1}^n \rho_i^{k_{i-1} + k_{i+1}+1}e^{\ii (k_{i-1}-k_i)\phi_i}\\
&=&\frac{2^n(-1)^{n(n-1)/2}}{(\ii)^n (2\pi)^{n}}\int_0^1...\int_0^1 d\rho_1...d\rho_n
\sum_{k=0}^\infty (\rho_1\ldots \rho_n)^{2k+1}\\
&=&\frac{2^n(-1)^{n(n-1)/2}}{(\ii)^n (2\pi)^{n}}\int_0^1...\int_0^1 \frac{ \rho_1\ldots \rho_n d\rho_1...d\rho_n}
{1-(\rho_1\ldots \rho_n)^{2}}\\
&=&\frac{(-1)^{n(n-1)/2}}{(\ii)^n (2\pi)^{n}}\int_0^1...\int_0^1 \frac{ dx_1...dx_n}
{1- x_1... x_n}\\
&=&(-1)^{n(n-1)/2}\frac{\zeta(n)}{(\ii)^n (2\pi)^{n}}.
\Eeqrn
\hfill $\Box$
\mip

Note that $\zeta(2n)=\frac{(-1)^{n+1}(2\pi)^{2n}B_{2n}}{2(2n)!}$ so that for $n$ even
$$
C_{W_{n}}=-(-1)^{n(n-1)/2}\frac{1}{2n!}B_{n}=-(-1)^{n(n-1)/2}n \cB_n.
$$

\bip

\begin{center}
{\sc Appendix 2: Example of a boundary strata}
\end{center}
Consider a pair
$$
G=\hspace{-4mm}\Ba{c}
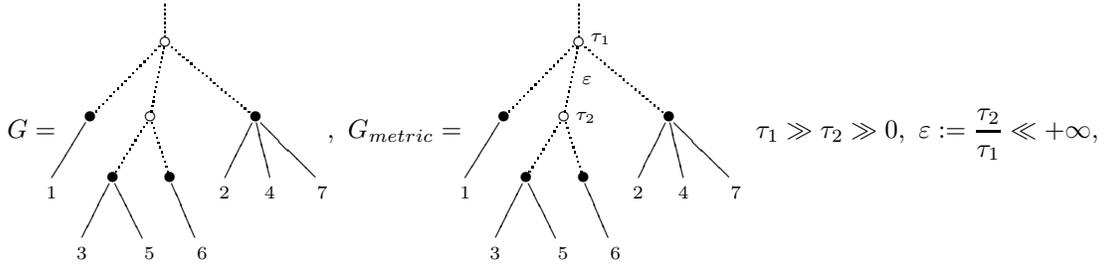
\xy
(-15,-10)*{_1},
(-11,-18)*{_3},
(-2,-18)*{_5},
(5,-18)*{_6},
(8,-10)*{_2},
(14,-10)*{_4},
(21,-10)*{_7},
(0,15)*{}="0",
 (0,10)*{\circ}="a",
(-10,0)*{\bullet}="b_1",
(-2,0)*{\circ}="b_2",
(12,0)*{\bullet}="b_3",
(-15,-8)*{}="c_0",
(0.6,-8)*{\bullet}="c_1",
(-7,-8)*{\bullet}="c_2",
(8,-8)*{}="c_3",
(14,-8)*{}="c_4",
(20,-8)*{}="c_5",
(-11,-16)*{}="d_1",
(-3,-16)*{}="d_2",
(4,-16)*{}="d_3",
\ar @{.} "a";"0" <0pt>
\ar @{.} "a";"b_1" <0pt>
\ar @{.} "a";"b_2" <0pt>
\ar @{.} "a";"b_3" <0pt>
\ar @{-} "b_1";"c_0" <0pt>
\ar @{.} "b_2";"c_1" <0pt>
\ar @{.} "b_2";"c_2" <0pt>
\ar @{-} "b_3";"c_3" <0pt>
\ar @{-} "b_3";"c_4" <0pt>
\ar @{-} "b_3";"c_5" <0pt>
\ar @{-} "c_2";"d_1" <0pt>
\ar @{-} "c_2";"d_2" <0pt>
\ar @{-} "c_1";"d_3" <0pt>
\endxy
\Ea \hspace{-2mm}, \
G_{metric}=\hspace{-3mm}
\Ba{c}
\xy
(-15,-10)*{_1},
(-11,-18)*{_3},
(-2,-18)*{_5},
(5,-18)*{_6},
(8,-10)*{_2},
(14,-10)*{_4},
(21,-10)*{_7},
(1.2,5)*{_{\var}},
(3,10)*{^{\tau_1}},
(1,0)*{_{\tau_2}},
(0,15)*{}="0",
 (0,10)*{\circ}="a",
(-10,0)*{\bullet}="b_1",
(-2,0)*{\circ}="b_2",
(12,0)*{\bullet}="b_3",
(-15,-8)*{}="c_0",
(0.6,-8)*{\bullet}="c_1",
(-7,-8)*{\bullet}="c_2",
(8,-8)*{}="c_3",
(14,-8)*{}="c_4",
(20,-8)*{}="c_5",
(-11,-16)*{}="d_1",
(-3,-16)*{}="d_2",
(4,-16)*{}="d_3",
\ar @{.} "a";"0" <0pt>
\ar @{.} "a";"b_1" <0pt>
\ar @{.} "a";"b_2" <0pt>
\ar @{.} "a";"b_3" <0pt>
\ar @{-} "b_1";"c_0" <0pt>
\ar @{.} "b_2";"c_1" <0pt>
\ar @{.} "b_2";"c_2" <0pt>
\ar @{-} "b_3";"c_3" <0pt>
\ar @{-} "b_3";"c_4" <0pt>
\ar @{-} "b_3";"c_5" <0pt>
\ar @{-} "c_2";"d_1" <0pt>
\ar @{-} "c_2";"d_2" <0pt>
\ar @{-} "c_1";"d_3" <0pt>
\endxy
\Ea  \tau_1\gg\tau_2\gg 0,\ \var:=\frac{\tau_2}{\tau_1}\ll +\infty,
$$
consisting of a graph $G$ from the face complex $C_\bullet(\widehat{C}_{7,0})$
 and an associated metric graph. The latter defines
a smooth coordinate chart,
$$
U(G)\simeq {C}_3\times C_2 \times \widehat{C}_{3,0}\times \widehat{C}_{2,0}\times \R\times \R,
$$
 near the face $G\subset \widehat{C}_{7,0}$
whose intersection with $C_{7,0}$ consists, by definition, of all those configurations, $p$,
  of 7 points in $\bbH$
 which result from the following four step construction:
 \Bi
\item[{\em Step 1}:] take an arbitrary hyperpositioned configuration, $p_0^{(1)}\in C_3$, of 3 points
labelled by $1$, and, say,
$z'$ and $z''$, and  magnify
it,  $p_0^{(1)}\rar \tau_1\cdot p_0^{(1)}$;
\item[{\em Step 2}:] take an arbitrary hyperpositioned configuration, $p_0^{(2)}\in C_2$, of 2 points
labelled by $6$, and, say,
$z'''$, magnify it,
  $p_0^{(2)}\rar \tau_2\cdot p_0^{(2)}$, and place the result at the position $z'$;
\item[{\em Step 3:}] take an arbitrary hyperpositioned
configuration of $3$ points,
$p^{(3)}\in C_{3,0}$, labelled  by  $\{2,4,7\}$ and place it at the position  $z''$;
\item[{\em Step 4}:] take an arbitrary hyperpositioned configuration, $p^{(4)}\in C_{2,0}$, of 2 points
labelled by $3$ and $5$ and place it at the position $z'''$.
\Ei
The final result is a hyperpositioned point $p=(z_1, z_3,z_5, z_6,z_2,z_4,z_7)$
in $C_{7,0}$ of the form
$$
 \Ba{c}
{
\xy
(-13,39)*{\bullet},
(-10,39)*{^{z_4}},
(-12.7,43.5)*{\bullet},
(-14,45)*{^{z_7}},
(10.0,26)*{\bullet}="z_3",
(8,26)*{_{z_3}},
(11.7,29.8)*{\bullet},
(14,31)*{_{z_5}},
(21.8,30)*{\bullet},
(25,30)*{_{z_6}},
(2,8)*{{\bullet}}="z_1",
(4.8,6.1)*{_{z_1}},
(-2.4,19.3)*{_{\times}},
(-4,18.7)*{_{z'}},
(4.4,13)*{_{\times}},
(7,12.7)*{_{z''}},
 (10.9,28.1)*{
\xycircle(2,2){.}};
(16.0,28)*{
\xycircle(6,6){.}};
%
%
 (-11,39.7)*{
\xycircle(4,4){.}};
(-9,36)*{\bullet}="z_2",
(-6.3,34.7)*{_{z_2}},
 (-30,0)*{}="a",
(40,0)*{}="b",
(3,0)*{}="c",
(3,50)*{}="d",
\ar @{->} "a";"b" <0pt>
\ar @{->} "c";"d" <0pt>
\ar @{.} "z_1";"z_2" <0pt>
\ar @{.} "z_1";"z_3" <0pt>
\endxy}
\Ea
$$

An equivalent coordinate chart  results from an alternative 5-step construction which uses
instead of the parameters $(\tau_1,\tau_2)$ another pair of independent parameters, $\tau_1$ and $\var$:
\Bi
\item[{\em Step 1}:] take an arbitrary hyperpositioned configuration, $p_0^{(1)}\in C_3$, of 3 points
labelled by $1$, and, say,
$z'$ and $z''$;
\item[{\em Step 2}:] take an arbitrary hyperpositioned configuration, $p_0^{(2)}\in C_2$, of 2 points
labelled by $6$, and, say,
$z'''$, $\var$-magnify it,
  $p_0^{(2)}\rar \var \cdot p_0^{(2)}$, and place the result at the position $z'$;
\item[{\em Step 3:}] magnify the resulting configuration of points $(1,6,z'',z''')$ by the factor $\tau_1$;
\item[{\em Step 4:}] take an arbitrary hyperpositioned
configuration of $3$ points,
$p^{(3)}\in C_{3,0}$, labelled  by  $\{2,4,7\}$ and place it at the position  $z''$;
\item[{\em Step 5}:] take an arbitrary hyperpositioned configuration, $p^{(4)}\in C_{2,0}$, of 2 points
labelled by $3$ and $5$ and place it at the position $z'''$.
\Ei
The $10$-dimensional face $G$ lies in the intersection of two $11$-dimensional
faces described by the following graphs,
$$
G_1=\hspace{-4mm}\Ba{c}
\xy
(-15,-10)*{_1},
(-7,-10)*{_3},
(0,-10)*{_5},
(4,-10)*{_6},
(8,-10)*{_2},
(14,-10)*{_4},
(21,-10)*{_7},
(0,15)*{}="0",
 (0,10)*{\circ}="a",
(-10,0)*{\bullet}="b_1",
(-3,0)*{\bullet}="b_2",
(3,0)*{\bullet}="b_3",
(12,0)*{\bullet}="b_4",
(-15,-8)*{}="c_0",
(4,-8)*{}="c_1",
(-7,-8)*{}="c_2",
(8,-8)*{}="c_3",
(14,-8)*{}="c_4",
(20,-8)*{}="c_5",
(-1,-8)*{}="c_6",
\ar @{.} "a";"0" <0pt>
\ar @{.} "a";"b_1" <0pt>
\ar @{.} "a";"b_2" <0pt>
\ar @{.} "a";"b_3" <0pt>
\ar @{.} "a";"b_4" <0pt>
\ar @{-} "b_1";"c_0" <0pt>
\ar @{-} "b_3";"c_1" <0pt>
\ar @{-} "b_2";"c_2" <0pt>
\ar @{-} "b_2";"c_6" <0pt>
\ar @{-} "b_4";"c_3" <0pt>
\ar @{-} "b_4";"c_4" <0pt>
\ar @{-} "b_4";"c_5" <0pt>
\endxy
\Ea
,
\ \ \ \
G_2=\hspace{-4mm}\Ba{c}
\xy
(-15,-10)*{_1},
(-7,-10)*{_3},
(-2,-10)*{_5},
(3,-10)*{_6},
(8,-10)*{_2},
(14,-10)*{_4},
(21,-10)*{_7},
(0,15)*{}="0",
 (0,10)*{\circ}="a",
(-10,0)*{\bullet}="b_1",
(-2,0)*{\bullet}="b_2",
(12,0)*{\bullet}="b_3",
(-15,-8)*{}="c_0",
(2,-8)*{}="c_1",
(-7,-8)*{}="c_2",
(8,-8)*{}="c_3",
(14,-8)*{}="c_4",
(20,-8)*{}="c_5",
(-3,-8)*{}="c_6",
\ar @{.} "a";"0" <0pt>
\ar @{.} "a";"b_1" <0pt>
\ar @{.} "a";"b_2" <0pt>
\ar @{.} "a";"b_3" <0pt>
\ar @{-} "b_1";"c_0" <0pt>
\ar @{-} "b_2";"c_1" <0pt>
\ar @{-} "b_2";"c_2" <0pt>
\ar @{-} "b_2";"c_6" <0pt>
\ar @{-} "b_3";"c_3" <0pt>
\ar @{-} "b_3";"c_4" <0pt>
\ar @{-} "b_3";"c_5" <0pt>
\endxy
\Ea
$$
and the above two constructions give us explicit descriptions
of the associated embeddings, $G\hook G_1$ and $G\hook G_2$.

\bip

\begin{center}
{\sc Appendix 3: $\caL eib_\infty$ automorphisms of Maurer-Cartan sets}
\end{center}

\mip
{\bf Proposition}. {\em Any $\caL eib_\infty$-automorphism,
$$
\left\{F_n:\ot \fg \rar\fg[2-2n]\right\}_{n\geq 1},
$$
 of a Lie algebra $(\fg,\, [\ ,\ ]: \odot^2\fg \rar \fg[-1])$
induces an automorphism of its set,
$$
\left\{\al\in \fg[[\hbar]]: [\al,\al]=0\ \& \ |\al|=2\right\},
$$ of Maurer-Cartan elements by
the formula}
$$
 \al \rar F^{Leib}(\al):=\sum_{n\geq 1}\frac{\hbar^{n-1}}{n!}F^{Leib}_n(\al^{\ot n})
$$

\Proof It follows from (\ref{Ch3: d on black corollas}) that, for $n\geq 2$,
\Beqrn
0 &=& \sum_{[n]=B_1\sqcup B_2 \atop \inf B_1< \inf B_2}
[F_{\# B_1}(\al^{\ot \#{B_1}}), F_{\# B_2}(\al^{\ot \#{B_2}}
)]\\
&=&
\sum_{[n]\setminus [1]=S_1\sqcup S_2
\atop \# S_1\geq 0, \# S_2\geq 1}
[F_{\# S_1+1}(\al^{\ot (\#{S_1}+1)}), F_{\# S_2}(\al^{\ot \#{S_2}}
)]\\
&=&
\sum_{k=0}^{n-2}\frac{(n-1)!}{k!(n-k-1)!}
[F_{k+1}(\al^{\ot (k+1)}), F_{n-k-1}(\al^{\ot (n-k-1)})
]\\
&=&\frac{1}{2}
\sum_{k=0}^{n-2}\left(\frac{(n-1)!}{k!(n-k-1)!}+ \frac{(n-1)!}{(n-k-2)!(k+1)!}\right)
[F_{k+1}(\al^{\ot (k+1)}), F_{n-k-1}(\al^{\ot (n-k-1)})
]\\
&=&
\frac{1}{2}
\sum_{k=0}^{n-2}\frac{n!}{(k+1)!(n-k-1)!}
[F_{k+1}(\al^{\ot (k+1)}), F_{n-k-1}(\al^{\ot (n-k-1)}).
]\\
&=&
\frac{n!}{2}
\sum_{n=p+q\atop p,q\geq 1}\frac{1}{p!q!}
[F_{p}(\al^{\ot p}), F_{q}(\al^{\ot q})
].
\Eeqrn
Then
$$
\hbar^2[F^{Leib}(\al),F^{Leib}(\al)]=\sum_{n\geq 2}\hbar^{n}\sum_{n=p+q\atop p,q\geq 1}\frac{1}{p!q!}
[F_{p}(\al^{\ot p}), F_{q}(\al^{\ot q})
]
=0.
$$
\hfill $\Box$

\bip
\begin{center}
{\sc Appendix 4: Weights of  4-vertex graphs in a de Rham field theory on $\overline{C}_\bu$}
\end{center}

\mip

{\bf A4.1. Lemma}. {\em For any fixed points $z_1,z_2\in \C$ with $\Im z_1 < \Im z_2$ the integral
$$
\int_{\Img z_1 <\Img z< \Img z_2} dArg(z-z_1)\wedge dArg (z-z_2)
$$
vanishes.}

\begin{proof}
Let $\circ$ be the middle point of the interval connecting $z_1$ to $z_2$,
$$
{\Ba{c}
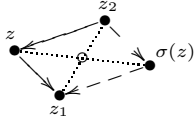
\xy
<0mm,-2.0mm>*{_{z_1}},
<-6.5mm,7.5mm>*{^{z}},
<6.5mm,11.5mm>*{^{z_2}},
<15.5mm,5.5mm>*{^{\sigma(z)}},
<3mm,5mm>*{\circ},
(0,0)*{\bullet}="1",
(-6,6)*{\bullet}="2",
(6,10)*{\bullet}="3",
(12,4)*{\bullet}="2'",
\ar @{->} "3";"2" <0pt>
\ar @{->} "2";"1" <0pt>
\ar @{.} "3";"1" <0pt>
\ar @{-->} "3";"2'" <0pt>
\ar @{-->} "2'";"1" <0pt>
\ar @{.} "2";"2'" <0pt>
\endxy\Ea}
$$
and let $\sigma: \C \rar \C$ be the reflection at $\circ$. This map preserves the subspace
$\Im z_1 <\Im z< \Im z_2\subset \C$ together with its natural orientation, but changes the sign
of the integrand,
$$
\sigma^*\left( dArg(z-z_1)\wedge dArg (z-z_2)\right)= dArg(z-z_2)\wedge dArg (z-z_1).
$$
Hence the claim.
\end{proof}

{\bf A4.2. Corollary}. {\em The weights, $\displaystyle c_\Ga= \int_{C_4}\wedge_{e\in \Ga}\frac{\pi_e^*(\om_K|_{outer\ circle})}{2\pi}$,   of the graphs
$$
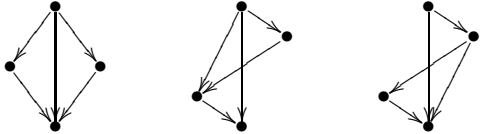
\Ba{c}\xy
(0,0)*{\bullet}="1",
(-6,8)*{\bullet}="2",
(6,8)*{\bullet}="3",
(0,16)*{\bullet}="4",
\ar @{<-} "2";"4" <0pt>
\ar @{<-} "3";"4" <0pt>
\ar @{->} "4";"1" <0pt>
\ar @{->} "2";"1" <0pt>
\ar @{->} "3";"1" <0pt>
\endxy\Ea \ \ \ \ \ \
 \Ba{c}\xy
(0,0)*{\bullet}="1",
(-6,4)*{\bullet}="2",
(6,12)*{\bullet}="3",
(0,16)*{\bullet}="4",
\ar @{<-} "2";"4" <0pt>
\ar @{->} "4";"3" <0pt>
\ar @{->} "4";"1" <0pt>
\ar @{->} "2";"1" <0pt>
\ar @{->} "3";"2" <0pt>
\endxy\Ea
\ \ \ \ \ \
 \Ba{c}\xy
(0,0)*{\bullet}="1",
(-6,4)*{\bullet}="2",
(6,12)*{\bullet}="3",
(0,16)*{\bullet}="4",
\ar @{<-} "1";"3" <0pt>
\ar @{->} "4";"3" <0pt>
\ar @{->} "4";"1" <0pt>
\ar @{->} "2";"1" <0pt>
\ar @{->} "3";"2" <0pt>
\endxy\Ea
$$
with respect to the propagator $\om_K|_{outer\ circle}$ vanish.
}

\mip

{\bf A4.3. Lemma}. {\em For any fixed points $z_1,z_2\in \C$ with $\Im z_1 < \Im z_2$ we have
$$
\int_{\Im z <\Img z_1} dArg(z-z_1)\wedge dArg (z-z_2)=-
\int_{\Im z >\Img z_2} dArg(z-z_1)\wedge dArg (z-z_2)= \frac{3\pi^2}{2}- \pi Arg(z_1-z_2).
$$
}
\begin{proof} Using the Fubini rule and the picture,
$$
{\Ba{c}\xy
<0mm,-2.0mm>*{_{z}},
<-9.5mm,9.3mm>*{^{z_1}},
<6.5mm,11.5mm>*{^{z_2}},
%
(0,0)*{\bullet}="1",
(-10,8)*{\bullet}="2",
(6,16)*{\bullet}="3",
(-25,16)*{}="u'",
(15,16)*{}="u''",
(-25,8)*{}="d'",
(15,8)*{}="d''",
(15,8)*{}="o",
\ar @{.} "3";"2" <0pt>
\ar @{->} "3";"1" <0pt>
\ar @{->} "2";"1" <0pt>
\ar @{--} "u'";"u''" <0pt>
\ar @{--} "d'";"d''" <0pt>
\endxy\Ea}
$$
we get, for example,
\Beqrn
\int_{\Im z <\Im z_1} dArg(z-z_1)\wedge dArg (z-z_2)&=& \int_{\pi}^{2\pi} dArg(z-z_1) \int_{Arg(z_1-z_2)}^{Arg(z-z_1)} dArg(z-z_2)\\
&=&  \int_{\pi}^{2\pi}\left(Arg(z-z_1) - Arg(z_1-z_2)\right) dArg(z-z_1) \\
&=&  \frac{3\pi^2}{2}- \pi Arg(z_1-z_2).
\Eeqrn
\end{proof}

{\bf A4.4. Corollary}. {\em The weights, $\displaystyle c_\Ga= \int_{C_4}\wedge_{e\in \Ga}\frac{\pi_e^*(\om_K|_{outer\ circle})}{2\pi}$,   of the graphs (see \S 4.2.4)
$$
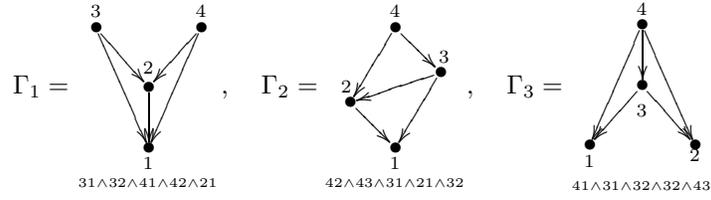
\Ga_1=  \underset{_{31\wedge 32\wedge 41\wedge 42\wedge21}}{\Ba{c}\xy
<0mm,-2.0mm>*{_{1}},
<0mm,10.0mm>*{^{2}},
<-7mm,17.5mm>*{^{3}},
<7mm,17.5mm>*{^{4}},
(0,0)*{\bullet}="1",
(-7,16)*{\bullet}="2",
(7,16)*{\bullet}="3",
(0,8)*{\bullet}="4",
\ar @{->} "2";"4" <0pt>
\ar @{->} "3";"4" <0pt>
\ar @{->} "4";"1" <0pt>
\ar @{->} "2";"1" <0pt>
\ar @{->} "3";"1" <0pt>
\endxy\Ea}, \ \ \
\Ga_2=  \underset{_{42\wedge 43\wedge 31\wedge 21\wedge 32}}{\Ba{c}\xy
<0mm,-2.0mm>*{_{1}},
<-6.5mm,7.5mm>*{^{2}},
<6.5mm,11.5mm>*{^{3}},
<0mm,17.5mm>*{^{4}},
(0,0)*{\bullet}="1",
(-6,6)*{\bullet}="2",
(6,10)*{\bullet}="3",
(0,16)*{\bullet}="4",
\ar @{->} "4";"3" <0pt>
\ar @{->} "4";"2" <0pt>
\ar @{->} "3";"2" <0pt>
\ar @{->} "2";"1" <0pt>
\ar @{->} "3";"1" <0pt>
\endxy\Ea}, \ \ \
\Ga_3=   \underset{_{41\wedge 31\wedge 32\wedge 32\wedge 43}}{ \Ba{c}\xy
<-7mm,-2.0mm>*{_{1}},
<7mm,-2.0mm>*{^{2}},
<0mm,4.0mm>*{^{3}},
<0mm,17.5mm>*{^{4}},
(0,16)*{\bullet}="1",
(-7,0)*{\bullet}="2",
(7,0)*{\bullet}="3",
(0,8)*{\bullet}="4",
\ar @{<-} "2";"4" <0pt>
\ar @{<-} "3";"4" <0pt>
\ar @{<-} "4";"1" <0pt>
\ar @{<-} "2";"1" <0pt>
\ar @{<-} "3";"1" <0pt>
\endxy\Ea}
$$
with respect to the propagator $\om_K|_{outer\ circle}$ are all equal to $\frac{1}{12}$.
}
\begin{proof}
By Lemma A4.3, all these weights are equal to the integral
$$
\frac{1}{\pi^5}\int_{\pi}^{2\pi} \left(\frac{3\pi^2}{2}- \pi x  \right)^2dx=\frac{1}{12}.
$$
\end{proof}

\bip
\begin{center}
{\sc Appendix 5: Wheeled prop of polyvector fields}
\end{center}

\bip
We refer to \cite{Me-Perm} for an elementary introduction into the language of
(wheeled) operads and props.
\sip

{\bf Definition}. {\em The wheeled prop of polyvector fields, $\cP oly$, is  a dg free
wheeled prop, $(\cF ree\left\langle E\right\rangle, \delta)$, generated by an $\bS$-bimodule
$E=\{E(m,n)\}_{m,n\geq 0}$,
$$
E(m,n)= \sgn_m\ot \id_n[m-2]=\mbox{span} \langle
\begin{xy}
 <0mm,0mm>*{\bullet};<0mm,0mm>*{}**@{},
 <0mm,0mm>*{};<-8mm,5mm>*{}**@{-},
 <0mm,0mm>*{};<-4.5mm,5mm>*{}**@{-},
 <0mm,0mm>*{};<-1mm,5mm>*{\ldots}**@{},
 <0mm,0mm>*{};<4.5mm,5mm>*{}**@{-},
 <0mm,0mm>*{};<8mm,5mm>*{}**@{-},
   <0mm,0mm>*{};<-8.5mm,5.5mm>*{^1}**@{},
   <0mm,0mm>*{};<-5mm,5.5mm>*{^2}**@{},
   <0mm,0mm>*{};<9.0mm,5.5mm>*{^m}**@{},
 <0mm,0mm>*{};<-8mm,-5mm>*{}**@{-},
 <0mm,0mm>*{};<-4.5mm,-5mm>*{}**@{-},
 <0mm,0mm>*{};<-1mm,-5mm>*{\ldots}**@{},
 <0mm,0mm>*{};<4.5mm,-5mm>*{}**@{-},
 <0mm,0mm>*{};<8mm,-5mm>*{}**@{-},
   <0mm,0mm>*{};<-8.5mm,-6.9mm>*{^1}**@{},
   <0mm,0mm>*{};<-5mm,-6.9mm>*{^2}**@{},
   <0mm,0mm>*{};<9.0mm,-6.9mm>*{^n}**@{},
 \end{xy}
\rangle
\vspace{-1mm}
$$
and equipped with the differential $\delta$ given on the generators by the formula
$$
\delta \begin{xy}
 <0mm,0mm>*{\bullet};<0mm,0mm>*{}**@{},
 <0mm,0mm>*{};<-8mm,5mm>*{}**@{-},
 <0mm,0mm>*{};<-4.5mm,5mm>*{}**@{-},
 <0mm,0mm>*{};<-1mm,5mm>*{\ldots}**@{},
 <0mm,0mm>*{};<4.5mm,5mm>*{}**@{-},
 <0mm,0mm>*{};<8mm,5mm>*{}**@{-},
   <0mm,0mm>*{};<-8.5mm,5.5mm>*{^1}**@{},
   <0mm,0mm>*{};<-5mm,5.5mm>*{^2}**@{},
   <0mm,0mm>*{};<9.0mm,5.5mm>*{^m}**@{},
 <0mm,0mm>*{};<-8mm,-5mm>*{}**@{-},
 <0mm,0mm>*{};<-4.5mm,-5mm>*{}**@{-},
 <0mm,0mm>*{};<-1mm,-5mm>*{\ldots}**@{},
 <0mm,0mm>*{};<4.5mm,-5mm>*{}**@{-},
 <0mm,0mm>*{};<8mm,-5mm>*{}**@{-},
   <0mm,0mm>*{};<-8.5mm,-6.9mm>*{^1}**@{},
   <0mm,0mm>*{};<-5mm,-6.9mm>*{^2}**@{},
   <0mm,0mm>*{};<9.0mm,-6.9mm>*{^n}**@{},
 \end{xy} =
 \sum_{[m]=I_1\sqcup I_2\atop {[n]=J_1\sqcup J_2\atop
 {|I_1|\geq 0, |I_2|\geq 1 \atop
 |J_1|\geq 1, |J_2|\geq 0}}
}\hspace{0mm} (-1)^{\sigma(I_1\sqcup I_2) + |I_1|(|I_2|+1)}
 \begin{xy}
 <0mm,0mm>*{\bullet};<0mm,0mm>*{}**@{},
 <0mm,0mm>*{};<-8mm,5mm>*{}**@{-},
 <0mm,0mm>*{};<-4.5mm,5mm>*{}**@{-},
 <0mm,0mm>*{};<0mm,5mm>*{\ldots}**@{},
 <0mm,0mm>*{};<4.5mm,5mm>*{}**@{-},
 <0mm,0mm>*{};<13mm,5mm>*{}**@{-},
     <0mm,0mm>*{};<-2mm,7mm>*{\overbrace{\ \ \ \ \ \ \ \ \ \ \ \ }}**@{},
     <0mm,0mm>*{};<-2mm,9mm>*{^{I_1}}**@{},
 <0mm,0mm>*{};<-8mm,-5mm>*{}**@{-},
 <0mm,0mm>*{};<-4.5mm,-5mm>*{}**@{-},
 <0mm,0mm>*{};<-1mm,-5mm>*{\ldots}**@{},
 <0mm,0mm>*{};<4.5mm,-5mm>*{}**@{-},
 <0mm,0mm>*{};<8mm,-5mm>*{}**@{-},
      <0mm,0mm>*{};<0mm,-7mm>*{\underbrace{\ \ \ \ \ \ \ \ \ \ \ \ \ \ \
      }}**@{},
      <0mm,0mm>*{};<0mm,-10.6mm>*{_{J_1}}**@{},
 <13mm,5mm>*{};<13mm,5mm>*{\bullet}**@{},
 <13mm,5mm>*{};<5mm,10mm>*{}**@{-},
 <13mm,5mm>*{};<8.5mm,10mm>*{}**@{-},
 <13mm,5mm>*{};<13mm,10mm>*{\ldots}**@{},
 <13mm,5mm>*{};<16.5mm,10mm>*{}**@{-},
 <13mm,5mm>*{};<20mm,10mm>*{}**@{-},
      <13mm,5mm>*{};<13mm,12mm>*{\overbrace{\ \ \ \ \ \ \ \ \ \ \ \ \ \ }}**@{},
      <13mm,5mm>*{};<13mm,14mm>*{^{I_2}}**@{},
 <13mm,5mm>*{};<8mm,0mm>*{}**@{-},
 <13mm,5mm>*{};<12mm,0mm>*{\ldots}**@{},
 <13mm,5mm>*{};<16.5mm,0mm>*{}**@{-},
 <13mm,5mm>*{};<20mm,0mm>*{}**@{-},
     <13mm,5mm>*{};<14.3mm,-2mm>*{\underbrace{\ \ \ \ \ \ \ \ \ \ \ }}**@{},
     <13mm,5mm>*{};<14.3mm,-4.5mm>*{_{J_2}}**@{},
 \end{xy}
$$
where $\sigma(I_1\sqcup I_2)$ is the sign of the permutation $[n]\rar I_1\sqcup I_2$.
}

\mip

Here $\sgn_n$ (resp.\ $\id_n$) is the 1-dimensional sign (resp.\ trivial) representation
of $\bS_n$. Representations, $f: \cP oly\rar \cE nd_V$, of this prop in a
finite-dimensional  $\Z$-graded vector
space $V$ are in one-to-one correspondence with formal $\Z$-graded Poisson structures on $V$
\cite{Me-lec}. The group, $Aut(\cP oly)$, of automorphisms of this prop consists of
all automorphisms, $F: \cF ree\left\langle E\right\rangle \rar \cF ree\left\langle E\right\rangle$,
of the free prop which respect the differential, $F\circ \delta=\delta\circ F$. Every such an
automorphism
is uniquely determined by its values on the generators,
$$
F\left(
\begin{xy}
 <0mm,0mm>*{\bullet};<0mm,0mm>*{}**@{},
 <0mm,0mm>*{};<-8mm,5mm>*{}**@{-},
 <0mm,0mm>*{};<-4.5mm,5mm>*{}**@{-},
 <0mm,0mm>*{};<-1mm,5mm>*{\ldots}**@{},
 <0mm,0mm>*{};<4.5mm,5mm>*{}**@{-},
 <0mm,0mm>*{};<8mm,5mm>*{}**@{-},
   <0mm,0mm>*{};<-8.5mm,5.5mm>*{^1}**@{},
   <0mm,0mm>*{};<-5mm,5.5mm>*{^2}**@{},
   <0mm,0mm>*{};<9.0mm,5.5mm>*{^m}**@{},
 <0mm,0mm>*{};<-8mm,-5mm>*{}**@{-},
 <0mm,0mm>*{};<-4.5mm,-5mm>*{}**@{-},
 <0mm,0mm>*{};<-1mm,-5mm>*{\ldots}**@{},
 <0mm,0mm>*{};<4.5mm,-5mm>*{}**@{-},
 <0mm,0mm>*{};<8mm,-5mm>*{}**@{-},
   <0mm,0mm>*{};<-8.5mm,-6.9mm>*{^1}**@{},
   <0mm,0mm>*{};<-5mm,-6.9mm>*{^2}**@{},
   <0mm,0mm>*{};<9.0mm,-6.9mm>*{^n}**@{},
 \end{xy}
\right)
=
\begin{xy}
 <0mm,0mm>*{\bullet};<0mm,0mm>*{}**@{},
 <0mm,0mm>*{};<-8mm,5mm>*{}**@{-},
 <0mm,0mm>*{};<-4.5mm,5mm>*{}**@{-},
 <0mm,0mm>*{};<-1mm,5mm>*{\ldots}**@{},
 <0mm,0mm>*{};<4.5mm,5mm>*{}**@{-},
 <0mm,0mm>*{};<8mm,5mm>*{}**@{-},
   <0mm,0mm>*{};<-8.5mm,5.5mm>*{^1}**@{},
   <0mm,0mm>*{};<-5mm,5.5mm>*{^2}**@{},
   <0mm,0mm>*{};<9.0mm,5.5mm>*{^m}**@{},
 <0mm,0mm>*{};<-8mm,-5mm>*{}**@{-},
 <0mm,0mm>*{};<-4.5mm,-5mm>*{}**@{-},
 <0mm,0mm>*{};<-1mm,-5mm>*{\ldots}**@{},
 <0mm,0mm>*{};<4.5mm,-5mm>*{}**@{-},
 <0mm,0mm>*{};<8mm,-5mm>*{}**@{-},
   <0mm,0mm>*{};<-8.5mm,-6.9mm>*{^1}**@{},
   <0mm,0mm>*{};<-5mm,-6.9mm>*{^2}**@{},
   <0mm,0mm>*{};<9.0mm,-6.9mm>*{^n}**@{},
 \end{xy}\  + \
  \sum_{k\geq 2}\sum_{\Ga\in G_{k,2k-2}^\circlearrowright(m,n)} c_\Gamma \Ga, \ \ c_\Ga\in\C,
$$
which, for purely degree reasons, must be a sum over a family $G_{k,2k-2}^\circlearrowright(m,n)$
of
 graphs $\Ga$ which are built
from
 the generating corollas
$$
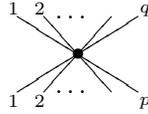
\begin{xy}
 <0mm,0mm>*{\bullet};<0mm,0mm>*{}**@{},
 <0mm,0mm>*{};<-8mm,5mm>*{}**@{-},
 <0mm,0mm>*{};<-4.5mm,5mm>*{}**@{-},
 <0mm,0mm>*{};<-1mm,5mm>*{\ldots}**@{},
 <0mm,0mm>*{};<4.5mm,5mm>*{}**@{-},
 <0mm,0mm>*{};<8mm,5mm>*{}**@{-},
   <0mm,0mm>*{};<-8.5mm,5.5mm>*{^1}**@{},
   <0mm,0mm>*{};<-5mm,5.5mm>*{^2}**@{},
   <0mm,0mm>*{};<9.0mm,5.5mm>*{^q}**@{},
 <0mm,0mm>*{};<-8mm,-5mm>*{}**@{-},
 <0mm,0mm>*{};<-4.5mm,-5mm>*{}**@{-},
 <0mm,0mm>*{};<-1mm,-5mm>*{\ldots}**@{},
 <0mm,0mm>*{};<4.5mm,-5mm>*{}**@{-},
 <0mm,0mm>*{};<8mm,-5mm>*{}**@{-},
   <0mm,0mm>*{};<-8.5mm,-6.9mm>*{^1}**@{},
   <0mm,0mm>*{};<-5mm,-6.9mm>*{^2}**@{},
   <0mm,0mm>*{};<9.0mm,-6.9mm>*{^p}**@{},
 \end{xy}
$$
by taking their disjoint unions and then gluing  some output legs with with the same number of input legs
and which satisfy three conditions: $\Ga$ has $k$ vertices, $2k-2$ edges,
 $n$ input legs and  $m$ output
legs (cf.\ \cite{Me-Perm}). The main result of our paper can be restated
as follows: any de Rham field theory on $\overline{C}\sqcup \widehat{C}$
defines an exotic automorphism of $(\cP oly, \delta)$ with weights $c_\Gamma$ given
by (\ref{Ch4: weight C_Ga}).

\bip

{\em Acknowledgements}. {\small It is a pleasure to thank Vasilij Dolgushev, Giovanni Felder,
Johan Gran\aa ker, Pascal Lambrechts,  Sergei Shadrin and especially
Thomas Willwacher  for very useful discussions and insightful
comments. I am also grateful to Anton Alekseev and Charles Torossian for showing me a preliminary version
of their work \cite{AT}. Finally, I thank an anonymous referee for useful criticism.}

\newpage

\def\cprime{$'$}

\end{document}